\documentclass[11pt,leqno]{amsart}

\usepackage[utf8x]{inputenc}
\usepackage[T1]{fontenc}
\usepackage[english]{babel}
\usepackage{textcmds, cite}

\usepackage{amsmath,amssymb,amsthm}
\usepackage{mathabx}
\usepackage{yfonts}
\usepackage{mathrsfs}
\usepackage{faktor}
\usepackage{amsrefs}
\usepackage[dvipsnames]{xcolor}
\usepackage[colorlinks=true, linkcolor=MidnightBlue, citecolor=OliveGreen]{hyperref}
\usepackage{multicol}
\usepackage{multirow}
\usepackage{enumerate}
\usepackage{float}

\usepackage{graphicx}
\usepackage{tikz}
\usepackage{tikz-cd}
\usetikzlibrary{arrows, automata, positioning, calc}

\theoremstyle{definition}
\newtheorem{defn}{{\bf Definition}}[section]
\newtheorem{notation}[defn]{{\bf Notation}}
\newtheorem{eg}[defn]{{\bf Example}}

\newtheorem{qst}[defn]{{\bf Question}}
\theoremstyle{theorem}
\newtheorem{lemma}[defn]{{\bf Lemma}}
\newtheorem{theorem}[defn]{{\bf Theorem}}
\newtheorem{coro}[defn]{{\bf Corollary}}
\newtheorem{prop}[defn]{{\bf Proposition}}

\newcommand*{\defeq}{\mathrel{\vcenter{\baselineskip0.5ex \lineskiplimit0pt\hbox{\scriptsize.}\hbox{\scriptsize.}}}=}
\newcommand*{\eqdef}{=\mathrel{\vcenter{\baselineskip0.5ex \lineskiplimit0pt \hbox{\scriptsize.}\hbox{\scriptsize.}}}}

\DeclareMathOperator{\Q}{\mathbb{Q}}
\DeclareMathOperator{\Z}{\mathbb{Z}}
\DeclareMathOperator{\N}{\mathbb{N}_0}
\DeclareMathOperator{\Nplus}{\mathbb{N}_+}
\DeclareMathOperator{\id}{id}

\DeclareMathOperator{\Sym}{\mathrm{Sym}}
\DeclareMathOperator{\Aut}{Aut}

\DeclareMathOperator{\Autf}{Aut_{fin}}
\DeclareMathOperator{\st}{st}
\DeclareMathOperator{\St}{St}

\DeclareMathOperator{\Rst}{Rist}
\DeclareMathOperator{\ab}{ab}
\DeclareMathOperator{\syl}{syl}
\DeclareMathOperator{\dimH}{dim_H}
\newcommand{\LD}{\langle}
\newcommand{\RD}{\rangle}

\newcommand{\intseg}[2]{{[#1, #2]}} 

\DeclareMathOperator{\bp}{Bas}
\newcommand{\cref}[3][]{\hyperref[#3]{#2~\ref*{#3}#1}}

\title{On the Basilica Operation}

\author[J. M. Petschick]{Jan Moritz Petschick} \address{Jan Moritz Petschick: Mathematisches
  Institut, Heinrich-Heine-Universit\"at, 40225 D\"usseldorf, Germany}
\email{jan.petschick@hhu.de}

\author[K. Rajeev]{Karthika Rajeev}  
\address{Karthika Rajeev: Mathematisches Institut, Heinrich-Heine-Universit\"at, 40225 D\"usseldorf, Germany}
\email{rajeev@uni-duesseldorf.de}

\thanks{The research was funded by the Deutsche Forschungsgemeinschaft (DFG, German Research Foundation) — 380258175. The research was also partially conducted in the framework of the DFG-funded research training group “GRK 2240: Algebro-Geometric Methods in Algebra, Arithmetic and Topology”.}



\date{\today}

\begin{document}

\begin{abstract}
 	Inspired by the Basilica~group $\mathcal B$, we describe a general construction which allows us to associate to any group of automorphisms $G \leq \Aut(T)$ of a rooted tree $T$ a family of Basilica groups $\bp_s(G), s \in \Nplus$. For the dyadic odometer $\mathcal O_2$, one has $\mathcal B = \bp_2(\mathcal O_2)$. We study which properties of groups acting on rooted trees are preserved under this operation. Introducing some techniques for handling $\bp_s(G)$, in case $G$ fulfills some branching conditions, we are able to calculate the Hausdorff dimension of the Basilica groups associated to certain $\mathsf{GGS}$-groups and of generalisations of the odometer, $\mathcal O_m^d$. Furthermore, we study the structure of groups of type $\bp_s(\mathcal O_m^d)$ and prove an analogue of the congruence subgroup property in the case $m = p$, a prime.
\end{abstract}

\maketitle

\section{Introduction}\label{sec:intro}

\noindent Groups acting on rooted trees play an important role in various areas of group theory, for example in the study of groups of intermediate growth, just infinite groups and groups related to the Burnside problem. Over the years, many groups of automorphisms of rooted trees have been defined and studied. Often they can be regarded as generalisations of early constructions to wider families of groups with similar properties.

In this paper, we consider an operation on the subgroups of the automorphism group $\Aut(T)$ of a rooted tree $T$ with degree $m \geq 2$. It is inspired by the \emph{Basilica~group} $\mathcal B$, a group acting on the binary rooted tree, which was introduced by Grigorchuk and Żuk in \cite{GZ01} and \cite{GZ02}. The Basilica~group $\mathcal B$ is a particularly interesting example in its own right: it is a self-similar torsion-free weakly branch group, just-(non-soluble) and of exponential word growth. It was the first group known to be not sub-exponentially amenable \cite{GZ02}, but amenable \cites{BV05, BKN10}. Furthermore, it is the iterated monodromy group of $z^2-1$ \cites{GZ01, Pin13}, and it has the $2$-congruence subgroup property \cite{GUA19}.

The Basilica group $\mathcal B$ is usually defined as the group generated by two automorphisms
\[
	a = (b, \id) \text{ and } b = (0\; 1) (a, \id),
\]
acting on the binary rooted tree (in \cite{GZ02} the elements are defined with $\id$ on the left, which is merely notational). We point out the similarities between these two generators and the single automorphism generating the \emph{dyadic odometer}. The latter provides an embedding of the infinite cyclic group into the automorphism group $\Aut(T)$ of the binary rooted tree $T$, given by
\[
	c = (0\; 1)(c, \id).
\]
We can regard $b$ as a delayed version of $c$, that takes an intermediate step acting as $a$, before returning to itself. Considering the automata defining the generators of both groups (cf.\ \cref{Figure}{fig:automata bas}), the relationship is even more apparent. We obtain the automaton defining $b$ from the automaton defining $c$ by replacing every edge that does not point to the state of the trivial element with an edge pointing to a new state, which in turn points to the old state upon reading $0$ and to the state of the trivial element upon reading any other letter. See \cref{Figure}{fig:Replacement rule for edges in an automaton} for an illustration of this replacement rule.

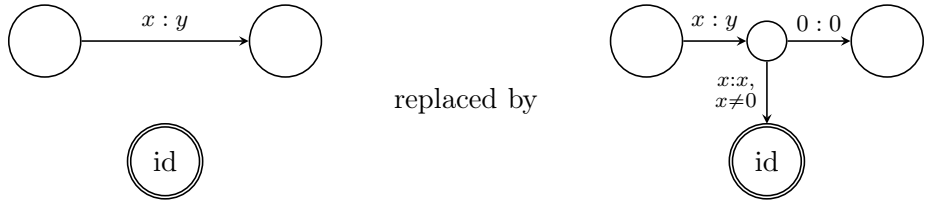
\begin{figure}[H]
\centering
\begin{tikzpicture}[>=stealth,semithick,->,scale=0.8]
	\node[state] (b) at (0,0) {};
	\node[state] (a) at (4,0) {};
	\node[state, accepting] (id) at (2,-2) {$\id$};
			
	\path
		(b) edge node[font=\footnotesize, above]{$x:y$} (a);
	
	\node at (7, -1) {replaced by};
	
	\node[state, accepting] (id) at (12,-2) {$\id$};
	\node[state] (b) at (10,0) {};
	\node[state] (a) at (14,0) {};
	\node[state, inner sep=0pt,minimum size=15pt,font=\footnotesize] (A) at (12,0) {};
	
	\path
		(b) edge node[font=\footnotesize, above]{$x:y$} (A)
		(A) edge node[font=\footnotesize, above]{$0:0$} (a)
		(A) edge node[left]{$\substack{x:x, \\x \neq 0}$} (id);
\end{tikzpicture}
\caption{Replacement rule for edges.}\label{fig:Replacement rule for edges in an automaton}
\end{figure}

The same can be done for any automorphism of $T$ and any number $s$ of intermediate states. For any group of automorphisms $G$, this operation yields a new group of tree automorphisms defined by the automaton with $s$ intermediate steps, which we call $\bp_s(G)$, the $s$-th Basilica group of $G$. A precise, algebraic definition that does not refer to automata will be given in \cref{Definition}{def:bp}. \cref{Figure}{fig:automata bas} depicts for example the automaton defining $\bp_8(\mathcal O_2)$, while \cref{Figure}{fig:gs3} depicts the automaton defining the generators of the Gupta--Sidki $3$-group $\ddot\Gamma$ and the corresponding automaton obtained by the operation $\bp_2$.

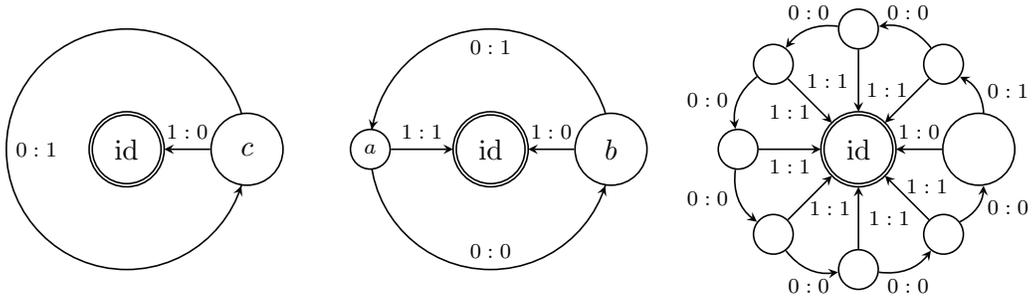
\begin{figure}[H]
	\minipage{0.33\textwidth}
	\centering
	\begin{tikzpicture}[>=stealth,semithick,->,scale=0.8]
		\draw (0,0) circle (2);
	
		\node[state, accepting] (I) at (0,0) {$\id$};
		\node[state, fill=white] (A) at (2,0) {$c$}; 
	
		\node[font=\scriptsize] at (-1.5, 0) {$0:1$};
		\path
			(A) edge[above] node[font=\scriptsize]{$1:0$} (I)
			(1.87,-0.7) edge (1.905,-0.61);	
	\end{tikzpicture}
	\endminipage\hfill
	\minipage{0.33\textwidth}
	\centering
	\begin{tikzpicture}[>=stealth,semithick,->,scale=0.8]
		\draw (0,0) circle (2);
		
		\node[state, accepting] (I) at (0,0) {$\id$};
		\node[state, fill=white] (b) at (2,0) {$b$}; 
		\node[state, fill=white, inner sep=0pt,minimum size=15pt,font=\scriptsize] (a) at (-2,0) {$a$}; 
		
		\node[font=\scriptsize] at (0, 1.7) {$0:1$};
		\node[font=\scriptsize] at (0, -1.7) {$0:0$};
		\path
			(b) edge[above] node[font=\scriptsize]{$1:0$} (I)
			(a) edge[above] node[font=\scriptsize]{$1:1$} (I)
			(1.87,-0.7) edge (1.905,-0.61)
			(-1.96,0.42) edge (-1.98,0.33);
	\end{tikzpicture}
	\endminipage\hfill
	\minipage{0.34\textwidth}
	\centering
	\begin{tikzpicture}[>=stealth,semithick,->,scale=0.8]
		\node[state, accepting] (I) at (0,0) {$\id$};
		
		\node[state, fill=white] (b) at (2,0) {};
		
		\node[state, fill=white, inner sep=0pt,minimum size=15pt,font=\scriptsize] (1) at ({sqrt(2)},{sqrt(2)}) {};
		\node[state, fill=white, inner sep=0pt,minimum size=15pt,font=\scriptsize] (2) at (0,2) {};
		\node[state, fill=white, inner sep=0pt,minimum size=15pt,font=\scriptsize] (3) at ({-sqrt(2)},{sqrt(2)}) {};
		\node[state, fill=white, inner sep=0pt,minimum size=15pt,font=\scriptsize] (4) at (-2,0) {};
		\node[state, fill=white, inner sep=0pt,minimum size=15pt,font=\scriptsize] (5) at ({-sqrt(2)},{-sqrt(2)}) {};
		\node[state, fill=white, inner sep=0pt,minimum size=15pt,font=\scriptsize] (6) at (0,-2) {};
		\node[state, fill=white, inner sep=0pt,minimum size=15pt,font=\scriptsize] (7) at ({sqrt(2)},{-sqrt(2)}) {};
		
		\path
			(b) edge node[font=\scriptsize, above]{$1:0$} (I)
			(1) edge node[font=\scriptsize, near start, left]{$1:1$} (I)
			(2) edge node[font=\scriptsize, left]{$1:1$} (I)
			(3) edge node[font=\scriptsize, near end, left]{$1:1$} (I)
			(4) edge node[font=\scriptsize, below]{$1:1$} (I)
			(5) edge node[font=\scriptsize, near start, right]{$1:1$} (I)
			(6) edge node[font=\scriptsize, right]{$1:1$} (I)
			(7) edge node[font=\scriptsize, near end, right]{$1:1$} (I);
		\path[bend right]
			(b) edge node[font=\scriptsize, right]{$0:1$} (1)
			(1) edge node[font=\scriptsize, above]{$0:0$} (2)
			(2) edge node[font=\scriptsize, above]{$0:0$} (3) 
			(3) edge node[font=\scriptsize, left]{$0:0$} (4) 
			(4) edge node[font=\scriptsize, left]{$0:0$} (5) 
			(5) edge node[font=\scriptsize, below]{$0:0$} (6) 
			(6) edge node[font=\scriptsize, below]{$0:0$} (7) 
			(7) edge node[font=\scriptsize, right]{$0:0$} (b);
	\end{tikzpicture}
	\endminipage
	\caption{Automata for the dyadic odometer $\mathcal O_2$, the Basilica group $\mathcal{B}=\bp_2(\mathcal O_2)$, and $\bp_8(\mathcal O_2)$.
	}\label{fig:automata bas}
\end{figure}

\begin{figure}[H]
	\minipage{0.5\textwidth}
	\centering
	\begin{tikzpicture}[>=stealth,semithick,->,scale=0.8]
		\node[state, accepting] (id) at (0,0) {$\id$};
		\node[state] (b) at (4,0) {};
		\node[state] (a) at (0,-2) {};
		\node[state] (aa) at (0,2) {};
				
		\path
			(b) edge (a)
			(b) edge (aa)
			(a)	edge node[font=\scriptsize, left]{$x: \sigma^2(x)$} (id)
			(aa)edge node[font=\scriptsize, left]{$x: \sigma(x)$} (id);
		
		\node[font=\scriptsize] at (2, -1.4) {$2:2$};
		\node[font=\scriptsize] at (2, 1.4) {$1:1$};
		
		\path[->,min distance=2cm] (b)edge[in=190,out=170,left] node[font=\scriptsize] {$0:0$}(b);
	\end{tikzpicture}
	\endminipage\hfill
	\minipage{0.5\textwidth}
	\centering
	\begin{tikzpicture}[>=stealth,semithick,->,scale=0.8]
		\node[state, accepting] (id) at (0,0) {$\id$};
		\node[state] (b) at (4,0) {};
		\node[state] (a) at (0,-2) {};
		\node[state] (aa) at (0,2) {};
		
		\node[state, inner sep=0pt,minimum size=15pt,font=\scriptsize] (B) at (2,0) {};
		\node[state, inner sep=0pt,minimum size=15pt,font=\scriptsize] (A) at (2,-1) {};
		\node[state, inner sep=0pt,minimum size=15pt,font=\scriptsize] (AA) at (2,1) {};
		
		\path
			(b) edge[bend left] node[font=\scriptsize, below]{$2:2$} (A)
			(b) edge[bend right] node[font=\scriptsize, above]{$1:1$} (AA)
			(b) edge[bend left] node[font=\scriptsize, below, near end]{$0:0$} (B)
			(B) edge[bend left] node[font=\scriptsize, above, near start]{$0:0$} (b)
			(B) edge node[font=\scriptsize, near start]{$\substack{x \neq 0,\\ \vspace{-0.3em}\\ x:x}$} (id)
			(A) edge (a)
			(A) edge node[font=\scriptsize, left, near start]{$\substack{\\\vspace{-0.3em}\\x \neq 0,\\ x:x}$} (id)
			(a)	edge node[font=\scriptsize, left]{$x: \sigma^2(x)$} (id)
			(AA)edge (aa)
			(AA)edge node[font=\scriptsize, above, near end]{$\substack{x \neq 0,\\ x:x}$} (id)
			(aa)edge node[font=\scriptsize, left]{$x: \sigma(x)$} (id);
		\node[font=\scriptsize] at (1.2, -1.8) {$0:0$};
		\node[font=\scriptsize] at (1.2, 1.8) {$0:0$};
	\end{tikzpicture}
	\endminipage
	\caption{Automata for the Gupta--Sidki $3$-group $\ddot\Gamma$ and $\bp_2(\ddot\Gamma)$, where $\sigma$ is a cyclic permutation.
	}\label{fig:gs3} 
\end{figure}
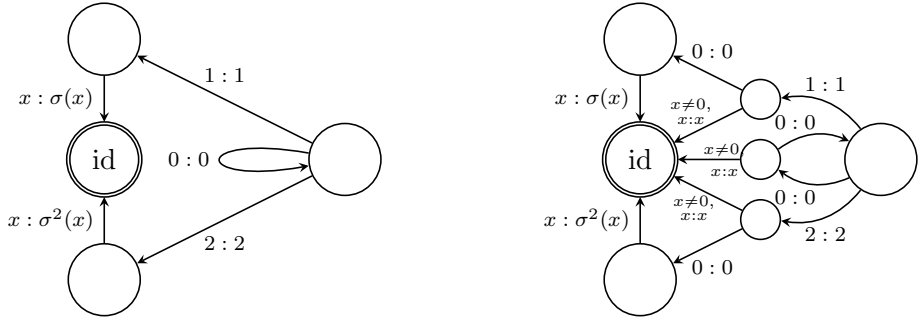

We prove that many of the desirable properties of the original Basilica group $\mathcal{B}$ are a consequence of the fact that the binary odometer $\mathcal O_2$ has those properties and that the properties are preserved under the Basilica operation. We summarise results of this kind for the general Basilica operation in the following theorem.

\begin{theorem}\label{thm:transfered_prop}
	Let $G$ be a group of automorphisms of a regular rooted tree. Let $P$ be a property from the list below. Then, if $G$ has $P$, the $s$-th Basilica group $\bp_s (G)$ of $G$ has $P$ for all $s \in \Nplus$.
	\begin{multicols}{2}
	\begin{enumerate}
		\item spherically transitive
		\item self-similar
		\item (strongly) fractal
		\item contracting
		\item weakly branch
		\item generated by finite-state bounded automorphisms
	\end{enumerate}
	\end{multicols}
\end{theorem}

As a consequence we derive conditions for $\bp_s(G)$ to have solvable word problem and to be amenable. Furthermore, we provide a condition for $\bp_s(G)$ to be a weakly regular branch group given that $G$ satisfies a group law. This enables us to construct a weakly regular branch group branching over a prescribed verbal subgroup.

The class of spinal groups, defined in \cite{BS01}, is another important class of groups acting on $T$; it contains the Grigorchuk group and all $\mathsf{GGS}$-groups, see \cref{Definition}{def:spinal}. It is not true that the Basilica operation preserves being spinal, however groups obtained from spinal groups act as spinal groups on another tree $\delta_s T$, obtained by deleting layers from $T$.

\begin{theorem}\label{thm:spinal}
	Let $G$ be a spinal group (resp.\ a $\mathsf{GGS}$-group) acting on $T$. Then $\bp_s G$ is a spinal (resp.\ a $\mathsf{GGS}$-group) acting on $\delta_s T$ for all $s \in \Nplus$.
\end{theorem}

In contrast to \cref{Theorem}{thm:transfered_prop}, the exponential word growth of the original Basilica group $\mathcal{B}$ is not a general feature of groups obtained by the Basilica operation. In fact, the situation appears to be chaotic, for which we provide some examples, see \cref{Proposition}{prop:poly_growth_to_exp_and_poly_growth} and \cref{Proposition}{prop:inter_growth_to_exp}.

Next we turn our attention to a class of groups $G$ whose Basilica groups $\bp_s(G)$ more closely resemble the original Basilica group. For this, we introduce the concept of the group $G$ being \textit{$s$-split} (see \cref{Definition}{def:s-split}). An $s$-split group decomposes by definition as a semi-direct product, algebraically modelling the property that the image of a delayed automorphism can be detected by observing the layers on which it has trivial labels. We prove that all abelian groups acting locally regular are $s$-split for all $s \in \Nplus$, and that conversely, all $s$-split groups acting spherically transitive are abelian. Furthermore we obtain the following.

\begin{theorem}\label{thm:split_transfer}
	Let $s > 1$ and let $G$ be an $s$-split self-similar group of automorphisms of a regular rooted tree acting spherically transitively. If $G$ is torsion-free, then $\bp_s(G)$ is torsion-free. Furthermore
	\(
		\bp_s(G)^{\mathrm{ab}} \cong G^s.
	\)
\end{theorem}

The \emph{$(s-1)$-th splitting kernel} $K_{s-1}$ is a normal subgroup of $G$ measuring the failure of $G$ to be $s$-split. A rigorous definition is found in \cref{Definition}{def:s-split}. If $G$ is weakly regular branch over $K_{s-1}$ (allowing $K_{s-1}$ to be trivial, hence including $s$-split groups), we obtain a strong structural description of the layer stabilisers of $\bp_s(G)$. The maps $\beta_i$ are the algebraic analogues of the various added steps delaying an automorphism, defined in \cref{Definition}{def:beta}.
\begin{theorem}\label{thm:stab}
	Let $G$ be a self-similar and very strongly fractal group of automorphisms of a regular rooted tree.
	Assume that $G$ is weakly regular branch over $K_{s-1}$. Let $n\in \N$. Write $n=sq+r$ with $q \geq 0$ and $0 \leq r \leq s-1$. Then for all $s > 1$
	\[
		\St_{\bp_s(G)}(n) = \LD \beta_i(\St_G( q + 1 )), \beta_{j}(\St_G( q )) \mid 0 \leq i < r \leq j < s\RD^{\bp_s(G)}.
	\]
\end{theorem}

This description allows us to provide an exact relationship between the Hausdorff dimension of a group $G$ fulfilling the conditions of \cref{Theorem}{thm:stab} and its Basilica groups $\bp_s(G)$. The precise description makes use of the \textit{series of obstructions} of $G$, a tailor-made technical construction, see \cref{Subsection}{sub:hausdorff_dimension} for details. Observing this series, we prove that the Hausdorff dimension of $\bp_s(G)$ is bounded below by the Hausdorff dimension of $G$ for all $s>1$.

\begin{coro}\label{cor:hausdorff_bounded}
	Let $G \leq \Aut T$ be very strongly fractal, self-similar, weakly regular branch over $K_{s-1}$, with $\dimH G < 1$. Then for all $s > 1$
	\[
		\dimH G < \dimH \bp_s(G).
	\]
\end{coro}
Here we define the Hausdorff dimension of $G \leq \Gamma$ as the Hausdorff dimension of its closure in $\Gamma$, where $\Gamma$ is the subgroup of all automorphisms acting locally by a power of a fixed $m$-cycle. This subgroup is isomorphic to
  \[
     \Gamma \cong \varprojlim \limits_{n \in \Nplus}{ \mathrm C_m \wr \overset{n}\cdots \wr \mathrm C_m }.
  \]
If $m =p$, a prime, then $\Gamma$ is a Sylow pro-$p$ subgroup of $\Aut(T)$. The notion of Hausdorff dimension in the profinite setting as above was initially studied by Abercrombie \cite{Abe94} and subsequently by Barnea and Shalev \cite{BS97}. It is analogous to the Hausdorff dimension defined as usual over $\mathbb{R}$.

In the second half of this paper we study the class of \emph{generalised Basilica groups} $\bp_s(\mathcal O_m^d)$, for $d,\,m,\,s \in \Nplus$ with $m,\,s \geq 2$, defined by applying $\bp_s$ to the free abelian group of rank $d$ with a self-similar action derived from the $m$-adic odometer.  We remark that the above generalisation of the original Basilica group $\mathcal{B}$ is different from the one given in \cite{BN08}, but it includes the class of $p$-Basilica groups, where $p$ is a prime, studied recently in \cite{DNT21}. For every odd prime $p$, we obtain the $p$-Basilica group by setting $d=1, m=p$ and $s=2$ in $\bp_s(\mathcal O_m^d)$. Our construction also includes special cases, $d=1$ and $m = s = p$, studied by Hanna Sasse in her master's thesis supervised by Benjamin Klopsch. We record the properties of the generalised Basilica groups in the following theorem. 

\begin{theorem}\label{thm:gen_bas_properties}
	Let $d,\,m,\,s \in \Nplus$ with $m,\,s\geq 2$. Let $B = \bp_s(\mathcal{O}_m^{d})$ be the generalised Basilica group. The following assertions hold:
	\begin{enumerate}[(i)]
		\item $B$ acts spherically transitively on the corresponding $m$-regular rooted tree,
		\item $B$ is self-similar and strongly fractal,
		\item $B$ is contracting, and has solvable word problem,
		\item The group $\mathcal{O}_m^{d}$ is $s$-split, and $B^\mathrm{ab} \cong \Z^{ds}$,
		\item $B$ is torsion-free,
		\item $B$ is weakly regular branch over its commutator subgroup,
		\item $B$ has exponential word growth.
	\end{enumerate}
\end{theorem}
\cref[(i)]{Theorem}{thm:gen_bas_properties} to \cref[(vi)]{Theorem}{thm:gen_bas_properties} are obtained by direct application of \cref{Theorem}{thm:transfered_prop} and \cref{Theorem}{thm:split_transfer}.  The proof of \cref[(vii)]{Theorem}{thm:gen_bas_properties} is analogous to that of the original Basilica group $\mathcal{B}$ and can easily be generalised from \cite[Proposition 4]{GZ02}. Nevertheless, one can prove \cref{Theorem}{thm:gen_bas_properties} directly by considering the action of the group on the corresponding rooted tree, see \cite{Sas18}.
 
We explicitly compute the Hausdorff dimension of $\bp_s(\mathcal O_m^d)$, which turns out to be independent of the rank $d$ of the free abelian group $\mathcal O_m^d$:

\begin{theorem}\label{thm:gen_bas_haus}
	For all $d,\,m,\,s \in \Nplus$ with $m,\,s \geq 2$
	\[
		\dimH(\bp_s(\mathcal O_m^d)) = \frac{m(m^{s-1}-1)}{m^s-1}.
	\]
\end{theorem}

The above equality agrees with the formula of the Hausdorff dimension of $p$-Basilica groups given by \cite{DNT21}, and also with the Hausdorff dimension of the original Basilica group $\mathcal{B}$ given in \cite{Bar06}.

\begin{theorem}\label{thm:gen_bas_Lpres}
    Let $d,\,m,\,s \in \Nplus$ with $m,\,s \geq 2$. The generalised Basilica group $\bp_s(\mathcal O_m^d)$ admits an $L$-presentation 
    \[
     L = \LD Y \mid Q \mid \Phi \mid R \RD
    \]
    where the data $Y,\, Q, \, R$ and $\Phi$ are specified in \cref{Section}{sec:gen_bas_Lpres}. 
\end{theorem}

The concrete $L$-presentation requires unwieldy notation, whence it is not given here. It is analogous to the $L$-presentation of the original Basilica group $\mathcal{B}$ \cite{GZ02}. The name $L$-presentation stands as a tribute to Igor Lysionok who obtained such a presentation for the Grigorchuk group in \cite{Lys85}. It is now known that, every finitely generated, contracting, regular branch group admits a finite $L$-presentation but it is not finitely presentable (cf. \cite{Bar03}). Unfortunately, this result is not applicable to generalised Basilica groups as they are merely weakly branch. Also, the $L$-presentation of the generalised Basilica group is not finite as the set of relations is infinite. Nonetheless, akin to \cite[Proposition 11]{GZ02}, we can introduce a set of endomorphisms of the free group on the set of generators of the generalised Basilica group and obtain a finite $L$-presentation, see \cref{Definition}{def:L-ptn}, as defined in \cite{Bar03}.

Using the concrete $L$-presentation of a generalised Basilica group, we obtain the following structural result.

\begin{theorem} \label{thm:gamma}
    Let $d,\,m,\,s \in \Nplus$ with $m,\,s\geq 2$ and let $B$ be the generalised Basilica group $\bp_s(\mathcal{O}_m^{d})$.  We have:
 \begin{enumerate}
    \item[(i)] For $s = 2$, the quotient group $\gamma_2(B)/\gamma_3(B) \cong \Z^{d^2}$,
    \item[(ii)] For $s > 2$, the quotient group $\gamma_2(B)/\gamma_3(B) \cong \mathrm C_m^{ds-2} \times \mathrm C_{m^2}$. 
 \end{enumerate}
\end{theorem}

This implies that the quotients $\gamma_i(B)/\gamma_{i+1}(B)$ of  consecutive terms of the lower central series of a generalised Basilica group for $s>2$ are finite for all $i \geq 2$, whereas a similar behaviour happens for the original Basilica group $\mathcal{B}$ from $i \geq 3$, see \cite{BEH08} for details.

For a group $G$ of automorphisms of an $m$-regular rooted tree, we say that $G$ has the congruence subgroup property (CSP) if every subgroup of finite index in $G$ contains some layer stabiliser in $G$. The congruence subgroup property of branch groups has been studied comprehensively over the years, see \cite{BSZ12}, \cite{Gar16}, 
\cite{FAGUA17}. The generalised Basilica group $\bp_s(\mathcal{O}_m^{d})$ does not have the CSP as its abelianisation is isomorphic to {$\Z^{ds}$} (\cref{Theorem}{thm:gen_bas_properties}). However, the quotients of $\bp_s(\mathcal{O}_m^{d})$ by the layer stabilisers are isomorphic to subgroups of $C_m \wr \overset{n}{\cdots} \wr C_m$, for suitable $n \in \N$. If $m=p$, a prime, then these quotients are, in particular, finite $p$-groups. The class of all finite $p$-groups is a well-behaved class, i.e., it is closed under taking subgroups, quotients, extensions and direct limits. In light of this, we prove that $\bp_s(\mathcal{O}_p^{d})$ has the $p$-congruence subgroup property ($p$-CSP), a weaker version of CSP introduced by Garrido and Uria-Albizuri in \cite{GUA19}.  The group $G$ has the $p$-CSP if every subgroup of index a power of $p$ in $G$ contains some layer stabiliser in $G$. In \cite{GUA19} one finds a sufficient condition for a weakly branch group to have the $p$-CSP and it is also proved that the original Basilica group $\mathcal{B}$ has the $2$-CSP. This argument is generalised by Fernandez-Alcober, Di Domenico, Noce and Thillaisundaram to see that the $p$-Basilica groups have the $p$-CSP. We further generalise these result.

\begin{theorem} \label{thm:gen_bas_mCSP}
	For all $d,\,s \in \Nplus$ with $s > 2$, and all primes $p$, the generalised Basilica group $\bp_s(\mathcal{O}_p^{d})$ has the $p$-congruence subgroup property.
\end{theorem}

Even though we follow the same strategy as in \cite{GUA19}, the arguments differ significantly because of \cref{Theorem}{thm:gamma}. Here we make use of \cref{Theorem}{thm:stab} to obtain a normal generating set for the layer stabilisers of the generalised Basilica groups (\cref{Theorem}{thm:gen_bas_stab}). We remark that the result of Fernandez-Alcober, Di Domenico, Noce and Thillaisundaram on $p$-Basilica groups can be generalised to all $d \geq 2$ with additional work. 

The organisation of the paper is as follows: In \cref{Section}{sec:prelim}, we introduce the basic theory of groups acting on rooted trees and give the formal definition of the Basilica operation, together with important examples. The proofs of \cref{Theorem}{thm:transfered_prop} and \cref{Theorem}{thm:spinal} are given in \cref{Section}{sec:basics}. \cref{Theorem}{thm:split_transfer} and related results for $s$-split groups are contained in \cref{Section}{sec:split}, as well as the proofs of \cref{Theorem}{thm:stab} and \cref{Theorem}{thm:gen_bas_haus}. \cref{Section}{sec:gen_bas_Lpres} contains the proof of \cref{Theorem}{thm:gen_bas_Lpres}, while \cref{Section}{sec:gen_bas_prep} and \cref{Section}{sec:gen_bas_mCSP} contain the proofs of
\cref{Theorem}{thm:gamma} and 
\cref{Theorem}{thm:gen_bas_mCSP}.

\subsection*{Acknowledgements}
We carried out the work for this paper as PhD students under the supervision of Benjamin Klopsch at the Heinrich-Heine-Universität Düsseldorf. We thank him for introducing us to the subject and are grateful for his comments on earlier versions. We wish to express our thanks to Hanna Sasse for providing access to her Master thesis. We thank Alejandra Garrido and Anitha Thillaisundaram for helpful discussions, and we are grateful to the anonymous referee for the valuable suggestions.

\section{Preliminaries and Main Definitions}\label{sec:prelim}

\noindent For any two integers $i, j$, let $\intseg{i}{j}$ denote the interval in $\Z$. From here on, \(T_m = T\) denotes the \(m\)-regular rooted tree for an arbitrary but fixed integer $m > 1$. The vertices of \(T\) are identified with the elements of the free monoid \(X^*\) on \(X = \intseg{0}{m-1}\) by labeling the vertices from left-to-right. We denote the empty word by \(\epsilon\). For $n \in \N$, the \(n\)\emph{-th layer of} \(T\) is the set $X^n$ of vertices represented by words of length \(n\).

Every (graph) automorphism of \(T\) fixes \(\epsilon\) and moreover maps the \(n\)-th layer to itself for all $n \in \N$. The action of the full group of automorphisms \(\Aut(T)\) on each layer is transitive. A subgroup of $\Aut(T)$ with this property is called \emph{spherically transitive}. The stabiliser of a word \(u\) under the action of a group \(G\) of automorphisms of \(T\) is denoted by \(\st_G(u)\) and the intersection of all stabilisers of words of length \(n\) is called the \emph{\(n\)-th layer stabiliser}, denoted \(\St_G(n)\).

Let \(a \in \Aut(T)\) and let \(u, v\) be words. Since layers are invariant under \(a\), the equation
\[
	a(uv) = a(u)a|_u(v)
\]
defines a unique automorphism \(a|_u\) of \(T\) called the \emph{section of \(a\) at \(u\)}. This automorphism can be thought of as the automorphism induced by \(a\) by identifying the subtrees of \(T\) rooted at the vertices \(u\) and \(a(u)\) with the tree $T$. If \(G\) is a group of automorphisms, \(G|_u\) will denote the set of all sections of group elements at \(u\). The restriction of the action of the section \(a|_u\) to $X^1 = X$ is called the \emph{label of $a$ at $u$} and it will be written as \(a|^u\).

The following holds for all words \(u, v\) and all automorphisms \(a, b\):
\begin{align*}
	(a|_u)|_v &= a|_{uv},\\
	(ab)|_u &= a|_{b(u)}b|_u.
\end{align*}
The analogous identities hold for the labels \(a|^u\), so the action of \(a\) on any word \(x_0\dots x_{n-1}\) of length \(n\) is given by
\[
	a(x_0 \dots x_{n-1}) = a|^\epsilon(x_0)a|_{x_0}(x_1\dots x_{n-1}) = a|^\epsilon(x_0)a|^{x_0}(x_1)\dots a|^{x_0\dots x_{n-2}}(x_{n-1}).
\]
Hence every automorphism \(a\) is completely described by the label map $X^* \to \Sym(X),$ $ u \mapsto a|^u$, called the \emph{portrait of \(a\)}.

For $n \in \N$, the isomorphim
\begin{align*}
	\psi_n : \St(n) \to (\Aut(T))^{m^n}, \;g \mapsto (g|_x)_{x \in X^n},
\end{align*}
is called the \emph{$n$-th layer section decomposition}. We will shorten the notation of big tuples arising for example in this way by writing $g^{\ast k}$ for a sequence of $k$ identical entries $g$ in a tuple, implicitly ordering the vertices lexicographically.

We can uniquely describe an automorphism $g \in \Aut(T)$ by its label at $\epsilon$ and the first layer {section decomposition} of $(g|^\epsilon)^{-1}g$, i.e.\ by
\[
	g = g|^\epsilon \; (g|_x)_{x \in X}.
\]

Let $H \leq \Sym(X)$ be any subgroup of the symmetric group on $X$. Then denote by $\Gamma(H)$ the subgroup of $\Aut(T)$ defined as 
  \[
    \Gamma(H) = \LD a \in \Aut(T) \mid \forall u \in T, a|^u \in H \RD. 
  \]
If $H$ is a Sylow-$p$ subgroup of $\Sym(X)$, then $\Gamma(H)$ is a Sylow-pro-$p$ subgroup of $\Aut(T)$. We further fix $\sigma = (0\;1\;\dots\;m-1) \in \Sym(X)$ and write $\Gamma$ for $\Gamma(\LD\sigma\RD)$.

A group $G \leq \Aut(T)$ is called \emph{self-similar} if it is closed under taking sections at every vertex, i.e.\ if $G|_v \subseteq G$ for all $v \in T$. Self-similar groups correspond to certain automata modelling the behaviour of the section map: there is a state for every element $g \in G$, and an arrow $g \to g|_x$ labelled $x:g(x)$ for every $x \in X$ (for details see \cite{Nek05}).

We follow \cite{UA16} in the terminology for the first three of the following self-referential properties, and add a fourth one: A group $G \leq \Aut(T)$ acting spherically transitively is called
\begin{enumerate}
	\item \emph{fractal} if $\st_G(u)|_u = G$ for all $u \in T$.
	\item \emph{strongly fractal} if $\St_G(1)|_x = G$ for all $x \in X$.
	\item \emph{super strongly fractal} if $\St_G(n)|_u = G$ for all $n \in \N$ and $u \in X^n$.
	\item \emph{very strongly fractal} if $\St_G(n+1)|_x = \St_G(n)$ for all $n \in \N$ and $x \in X$.
\end{enumerate}

Notice that for every group $H$ acting regularly on $X$ and $G \leq \Gamma(H)$ the properties $(1)$ and $(2)$ coincide. The following lemma will be of great use.

\begin{lemma}\label{lem:fractal_branching}
	Let $G \leq \Aut(T)$ be fractal and self-similar, and let $x, y \in X$. For every $g \in G$ there exists an element $\tilde{g} \in G$ such that $\tilde{g}(x) = y$ and $\tilde{g}|_x = g$. Furthermore, if $H \leq G$ is any subgroup of $G$ such that $H \times \{\id\} \times \cdots \times \{\id\} \leq \psi_1(K)$ for some normal subgroup $K \trianglelefteq G$, then $(H^G)^m \leq \psi_1(K)$.
\end{lemma}

\begin{proof}
	Since $G$ is fractal, it is spherically transitive and in particular it is transitive on the first layer of $T$. Hence there exists some element $h \in G$ mapping $x$ to $y$. Also because $G$ is fractal and $h|_x \in G$ by self-similarity, there is some element $k \in \st_G(x)$ such that $k|_x = (h|_x)^{-1}g$. Now $\tilde{g} = hk$ fulfills both $\tilde{g}(x) = y$ and $\tilde{g}|_x = h|_xk|_x = g$.
	
	Assume further that $H \leq G$ and $H \times \{\id\} \times \cdots \times \{\id\} \leq \psi_1(K)$ for $K \trianglelefteq G$. Let $g \in G$. Choose an element $\tilde{g} \in G$ such that $\tilde{g}(x) = 0$ and $\tilde{g}|_x = g$. Then for every $h \in H$
	\[
		(\id^{\ast x}, h^g, \id^{\ast (m-x-1)}) = \psi_1((\tilde{g})^{-1}\psi_1^{-1}(h, \id, \dots, \id)\tilde{g}) \in \psi_1((\tilde{g})^{-1}K\tilde{g}) = \psi_1(K).\qedhere
	\]
\end{proof}
From this point on, we fix a positive integer $s$.

\begin{defn}\label{def:beta}
	There is a set of $s$ interdependent monomorphims $\beta^s_i:\Aut(T) \to \Aut(T)$ defined by
	\begin{align*}
		\beta^s_i(g) &= (\beta^s_{i-1}(g), \id, \dots, \id) &\text{ for } i \in \intseg{1}{s-1},\\
		\beta^s_0(g) &= g|^\epsilon (\beta^s_{s-1}(g|_0), \dots, \beta^s_{s-1}(g|_{m-1})).&
	\end{align*}
\end{defn}

\noindent We adopt the convention that the subscript for these maps is taken modulo $s$, whence $\beta^s_i(g)|_x \in \beta^s_{i-1}(\Aut(T))$ for all $i \in \intseg{0}{s-1}$ and $g \in \Aut(T)$. Whenever there is no reason for confusion, we drop the superscript $s$.

\begin{defn}\label{def:bp}
	Let $G \leq \Aut(T)$. The \emph{$s$-th Basilica~group of $G$} is defined as
	\[
		\bp_s(G) = \LD\beta^s_i(g)\mid g\in G, i \in \intseg{0}{s-1}\,\RD.
	\]
\end{defn}

\noindent Clearly, for $s = 1$ the homomorphism $\beta^1_0$ is the identity map and $\bp_1(G) = G$.
In the case of a self-similar group $G$, the $s$-th Basilica~group of $G$ can be equivalently defined as the self-similar closure of the group $\beta^s_0(G)$, i.e.\ the smallest self-similar group containing $\beta^s_0(G)$.
If $G$ is finitely generated by $g_1, \dots, g_r$, then $\bp_s(G)$ is generated by $\beta^s_i(g_j)$ with $i \in \intseg{0}{s-1}$ and $j \in \intseg{1}{r}$.

The operation $\bp_s$ is multiplicative in $s$, i.e.\ for $s,t \in \Nplus$ and $G \leq \Aut(T)$ we have $\bp_s\bp_t(G) = \bp_{st}(G)$. This is a consequence of
\[
	\beta^s_i(\beta^t_j (g)) = \beta^{st}_{i+sj}(g),
\]
which is an easy consequence of \cref{Definition}{def:beta}.

We now describe the monomorphisms $\beta^s_i$ for $i \in \intseg{0}{s-1}$ in terms of their portraits. We define a map $\omega_i: T \to T$. For every $k \in \N$ and every vertex $u \in X^k$, write $u = x_0\dots x_{k-1} \in X^k$, and define 
\[
	\omega_i(u) \defeq 0^{i}\prod_{j=0}^{k-2} \left( x_j0^{s-1} \right) x_{k-1}.
\]
Writing $\omega_i(T)$ for the subgraph of $T$ induced by the image of $\omega_i$, with edges inherited from paths in $T$, we again obtain an $m$-regular rooted tree.

\begin{lemma}\label{lem:portrait}
	Let $g \in \Aut(T)$ and $i \in \intseg{0}{s-1}$. Then the portrait of $\beta_i^s (g)$ is given by
	\begin{align*}
		\beta_i^s(g)|^u = \begin{cases}
			g|^{v}, &\text{ if }u = \omega_i(v),\\
			\id, &\text{ if } u \not \in \omega_i (T).
		\end{cases}
	\end{align*}
	In particular $\bp_s(G) \leq \Gamma(H)$, if $G \leq \Gamma(H)$ for some $H \leq \Sym(X)$.
\end{lemma}

\begin{proof}
	First suppose that $u = \omega_i(v)$ for $v = x_0\dots x_{k-1}$.
	From \cref{Definition}{def:beta} follows
	\[
		\beta_i^s(g)|^{\omega_i(x_0\dots x_{k-1})} = \beta_0^s(g)|^{\omega_0(x_0\dots x_{k-1})} = \beta_{s-1}^s(g|_{x_0})|^{\omega_{s-1}(x_1\dots x_{k-1})},
	\]
	and iteration establishes $\beta_i^s(g)|^u = g|^v$. Now, if $u = u_0\dots u_{k-1} \not\in \omega_i(T)$, there is some minimal number $n \not\equiv_s i$ such that $u_n \neq 0$. Thus $u = \omega_i(v)0^t u_n \dots u_{k-1}$ for $n \equiv_s t < i$ and some vertex $v$, hence
	\[
		\beta_i^s(g)|^u = \beta_i^s(g|_v)|^{0^t u_n \dots u_k} = \beta_{i-t}^s(g|_v)|^{u_n \dots u_k} = \id.\qedhere
	\]
\end{proof}

It is interesting to compare the effect of the Basilica operation with another method of deriving new self-similar groups from given ones described by Nekrashevych.

\begin{prop}[\cite{Nek05}*{Proposition 2.3.9}]\label{prop:direct_prod_of_self_sim}
	Let $G \leq \Aut(T)$ be a group and let $d$ be a positive integer. There is a set of $d$ injective endomorphisms of $\Aut(T)$ given by
	\begin{align*}
		\pi_0(g) &\defeq g|^\epsilon \; (\pi_{d-1}(g|_x))_{x \in X}, \\
		\pi_i(g) &\defeq (\pi_{i-1}(g))_{x \in X}  \text{ for } i \in \intseg{1}{d-1}.
	\end{align*}
	The group $D_d(G) \defeq \LD \pi_i(G) \mid i \in \intseg{0}{d-1} \RD$ is isomorphic to the direct product $G^d$.
\end{prop}

We combine both constructions to define a class of groups very closely resembling the original Basilica group $\mathcal{B}$.

\begin{defn}\label{def:gen_basilica}
	Let $d,\,m,\,s \in \Nplus$ with $m \geq 2$. The \emph{$m$-adic odometer} $\mathcal O_m$ is the infinite cyclic group generated by
	\[
		a = \sigma (a, \id, \dots, \id),
	\]
	where $\sigma$ is the $m$-cycle $(m-1\;m-2\;\dots\;1\;0)$.
	Write $\mathcal O_{m}^d$ for $D_d(\mathcal O_m)$, the $d$-fold direct product of $\mathcal O_m$ embedded into $\Aut(T)$ by the construction described in \cref{Proposition}{prop:direct_prod_of_self_sim}. We call the group $\bp_s(\mathcal O_{m}^d)$ the \emph{generalised~Basilica~group}.
\end{defn}

Clearly, $\mathcal B = \bp_2(\mathcal O_2)$ is the original Basilica group introduced by Grigorchuk and Żuk in \cite{GZ02}.

For illustration we depict explicitly the automaton defining the self-similar action of the dyadic odometer $\mathcal O_2$, the automaton defining the action of $D_8(\mathcal O_2)$ described above and the automaton defining $\bp_8(\mathcal O_2)$ in \cref{Figure}{fig:automatons}.

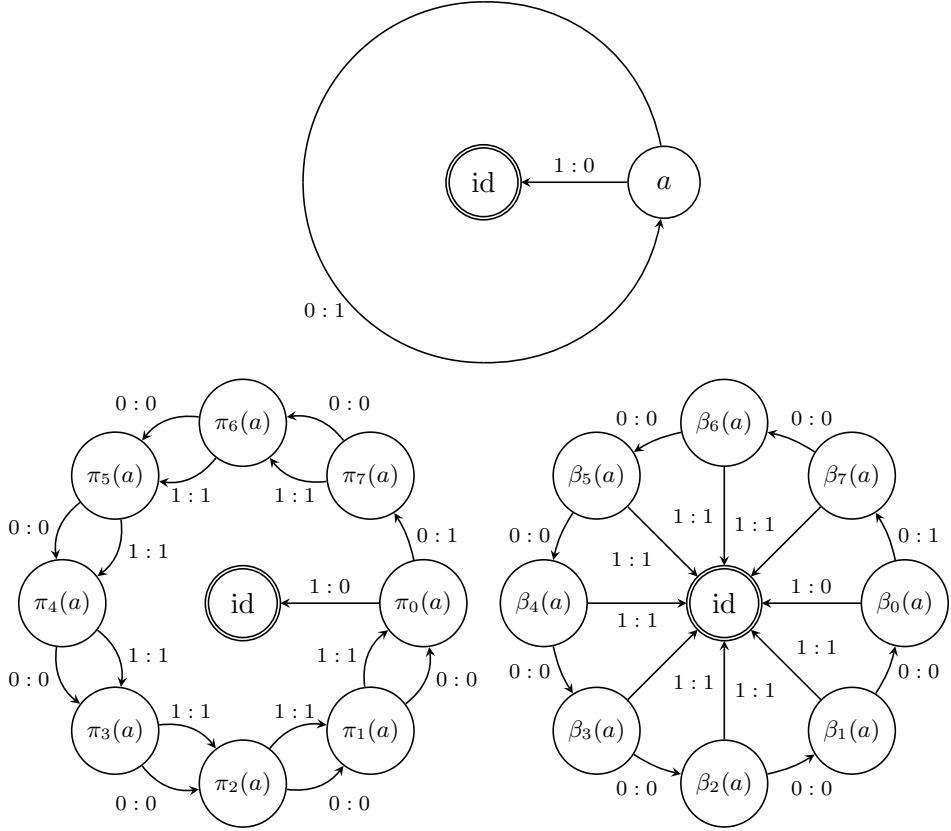
\begin{figure}
	\centering
	\begin{tikzpicture}[>=stealth,node distance=2.5cm,semithick,->,scale=0.8,font=\scriptsize]
	
		\node[state,accepting, font=\normalsize]	(I)		at (4,7)	{$\id$};
		\node[state, font=\normalsize]	(A)		at (4+3,7)	{$a$};
		
		\path	
				(A)		edge	[above]					node {$1:0$}			(I)
				(A)		edge	[bend right=40,-]				 			(4,7+3)
				(4,7+3)	edge	[bend right=45,-]				 			(4-3,7)
				(4-3,7)	edge	[bend right=45,-]		node [left] {$0:1$}	(4,7-3)
				(4,7-3)	edge	[bend right=40]			 					(A);
		
		\node[state,accepting, font=\normalsize]	(I)	 at (8,0)						{$\id$};
		\node[state, font=\footnotesize]	(A1) at (8+3,0)						{$\beta_0(a)$};
		\node[state, font=\footnotesize]	(A2) at ({8+3/sqrt(2)},{3/sqrt(2)})	{$\beta_7(a)$};
		\node[state, font=\footnotesize]	(A3) at (8,3)						{$\beta_6(a)$};
		\node[state, font=\footnotesize]	(A4) at ({8-3/sqrt(2)},{3/sqrt(2)})	{$\beta_5(a)$};
		\node[state, font=\footnotesize]	(A5) at (8-3,0)						{$\beta_4(a)$};
		\node[state, font=\footnotesize]	(A6) at ({8-3/sqrt(2)},{-3/sqrt(2)})	{$\beta_3(a)$};
		\node[state, font=\footnotesize]	(A7) at (8,-3)						{$\beta_2(a)$};
		\node[state, font=\footnotesize]	(A8) at ({8+3/sqrt(2)},{-3/sqrt(2)})	{$\beta_1(a)$};
		
		\path	
			(A1)	edge node	[above]		  {$1:0$}		(I)
			(A2)	edge node	[above left]  {$1:1$}		(I)
			(A3)	edge node	[left]		  {$1:1$}		(I)
			(A4)	edge node	[below left]  {$1:1$}		(I)
			(A5)	edge node	[below]		  {$1:1$}		(I)
			(A6)	edge node	[below right] {$1:1$}		(I)
			(A7)	edge node	[right]		  {$1:1$}		(I)
			(A8)	edge node	[above right] {$1:1$}		(I)
			(A1)	edge [bend right=10] node	[right]			{$0:1$}		(A2)
			(A2)	edge [bend right=10] node	[near end, above right]	{$0:0$}		(A3)
			(A3)	edge [bend right=10] node	[near start, above left]	{$0:0$}		(A4)
			(A4)	edge [bend right=10] node	[left]	{$0:0$}		(A5)
			(A5)	edge [bend right=10] node	[left]			{$0:0$}		(A6)
			(A6)	edge [bend right=10] node	[near end, below left]	{$0:0$}		(A7)
			(A7)	edge [bend right=10] node	[near start, below right]			{$0:0$}		(A8)
			(A8)	edge [bend right=10] node	[right]	{$0:0$}		(A1);

		\node[state,accepting, font=\normalsize]	(I)	 at (0,0)						{$\id$};
		\node[state, font=\footnotesize]	(A1) at (3, 0)						{$\pi_0(a)$};
		\node[state, font=\footnotesize]	(A2) at ({3/sqrt(2)},{3/sqrt(2)})	{$\pi_7(a)$};
		\node[state, font=\footnotesize]	(A3) at (0,3)						{$\pi_6(a)$};
		\node[state, font=\footnotesize]	(A4) at ({-3/sqrt(2)},{3/sqrt(2)})	{$\pi_5(a)$};
		\node[state, font=\footnotesize]	(A5) at (-3,0)						{$\pi_4(a)$};
		\node[state, font=\footnotesize]	(A6) at ({-3/sqrt(2)},{-3/sqrt(2)})	{$\pi_3(a)$};
		\node[state, font=\footnotesize]	(A7) at (0,-3)						{$\pi_2(a)$};
		\node[state, font=\footnotesize]	(A8) at ({3/sqrt(2)},{-3/sqrt(2)})	{$\pi_1(a)$};
		
		\path	
			(A1)	edge node	[above]		  {$1:0$}		(I)
			(A2)	edge [bend left] node [below] {$1:1$} (A3)
			(A3)	edge [bend left] node [below] {$1:1$} (A4)
			(A4)	edge [bend left] node [right] {$1:1$} (A5)
			(A5)	edge [bend left] node [right] {$1:1$} (A6)
			(A6)	edge [bend left] node [above] {$1:1$} (A7)
			(A7)	edge [bend left] node [above] {$1:1$} (A8)
			(A8)	edge [bend left] node [left] {$1:1$} (A1)
			(A1)	edge [bend right=10] node	[right]			{$0:1$}		(A2)
			(A2)	edge [bend right] node [above right] {$0:0$} (A3)
			(A3)	edge [bend right] node [above left] {$0:0$} (A4)
			(A4)	edge [bend right] node [left] {$0:0$} (A5)
			(A5)	edge [bend right] node [left] {$0:0$} (A6)
			(A6)	edge [bend right] node [below left] {$0:0$} (A7)
			(A7)	edge [bend right] node [below right] {$0:0$} (A8)
			(A8)	edge [bend right] node [right] {$0:0$} (A1);
		
	\end{tikzpicture}
		
	\caption{The automata defining the generators of $\mathcal O_2$, $D_8(\mathcal O_2)$ and $\bp_8(\mathcal O_2)$.}\label{fig:automatons}
\end{figure}

We shall prove in the following (cf.\ \cref{Section}{sec:gen_bas_Lpres}, \cref{Section}{sec:gen_bas_prep}, \cref{Section}{sec:gen_bas_mCSP}) that generalised Basilica~groups resemble the original Basilica~group in many ways, justifying the terminology.

\begin{prop}\label{prop:fin_aut}
	Let $\Autf(T)$ be the group of all \emph{finitary automorphisms}, i.e.\ the group generated by all automorphisms $g_{\tau, v}$ for $v \in T$, $\tau \in \Sym(X)$ that have label $\tau$ at $v$ and trivial label everywhere else. For any $s \in \Nplus$
	\[
		\bp_s(\Autf(T)) = \Autf(T).
	\]
	On the other hand $\bp_s(\Aut(T))$ is not of finite index in $\Aut(T)$ for all $s > 1$.
\end{prop}

\begin{proof}
	Define for every $n \in \N$ a map $\mu_n: \Aut(T) \to \N$ by 
	$$
		\mu_n(g) = |\{ u \in X^n \mid g|_u \neq \id \}|.
	$$
   \cref{Lemma}{lem:portrait} shows that $g_{\tau, v} = \beta_i(g_{\tau, \omega_i^{-1}(v)}) \in \bp_s(\Autf(T))$ for every $v \in \bigcup_{i = 0}^{s-1} \omega_i(T)$. Conjugation with suitable elements produces all other generators, hence $\Autf(T)$ is contained in $\bp_s(\Autf(T))$. On the other hand, $\sum_{n \in \N} \mu_n(g) < \infty$ for any $g \in \Autf(T)$, implying that the same holds for all generators (and hence, all elements) of $\bp_s(\Autf(T))$. Thus, $\bp_s(\Autf(T)) = \Autf(T)$.
	
	For any $g \in \Aut(T)$ we have $\mu_n(g) \leq |X^n| = m^n$. But for all generators $\beta_i(g)$ of $\bp_s(\Aut(T))$ the stronger inequality
	\(
		\mu_{sn+i}(\beta_i(g)) \leq m^n
	\)
	holds, since $\beta_i(g)$ has trivial label at all vertices outside of $\omega_i(T)$.
	Let $g \in \Aut(T)$ and $q(g) \in \Q_+$ be the infimum of all numbers $r$ such that
	\[
		\limsup_{n \to \infty}\frac{\mu_{sn}(g)}{m^{(1+r)n}} = \infty.
	\]
	Then $g$ cannot be in $\bp_s(\Aut(T))$, since the inequality $\mu_n(ab) \leq \mu_n(a) + \mu_n(b)$ for $a,b \in \Aut(T)$ implies that it cannot be a finite product of the generators of $\bp_s(\Aut(T))$. By the same reason, all elements with different $q(g)$ are in different cosets. Since $q(\Aut(T)) = (0, s-1) \cap \Q$, the second statement follows.
\end{proof}

\begin{qst}
	In view of \cref{Proposition}{prop:fin_aut} and the original Basilica group $\mathcal B$ it seems plausible that the operation $\bp_s$ makes (in some vague sense) big groups smaller and small groups bigger. Let $H \leq \Sym(X)$ be a transitive subgroup. Write $\Gamma_{\mathrm{fin}}(H) = \Autf(T)\cap\Gamma(H)$. Replacing $\Autf(T)$ with $\Gamma_{\mathrm{fin}}(H)$ in the proof of \cref{Proposition}{prop:fin_aut} we obtain $\bp_s(\Gamma_{\mathrm{fin}}(H)) = \Gamma_{\mathrm{fin}}(H)$.
	
	\emph{Is there a group $G$ not of the form $\Gamma_{\mathrm{fin}}(H)$ such that $\bp_s(G) = G$?}
\end{qst}

\section{Properties inherited by Basilica groups}\label{sec:basics}

We recall our standing assumptions: $m$ and $s$ are positive integers with $m \neq 1$, $X = \intseg{0}{m-1}$, and $T$ the $m$-regular rooted tree. The subscript of the maps $\beta_i^s$ is taken modulo $s$, and we will drop the superscript $s$ from now on.

\subsection{Self-similarity and fractalness}\label{sub:self-sim_and_frac}

\begin{lemma}\label{lem:levtra}
	Let $G \leq \Aut(T)$ act spherically transitively on $T$. Then $\bp_s(G)$ acts spherically transitively on $T$.
\end{lemma}

\begin{proof}
	It is enough to prove that for any number $n = qs + r \in \Nplus$ with $r \in \intseg{0}{s-1}$ and $q \geq 0$, and $y \in X$ there is an element $b \in \bp_s(G)$ such that $b(0^n0) = 0^ny$. Let $g \in G$ be such that $g(0^q0) = 0^{q}y$ and observe that $\beta_r(g)$ stabilises $0^n$. By \cref{Lemma}{lem:portrait} it follows
	\[
		\beta_r(g)(0^n0) = 0^n\beta_0(g|_{0^q})(0) = 0^ny.\qedhere
	\]
\end{proof}

\begin{lemma}\label{lem:self-sim}
	Let $G \leq \Aut(T)$ be self-similar. Then $\bp_s(G) \leq \Aut(T)$ is self-similar.
\end{lemma}

\begin{proof}
	We check that $\beta_i(g)|_v$ is a member of $\bp_s(G)$ for all $v \in T$. This holds by \cref{Definition}{def:beta} for words $v$ of length $1$, and follows from $g|_x|_y = g|_{xy}$ by induction for words of any length.
\end{proof}

\begin{lemma}\label{lem:fractality}
	Let $G \leq \Aut(T)$ be self-similar, and fractal (resp.\ strongly fractal). Then
	\begin{enumerate}
		\item[(i)] The group $B = \bp_s(G) \leq \Aut(T)$ is fractal (resp.\ strongly fractal).
		\item[(ii)] For all $b \in B$ there is an element $c \in \st_B(0)$ (resp.\ $c \in \St_B(1)$) such that $c|_0 = b$ and $c|_x \in \beta_{s-1}(G)$ for all $x \in \intseg{1}{m-1}$.
	\end{enumerate}
\end{lemma}

\begin{proof}
	\cref{Lemma}{lem:levtra} shows that $B$ acts spherically transitively, and by \cref{Lemma}{lem:self-sim} the group $B$ is self-similar. First suppose that $G$ is fractal. Since the statement (ii) implies the statement (i), it is enough to prove (ii).
	
	Observe that 
	\[
		H = \{g \in \st_B(0) \mid g|_x \in \beta_{s-1}(G) \text{ for all } x \in \intseg{1}{m-1}\}
	\]
	is a subgroup since $h(x) \neq 0$ and $(gh)|_x = g|_{h(x)}h|_x \in \beta_{s-1}(G)$ for all $g, h \in H, x \in \intseg{1}{m-1}$. Thus it is enough to show that $\beta_i(G) \leq H|_0$ for all $i \in \intseg{0}{s-1}$.
	
	It is easy to see that $\beta_{i}(G) \leq H$ for $i \neq 0$, hence since $\beta_i(G)|_{0} = \beta_{i-1}(G)$ we have $\beta_i(G) \leq H|_0$ for $i \neq s-1$. But also $\beta_{0}(\st_G(0)) \leq H$. Note that, since $G$ is fractal, we have $\st_G(0)|_0 = G$. Hence $\beta_{s-1}(G) \leq \beta_{0}(\st_G(0))|_0 \leq H|_0$.
	
	If $G$ is strongly fractal, we may replace $H$ by its intersection with $\St_B(1)$ and $\st_G(0)$ by $\St_G(1)$ to obtain a proof for the analogous statement.
\end{proof}

\cref{Lemma}{lem:levtra}, \cref{Lemma}{lem:fractality} and \cref{Lemma}{lem:self-sim} yield proofs for the statements $(1), (2)$ and $(3)$ of \cref{Theorem}{thm:transfered_prop}.

\subsection{Amenability}\label{sub:amenability}

The original Basilica group $\mathcal B$ was the first example of an amenable, but not subexponentially amenable group. This had been conjectured already in \cite{GZ02}, where non-subexponentially amenability of $\mathcal B$ was proven. Amenability was proven by Bartholdi and Virág in \cite{BV05}. Later, Bartholdi, Kaimanovich and Nekrashevych proved that all groups generated from bounded finite-state automorphisms are amenable \cite{BKN10}, which includes $\mathcal B$. We recall the relevant definitions and then apply the result of  Bartholdi, Kaimanovich and Nekrashevych to a wider class of groups produced by the Basilica operation.

\begin{defn}\label{def:finite-state_et_bounded}
	An automorphism $f \in \Aut(T)$ is called
	\begin{enumerate}
		\item \emph{finite-state} if the set $\{f|_u \mid u \in T\}$ is finite, and
		\item \emph{bounded} if the sequence $\mu_n(f) \defeq |\{u \in X^n \mid f|_u \neq \id \}|$ is bounded.
	\end{enumerate}
\end{defn}

\begin{prop}\label{prop:finite_state_bounded_invariant}
	Let $G \leq \Aut(T)$ be generated from finite-state bounded automorphisms. Then $\bp_s(G)$ is also generated from finite-state bounded automorphisms.
\end{prop}

\begin{proof}
	It is enough to prove that for every finite-state bounded $f \in \Aut(T)$ and $i \in \intseg{0}{s-1}$ the element $\beta_i(f)$ is again finite-state and bounded. Notice that all sections of $f$ are of the form $\beta_j(f|_u)$ for some $u \in T$, hence there are only finitely many candidates and $\beta_i(f)$ is finite-state.
	Moreover, by \cref{Definition}{def:beta} $\mu_{n}(\beta_i(f)) = \mu_{\lfloor\frac {n-i} s \rfloor}(f)$, bounding $\mu_n(\beta_i(f))$.
\end{proof}

This proves statement $(6)$ of \cref{Theorem}{thm:transfered_prop}, and we use \cite{BKN10} to conclude:

\begin{coro}\label{coro:amenability}
	Let $G \leq \Aut(T)$ be generated by finite-state bounded automorphisms. Then $\bp_s(G)$ is amenable.
\end{coro}

\subsection{Spinal Groups}\label{sub:spinality}
A well-known class of subgroups of $\Aut(T)$ containing most known branch groups is the class of \emph{spinal groups}, containing both the first and the second Grigorchuk group, and all  $\mathsf{GGS}$\emph{-groups}. We use, with modifications for $\mathsf{GGS}$-groups, the definition given in \cite{BGS03}.

\begin{defn}[cf.\ \cite{BGS03}*{Definition 2.1
	}]\label{def:spinal}
	Let $R \leq \Sym(X)$, let $D$ be a finite group and let
	\[
		\omega = (\omega_{i,j})_{i \in \Nplus, j \in \intseg{1}{m-1}}
	\]
	be a family of homomorphisms $\omega_{i,j}: D \to \Sym(X)$. Identify $R$ with $\{ r (\id, \dots, \id) \mid r \in R\} \leq \Aut(T)$ and identify each $d \in D$ with the automorphism of $T$ given by
	\begin{align*}
		d|^w \defeq \begin{cases}
			\omega_{i, j}(d) &\text{if }w = 0^{i-1}j \text{ for } i\in \Nplus, j \in \intseg{1}{m-1},\\
			\id &\text{otherwise.}
		\end{cases}
	\end{align*}
	Suppose that the following holds:
	\begin{enumerate}
		\item The group $R$ and all groups $\LD \omega_{n, j}(D) \mid j \in \intseg{1}{m-1} \RD$, for $n \in \Nplus$, act transitively on $X$.
		\item For all $n \in \Nplus$,
			\[
				\bigcap_{i = n}^\infty\bigcap_{j = 1}^{m-1} \ker\omega_{i,j} = 1.
			\]
	\end{enumerate}
	Then $\langle R, D\rangle \leq \Aut(T)$ is called the \emph{spinal group acting on $T$ with defining triple} $(R, D, \omega)$.	
	The spinal group with defining triple $(R, D, \omega)$ is called a \emph{$\mathsf{GGS}$-group acting on $T$} if $\omega_{n, j} = \omega_{k, j}$ for all $n, k \in \Nplus$ and $j \in \intseg{1}{m-1}$.
\end{defn}

We now describe the Basilica groups of spinal groups. For this, we record the following lemma.

\begin{lemma}\label{lem:spinal_commuting}
	Let $i, j \in \intseg{0}{s-1}$ with $i \neq j$. Denote by $\st(\overline 0)$ the stabiliser of the infinite ray $\overline{0} \defeq \{0^i \mid i \in \N\}$ in $\Aut(T)$ (a so-called \textit{parabolic subgroup}). Then
	\(
		[\beta_i(\st(\overline 0)), \beta_j(\st(\overline 0))] = 1.
	\)
\end{lemma}

\begin{proof}
	We prove that for all $g_0, g_1 \in \st(\overline 0)$ the images $b_0 = \beta_i(g_0)$ and $b_1 = \beta_j(g_1)$ commute, using the fact that $\st(\overline 0)|_0 = \st(\overline 0)$.
	Assume without loss of generality that either $j > i > 0$ or $i = 0$. In the first case both $b_0$ and $b_1$ stabilise the $i$-th layer, we can consider
	\[
		\psi_{i}([b_0, b_1]) = ([b_0|_{0^i}, b_1|_{0^i}], \id^{\ast (m^i-1)})=([\beta_0(g_0), \beta_{j-i}(g_1)], \id^{\ast (m^i-1)}),
	\]
	and thus reduce to the second case. Suppose now that $i = 0$. Since the only non-trivial first layer section of $b_1$ is at the vertex $0$ and by assumption $b_0$ fixes this vertex,
	\[
		\psi_1([b_0, b_1]) = ([b_0|_0, b_1|_0], \id^{\ast (m-1)}).
	\]
	Since $b_0|_0, b_1|_0 \in \st(\overline 0)$, we conclude by infinite descent that $[b_0, b_1]$ fixes all vertices outside the ray $\overline 0$, thus acts trivially on the entire tree $T$.
\end{proof}

The elements $d \in D$ of a spinal group defined by $(R, D, \omega)$ can be characterised by the fact that they stabilise the infinite ray (or \qq{spine}) $\overline 0$ and $d|^x \neq \id$ implies that $x$ has distance precisely $1$ from $\overline 0$. Therefore it is easy to see that a Basilica group $B = \bp_s(G)$ of a spinal group $G$ acting on $T$ cannot act as a spinal group on $T$, as the elements $\beta^s_i(d)$ have non-trivial labels at vertices of distance $s$ from the ray $\overline 0$. However, the group $B$ acts as a spinal group on a tree obtained from $T$ by deletion of layers. 

Motivated from \cref{Examples}{eg:gupta_sidki} and \ref{eg:grigorchuks} below, we introduce the following notations.

There is an injection $\iota_s: (X^s)^\ast \to X^\ast$ given by
\[
	(x_{0,0} \cdots x_{0, s-1}) \cdots (x_{n-1,0} \cdots x_{n-1, s-1}) \mapsto x_{0,0} \cdots x_{n-1, s-1},
\]
whose image is the union $\bigcup_{n \in \N} X^{sn}$. The restriction map induces an injection
\[
	\iota^*_s: \Aut (X^\ast) \to \Aut ((X^s)^*),
\]
and clearly the image $\iota^*_s(\Aut(T))$ is
\[
	\Gamma(\Sym(X) \wr \dots \wr \Sym(X)) \leq \Aut((X^s)^*),
\]
where the permutational wreath product is iterated $s$ times. Recall that $\Gamma(H)$ for a permutation group $G$ denotes the subgroup of $\Aut(T)$ with every local action a member of $H$. Define for $i \in \intseg{0}{s-1}$
\begin{align*}
	\tau_i: \Sym(X) &\to \Sym(X) \wr \dots \wr \Sym(X)\\
	\rho &\mapsto \iota^*_s(g_{\rho, 0^i})|^\epsilon,
\end{align*}
where $g_{\rho, 0^i}$ is the automorphism with $g|^{0^i} = \rho$ and $g|^x = \id$ everywhere else. It is easy to see that for every transitive permutation group $H \leq \Sym(X)$ the group $\LD \tau_k(H) \mid k \in \intseg{0}{s-1} \RD$ is isomorphic to the $s$-fold iterated permutational wreath product $H \wr \dots \wr H$.

Now given a family of homomorphisms $(\omega_{i,j}: D \to \Sym(X))_{i\in\Nplus, j \in X\setminus\{0\}}$ we define a new family $\tilde\omega = (\tilde\omega_{i,j}: D^s \to \Sym(X^s))_{i\in\Nplus, j\in X^s \setminus\{0^s\} }$ by
\begin{align*}
	\tilde\omega_{n,j} = \begin{cases}
		\tau_i \circ \omega_{n, x} \circ \pi_i, &\text{ if } j = 0^i x 0^{s-i-1} \text{ for some } x \in \intseg{1}{m-1} \text{ and } i \in \intseg{0}{s-1},\\
		d \mapsto \id, d \in D^s, &\text{ otherwise,} 
	\end{cases}
\end{align*}
where $\pi_i: D^s \to D$ denotes the projection to the $(i + 1)$-th factor.

\begin{prop}\label{prop:spinal}
	Let $G$ be the spinal group on $T$ with defining triple $(R, D, \omega)$. Then $\iota_s^*(\bp_s(G))$ is the spinal group on $(X^s)^*$ with defining triple $(R \wr \dots \wr R, D^s, \tilde\omega)$, by the action of $\bp_s(G)$ on the $m^s$-regular tree $\delta_s T$ defined by the deletion of layers.
	
	If furthermore $G$ is a $\mathsf{GGS}$-group on $T$, $\iota_s^*(\bp_s (G))$ is a $\mathsf{GGS}$-group on $(X^s)^*$.
\end{prop}

\begin{proof}
	First consider the elements of the form $\beta_k(a)$, for $a \in R$, $k \in \intseg{0}{s-1}$. On $(X^s)^*$ this element acts as $\tau_k(a)$. Since $R$ is transitive, the images of $R$ generate $R \wr \dots \wr R$, and the first entry of the defining triple is described.
	
	We deal in a similar way with the sections $\beta_i(d|_{0^ky})$ of a directed element for every $d \in D, i \in \intseg{0}{s-1}, k \in \N, y \in X\setminus\{ 0 \}$. To obtain the first section decomposition of the action of $\beta_i(d|_{0^k})$ on $\delta_s T$ (which stabilises the first layer) we have to take sections of $\beta_i(d|_{0^k})$ at words $x = x_0 \dots x_{s-1}$ of length $s$ in $T$. Now by \cref{Lemma}{lem:portrait},
	\[
		\beta_i(d|_{0^k})|_{x} = \begin{cases}
			\beta_i(d|_{0^{k+1}}) &\text{ if }x = 0^s,\\
			\beta_i(\omega_{k+1, x_i}(d)) = \tau_i \omega_{k+1, x_i}(d) &\text{ if }x = 0^{i}x_i0^{s-i-1}, x_i \neq 0,\\
			\id &\text{ otherwise.}
		\end{cases}
	\]
	By \cref{Lemma}{lem:spinal_commuting} all pairs $\beta_i(d_1), \beta_j(d_2)$ with $d_1, d_2 \in D$, $i, j \in \intseg{0}{s-1}$ and $i \neq j$ commute. We identify $\beta_i(D)$ with the $(i+1)$-th direct factor of $D^s$. Thus $\bp_s(G)$ is generated by $R \wr \dots \wr R$ and $\langle \beta_i(D) \mid i \in \intseg{0}{s-1}\rangle \cong D^s$, where $(\id, \dots, \id, d_i, \id, \dots, \id) \in D^s$ acts on $\delta_s T$ by
	\[
		(\id, \dots, \id, d_i, \id, \dots, \id)|_{0^{ks}x} = \beta_i(d|_{0^k})|_{x},
	\]
	thus, the elements of $D^s$ are defined by the family $\tilde\omega$ of homomorphisms. 
	
	It remains to establish the two defining properties of spinal groups. Property (1) holds by the observation that 
	\[
		\LD \tilde\omega_{i,j}(D^s) \mid j \in \intseg{1}{m^s-1} \RD
	\]
	acts as $\LD \tau_k(\omega_{i,j}(D)) \mid j \in \intseg{1}{m-1}, k \in \intseg{0}{s-1} \RD$, hence $\LD \tilde\omega_{i,j}(D^s) \mid j \in \intseg{1}{m^s-1} \RD$ acts as the $s$-fold wreath product of $\LD \omega_{i,j}(D) \mid j \in \intseg{1}{m-1} \RD$, in particular, transitively on the first layer of $\delta_s T$.
	
	For (2) consider 
	\[
		\ker \tilde\omega_{n, j} = \begin{cases}
			\ker (\omega_{n, x}\circ\pi_i), &\text{ if } j = 0^i x 0^{s-i-1}, \text{ for some } x \in \intseg{1}{m-1}, i \in \intseg{0}{s-1}\\
			D^s, &\text{ else,}
		\end{cases}
	\]
	hence
	\[
		\bigcap_{j \in X^s\setminus\{0^s\}} \ker \tilde\omega_{n, j} = \left(\bigcap_{j \in X\setminus\{0\}} \ker \omega_{n,j}\right) \times \dots \times \left(\bigcap_{j \in X\setminus\{0\}} \ker \omega_{n,j}\right).
	\]
	Therefore we see that since (2) holds for $G$, (2) holds for $\bp_s(G)$.
	
	The statement regarding $\mathsf{GGS}$-groups follows directly from the description of the defining triple of $\bp_s (G)$.
\end{proof}

\cref{Proposition}{prop:spinal} yields \cref{Theorem}{thm:spinal}.

\begin{eg}\label{eg:gupta_sidki}
	One of the eponymous examples of a $\mathsf{GGS}$-group is the family of the Gupta--Sidki $p$-groups acting on the $p$-adic tree. In the language of spinal groups they are defined by the triple
	\[
		(\LD \sigma \RD, \LD \sigma \RD, (\sigma \mapsto \sigma, \sigma \mapsto \sigma^{-1}, \sigma \mapsto \id, \dots, \sigma \mapsto \id)_{i\in\Nplus}),
	\]
	or in usual notation by the generators $a = \sigma ( \id, \dots, \id), b = (b, a, a^{-1}, \id, \dots, \id)$. We can describe the generators of the second Basilica~group of the Gupta--Sidki $3$-group $\ddot \Gamma$ by
	\begin{align*}
			\beta^2_0(a) &= \sigma (\id, \id, \id) = a	& \beta^2_0(b) &= (\beta^2_1(b), \beta^2_1(a), \beta^2_1(a^{-1})), \\
			\beta^2_1(a) &= (a, \id, \id)				& \beta^2_1(b) &= (\beta^2_0(b), \id, \id).
	\end{align*}
	The automaton describing these generators is given explicitly in \cref{Figure}{fig:gs3}. By ordering $X^2$ reverse lexicographically, the action of the generators on $(X^2)^*$ is
	\begin{align*}
		\beta^2_0(a) &= (00\;10\;20)(01\;11\;21)(02\;12\;22) \\
		\beta^2_0(b) &= (\beta^2_0(b), \beta^2_0(a), \beta^2_0(a)^{-1}, \id, \dots, \id)\\
		\beta^2_1(a) &= (00\;01\;02)\\
		\beta^2_1(b) &= (\beta^2_1(b), \id, \id, \beta^2_1(a), \id, \id, \beta^2_1(a)^{-1}, \id, \id).
	\end{align*}
\end{eg}

\begin{eg}\label{eg:grigorchuks}
	The first Grigorchuk group $\mathcal G$ is the spinal group acting on the binary tree defined by $\mathrm C_2, \mathrm C_2^2$ and the sequence $\omega_{i,1}$ of (the three) monomorphisms $\mathrm C_2 \to \mathrm C_2^2$, where $\omega_{i,1} = \omega_{j,1}$ holds if and only if $i \equiv_3 j$. 
	Writing $a$ for the non-trivial rooted element and $b, c, d$ for the non-trivial directed elements, one has the descriptions
	\begin{align*}
		\begin{array}{llll}
			a = (0\;1) (\id, \id), &b = (c, a), &c = (d, a), &d = bc = (b, \id).
		\end{array}
	\end{align*}
	By \cref{Proposition}{prop:spinal} $\bp_2(\mathcal G)$ is a spinal group on the $4$-adic tree $(X^2)^*$, generated by the elements
	\begin{align*}
		\begin{array}{lllrll}
			\alpha &\defeq \beta^2_0(a) &= (0 \; 2)(1 \; 3), &
				\mathrm A &\defeq \beta^2_1(a) &= (0 \; 1), \\
			\beta &\defeq \beta^2_0(b) &= (\kappa,\alpha,\id,  \id), &
				\mathrm B &\defeq \beta^2_1(b) &= (\mathrm K,\id, \mathrm A,  \id), \\
			\kappa &\defeq \beta^2_0(c) &= (\delta, \alpha, \id,  \id), &
				\mathrm K &\defeq \beta^2_1(c) &= (\Delta, \id, \mathrm A, \id), \\
			\delta &\defeq \beta\kappa, &&
				\Delta &\defeq \mathrm {BK},&
		\end{array}
	\end{align*}
	where we identify $\intseg{0}{3}$ with $X^2$ by the reverse lexicographic ordering.
\end{eg}

\subsection{Contracting groups}\label{sub:contraction}

For this subsection we fix a self-similar group $G \leq \Aut(T)$ and some generating set $S$ of $G$, which yields a natural generating set $\bigcup_{i \in \intseg{0}{s-1}}\beta_i(S)$ for $B \defeq \bp_s(G)$.

The group $G \leq \Aut(T)$ is said to be \emph{contracting}, if there exists a finite set $\mathcal N \subset G$ (called a \emph{nucleus} of $G$) such that for all $g \in G$ there is an integer $k(g)$ such that $g|_v \in \mathcal N$ for all $v \in T$ with $|v| > k(g)$, where $|\cdot|$ denotes the word norm.
	
In this section we prove that a contracting group $G$ has contracting Basilica groups $B = \bp_s(G)$, considering the natural generating set for $B$. For this we define yet another length function,
the \emph{syllable length}, denoted by $\syl(b)$, of an element $b \in B$ as the word length w.r.t.\ the infinite generating set $\bigcup_{i \in \intseg{0}{s-1}}\beta_i(G)$, i.e.\ as
  \[
   \syl(b) \defeq \min \{ \ell \in \N \mid b = \prod_{j = 0}^{\ell - 1} \beta_{i_j}(g_j), \text{ with suitable }i_j \in \intseg{0}{s-1}, g_j \in G \},
  \]
where $\prod_{j = 0}^{\ell - 1} \beta_{i_j}(g_j)$ is a word representing $b$ in $B$ with respect to the generating set $\{ \beta_i(g) \mid i \in \intseg{0}{s-1}, g \in G \}$. Consequently, we will call a non-trivial element of the given generating set a \emph{syllable} and the corresponding index $i$ its \emph{type}. Since for every non-trivial element $b \in \beta_i(G)$ there is some $u \in X^{ns+i}$ for some $n \in \N$ such that $b|^u \neq \id$, while there is no $u \in T\setminus \bigcup_{n \in \N} X^{ns+i}$ such that $b|^u \neq \id$, the type of a syllable is unique. Since all sections of a syllable are either trivial or a syllable itself, the syllable length of a section of $b$ is at most $\syl(b)$.

We further define for every $g \in \Aut(T)$,
\[
	\mathfrak r(g) \defeq \begin{cases}
		\min \{ n \in \N \mid g|^{0^{n}}(0) \neq 0 \} &\text{ if $g$ does not stabilise } \overline{0} = \{0^n \mid n \in \N\},\\
		\infty &\text{ otherwise.}
	\end{cases}
\]

\begin{lemma}\label{lem:technical part of syllable reduction}
	Let $r \in \N$. Define
	\begin{align*}
		D_r \defeq \{ \beta_{a_1}(h_1)\beta_{a_2}(h_2)\beta_{a_3}(h_3) \mid\; &h_1, h_2, h_3 \in G\backslash \{1\},\\ &a_1, a_2, a_3 \in \intseg{0}{s-1} \text{ pairwise distinct, }\\ &\mathfrak r(\beta_{a_2}(h_2)) = r \}.
	\end{align*}
	Then $\syl(c|_u) < 3$ for $c \in D_r$ and all $u$ with $|u| > r$.
\end{lemma}

\begin{proof}
	Let $c = \beta_{a_1}(h_1)\beta_{a_2}(h_2)\beta_{a_3}(h_3) \in D_r$, where $a_1, a_2, a_3, h_1, h_2, h_3$ satisfy the conditions stated above. We use induction on $r$. First consider the case $r = 0$. From $\beta_{a_2}(h_2)(0) \neq 0$ we deduce that $a_2 = 0$. Calculate, for $x \in \intseg{0}{m-1}$,
	\begin{align*}
		c|_x &= \begin{cases}
			\beta_{s-1}(h_2|_0)\beta_{a_3-1}(h_3) &\text{ if }x = 0,\\
			\beta_{a_1-1}(h_1)\beta_{s-1}(h_2|_x) &\text{ if }x = h_2^{-1}(0),\\
			\beta_{s-1}(h_2|_x) &\text{ otherwise.}
		\end{cases}
	\end{align*}
	This shows that $c|_x$ and, by recursion, $c|_u$ for all $u$ with $|u| \geq 1$ have syllable length at most $2$.
	Now we assume that $r > 0$. We may reduce to the case that $0 \in \{ a_1, a_2, a_3 \}$. If $0 \notin \{ a_1, a_2, a_3 \}$, $c|_0 \in D_{r-1}$ and $c|_x = \id$ for all $x \in X, x \neq 0$. Therefore, by induction $\syl(c|_x|_u) < 3$ for $x \in X$ and $|u| > r-1$, hence $\syl(x|_u) < 3$ for all $|u| > r$.
	
	If $a_3 = 0$, respectively $a_1 = 0$, we have
	\begin{align*}
		c|_x &= \begin{cases}
			\beta_{a_1-1}(h_1)\beta_{a_2-1}(h_2)\beta_{s-1}(h_3|_x) \in D_{r-1} &\text{ if }x = h_3^{-1}(0),\\
			\beta_{s-1}(h_3|_x) &\text{ otherwise,}
		\end{cases}\\
		\text{respectively }
		c|_x &= \begin{cases}
			\beta_{s-1}(h_1|_0)\beta_{a_2-1}(h_2)\beta_{a_3-1}(h_3) \in D_{r-1} &\text{ if }x = 0,\\
			\beta_{s-1}(h_1|_x) &\text{ otherwise.}
		\end{cases}
	\end{align*}
	In both cases all but one section have length $< 3$ and the remaining section is contained in $D_{r-1}$, hence by induction $\syl(c|_{xu}) < 3$ for all $x \in X, |u| > r - 1$.
	
	The case $a_2 = 0$ remains. Now $r > 0$ implies $h_2^{-1}(0) = 0$ and we have $\mathfrak r (\beta_{s-1}(h_2|_0)) = r - 1$. Thus
	\begin{align*}
		c|_x &= \begin{cases}
			\beta_{a_1-1}(h_1)\beta_{s-1}(h_2|_0)\beta_{a_3-1}(h_3) \in D_{r-1} &\text{ if }x = 0,\\
			\beta_{s-1}(h_2|_x) &\text{ otherwise.}
		\end{cases}
	\end{align*}
	Hence we conclude that $\syl(c|_{xu}) < 3$ for all $u$ with $|u| \geq 1$ by induction as before.
\end{proof}

\begin{lemma}\label{lem:syllable_red}
	For every element $b \in B$ with $\syl(b) > s + 1$ there is a number $r \in \N$ such that for all sections $b|_u$ with $|u| > r$,
	\[
		\syl(b|_u) < \syl(b).
	\]
\end{lemma}

\begin{proof}
	Let $b \in B$ be an element with $\syl(b) > s+1$. If $b$ is minimally represented by a word $w$, it suffices to prove that there is a subword of $w$ representing an element which has a reduction of the syllable length upon taking sections.
	
	Since $\syl(b) > s + 1$ there must be at least one syllable type appearing twice, and there is a subword of $w$ that can be written in the form
	\begin{align*}
		\beta_i(\tilde g_1) b_0 \beta_i(\tilde g_2)b_1 \text{ or } b_1\beta_i(\tilde g_1) b_0 \beta_i(\tilde g_2),
	\end{align*}
	where $b_0, b_1$ are non-trivial and contain neither two syllables of the same type nor a syllable of type $i$. Passing to the inverse if necessary we restrict to the first case.
	
	Under the assumption of $w$ being minimal it is impossible that both $b_0$ and $\beta_i(\tilde g_2)$ fix the infinite ray $\overline{0}$, since if they did, they would commute by \cref{Lemma}{lem:spinal_commuting}, and consequently it would be possible to reduce the number of syllables.
	
	Thus there are syllables in $b_0\beta_i(\tilde g_2)$ that do not stabilise the ray $\overline{0}$. Among these we choose $k$ such that $r \defeq \mathfrak r(\beta_{j_k}(g_k))$ is minimal.
	
	Apply \cref{Lemma}{lem:technical part of syllable reduction} to the subword $\beta_{j_{k - 1}}(g_{k - 1}) \beta_{j_{k}}(g_k) \beta_{j_{k + 1}}(g_{k + 1})$ of $\beta_i(\tilde g_1) b_0 \beta_i(\tilde g_2)b_1$ consisting only of the syllable $\beta_{j_{k}}(g_k)$ and its direct neighbours, and obtain for all $u \in T, |u| > r$
	\[
		\syl(b|_u) < \syl(b).\qedhere
	\]	
\end{proof}

Although interesting in its own right we use \cref{Lemma}{lem:syllable_red} solely to prove the following proposition.

\begin{prop}\label{prop:contrating}
	Let $G \leq \Aut(T)$ be contracting. Then $B = \bp_s(G)$ is contracting.
\end{prop}

\begin{proof}
	Let $\mathcal N(G)$ be a nucleus of $G$. Define
	\[
		\mathcal N(B) \defeq \left\{ \prod_{i = 0}^{\ell} \beta_{j_i}(g_i) \mid \ell \leq s + 1, j_i \in \intseg{0}{s-1}, g_i \in \mathcal N(G) \right\}.
	\]
	Since $\mathcal N(G)$ is a finite set, $\mathcal N(B)$ is finite as well. We will prove that it is a nucleus of $B$. Let $b \in B$. If $\syl(b) > s+1$, by \cref{Lemma}{lem:syllable_red} there is a layer, from which onwards all sections of $b$ have syllable length $s + 1$ or smaller.
	
	Hence we can assume, that $\syl(b) \leq s + 1$. Write
	\(
		b = \prod_{i = 0}^{\syl(b)-1} \beta_{j_i}(g_i).
	\)
	Since $G$ is contracting, for every $g_i$ there is a number $k(g_i)$ such that $g_i|_u \in \mathcal N(G)$ for all $|u| \geq k(g_i)$. Set $K \defeq \max \{ k(g_i) \mid i \in \intseg{0}{\syl(b)-1} \}$, and observe that for $u$ with $|u| \geq sK$ the section $b|_u$ is a product of at most $\syl(b) \leq s + 1$ syllables of the form $\beta_i(g)$ with $g \in \mathcal N(G)$. Thus $b|_u$ is in $\mathcal N(B)$ and $B$ is contracting.
\end{proof}

\cref{Proposition}{prop:contrating} proves statement $(4)$ of \cref{Theorem}{thm:transfered_prop}.

As a consequence, the word problem for Basilica groups of self-similar and contracting groups is solvable, since it is solvable for self-similar and contracting groups \cite{Nek05}*{Proposition 2.13.8}.

\begin{coro}\label{cor:word_problem}
	Let $G$ be self-similar and contracting. Then $\bp_s(G)$ has solvable word problem.
\end{coro}

\begin{qst}
	Let $G \leq \Aut(T)$ be contracting. The fact that $\bp_s(G)$ is contracting implies the existence of constants $\lambda < 1, L, C \in \mathbb R_+$ such that for every $g \in G$, $u \in X^n$ with $n > L$ it holds
	\[
		|g|_u| < \lambda |g| + C.
	\]
	In \cite{GZ02} one set of constants is given for the original Basilica group $\mathcal B$, namely $\lambda = \frac 2 3$ and $L = C = 1$.
	
	\emph{Is there a general formula for the above constants valid for all contracting groups and their Basilica groups, yielding $\lambda = \frac 2 3$ for $\mathcal B$?}
\end{qst}

\subsection{Word growth}\label{sub:growth}

We now provide some examples of the possible growth types of Basilica groups. It is known that the original Basilica group $\mathcal B$ has exponential word growth, cf.\ \cite{GZ02}*{Proposition 4}. The same proof as the one given there also shows that $\bp_2(\mathcal O_m)$ is of exponential growth for all $m \geq 2$. This, however, is not a general phenomenon.

\begin{prop}\label{prop:poly_growth_to_exp_and_poly_growth}
	Let $a = (0\;1)(a, \id)$ be the generator of the dyadic odometer acting on the binary rooted tree. Then \(\bp_s(\LD(\id, a)\RD)\) is a free abelian group of rank $s$, and is of polynomial growth in particular.
\end{prop}

\begin{proof}
	The element $(\id, a)$ stabilises the ray $\overline 0$, thus by \cref{Lemma}{lem:spinal_commuting} we have
	\[
	[\beta_i(\LD (\id, a)\RD), \beta_j(\LD (\id, a)\RD)] = \id
	\]
	for distinct $i , \, j \in \intseg{0}{s-1}$. Also $\beta_i(\LD (\id, a)\RD) \cong \Z$ for all $i \in \intseg{0}{s-1}$.
\end{proof}

As another example, we prove that there is a group of intermediate word growth such that its second Basilica group has exponential word growth.

\begin{prop}\label{prop:inter_growth_to_exp}
	Let $G = \LD a = (1\;2\;3), b = (a, 1, b)\RD$ be the Fabrykowski--Gupta group \cite{FG85} acting on the ternary rooted tree, which is of intermediate growth according to \cite{BP09}. Then there exists an element $f \in \Aut(T)$ such that the group $\bp_2(G^f)$ is of exponential growth.
\end{prop}

\begin{proof}
	The Fabrykowski-Gupta group is a $\mathsf{GGS}$-group. In contrast to the Gupta--Sidki $3$-group it is not periodic: an example for an element of infinite order is $ab$, for which the relation
	\[
		(ab)^3 = (ab, ba, ba)
	\]
	holds. In view of the decomposition it is clear that $ab$ acts spherically transitively on $T$ and thus by a result of Gawron, Nekrashevych and Sushchansky \cite{GNS01} it is $\Aut(T)$-conjugate to the $3$-adic odometer group. Let $f \in \Aut(T)$ be an element such that $(ab)^f = (1\;2\;3)((ab)^f,1,1)$. Then the subgroup generated by $\beta_0((ab)^f)$ and $\beta_1((ab)^f)$ in $\bp_2(G^f)$ is isomorphic to the generalised Basilica~group $\bp_2(\mathcal O_3)$, which is of exponential growth by following the proof of \cite{GZ02}*{Proposition 4} (which is the same result for $\mathcal B$) replacing the $2$-cycle with a $3$-cycle corresponding to $a|^\epsilon$.
\end{proof}

The same idea can be used to obtain the following proposition.

\begin{prop}
	Let $G \leq \Aut(T)$ be a group containing an element acting spherically transitively on $T$. Then there is an $\Aut(T)$-conjugate $G^f$ of $G$ such that $\bp_s(G^f)$ has exponential word growth.
\end{prop}

\subsection{Weakly Branch Groups}\label{sub:branchness}

For every vertex $v \in T$ the \emph{rigid vertex stabiliser of $v$ in $G$} is the subgroup of all elements that fix all vertices outside the subtree rooted at $v$. For every $n \in \N$ the \emph{$n$-th rigid layer stabiliser} $\Rst_G(n)$ is the normal subgroup generated by all rigid vertex stabilisers of $n$-th layer  vertices. A group $G \leq \Aut(T)$ is called a \emph{weakly branch group}, if $G$ acts spherically transitively and all rigid layer stabilisers $\Rst_G(n)$ are non-trivial.
If there is a subgroup $H \leq G$ such that $\psi_1(\St_H(1)) \geq H \times \dots \times H$, the group $G$ is said to be \emph{weakly regular branch over $H$}. Clearly, a group that is weakly regular branch group over a non-trivial subgroup is a weakly branch group.

From \cref{Lemma}{lem:portrait}, it follows that elements of the rigid layer stabilisers of $G$ translate to elements of rigid layer stabilisers of $\bp_s(G)$.
\begin{lemma}\label{lem:rigid_stab}
	Let $n = qs + r \in \N$, with $r \in \intseg{0}{s-1}$ and $q \geq 0$. Let $B = \bp_s(G)$ for $G \leq \Aut(T)$. Then $\Rst_B(n)$ contains $\beta_i(\Rst_G(q + 1))$ and $\beta_j(\Rst_G(q))$ for $0 \leq i < r$ and  for $r \leq j < s$.
\end{lemma}

We immediately obtain the following proposition.

\begin{prop}\label{prop:(weakly)_branch_transfers}
	Let $G \leq \Aut(T)$ be a weakly branch group. Then $B \defeq \bp_s(G)$ is again weakly branch.
\end{prop}

This proves the statement $(5)$ of \cref{Theorem}{thm:transfered_prop}.

The group $\bp_s(G)$ can be weakly branch even when $G$ is not weakly branch. We recall that for any group $G$ and an abstract word $\omega$ on $k$ letters, the set of $\omega$-elements and the verbal subgroup associated to $\omega$ are
\begin{align*}
	G_\omega \defeq \{\omega(h_0, \dots, h_{k-1}) \mid h_0, \dots, h_{k-1} \in G\} \text{ and } \omega(G) \defeq \LD G_\omega\RD \text{ respectively}.
\end{align*}

\begin{prop}\label{prop:weakly_branch_by_verbal}
	Let $G \leq \Aut(T)$ be a self-similar strongly fractal group and let $B \defeq \bp_s(G)$. Let $\omega$ be a law in $G$, i.e.\ a word $\omega$ such that $\omega(G) = 1$, but let $\omega$ not be a law in $B$. Then $B$ is weakly regular branch over $\omega(B)$.
\end{prop}

\begin{proof}
	Let $b = \omega(b_0, \dots, b_{k-1}) \neq \id$ with $b_i \in B$ for $i \in \intseg{0}{k-1}$. By \cref{Lemma}{lem:fractality} there are elements $c_i \in \St_B(1)$ such that $c_i|_{0} = b_i$ and $c_i|_x \in \beta_{s-1}(G)$ for all $x \in X\setminus \{0\}$.
	
	For every $x \in X$, let $d_x \in B$ be an element such that $d_x|_x = \id$ and $d_x(x) = 0$ (cf.\ \cref{Lemma}{lem:fractal_branching}). Then $c_i^{d_x}$ stabilises the first layer and has sections $c_i^{d_x}|_x = b_i$ and $c_i^{d_x}|_y = {(c_i|_{d_{x(y)}})}^{d_x|_y} \in \beta_{s-1}(G)^{d_x|_y}$ for $y \neq x$.

	Since $c_i^{d_x}$ stabilises the first layer, the section maps are homomorphisms and
	\[
		\omega(c_0^{d_x}, \dots, c_{k-1}^{d_x})|_y = \omega(c_0^{d_x}|_y, \dots, c_{k-1}^{d_x}|_y) = \begin{cases}
			b, &\text{ if } y = x\\
			\id &\text{ else,}
		\end{cases}
	\]
	because in the second case we are evaluating $\omega$ in a group isomorphic to $G$. This shows that $B_\omega \times \dots \times B_\omega$ is geometrically contained in $B_\omega$, and thus the same holds for the verbal subgroups that are generated by these sets.
\end{proof}

We point out that, if $\omega$ is a law in $B$, then $B$ cannot be weakly branch as it satisfies an identity. \cref{Proposition}{prop:weakly_branch_by_verbal} allows to obtain examples of groups that are weakly branch over some prescribed verbal subgroup. We provide an easy example:

\begin{eg}
	The group $D \defeq \LD \sigma, b\RD$, with $\sigma = (0 \; 1)$ and $b = (b, \sigma)$, acting on the binary tree is isomorphic to the infinite dihedral group (hence metabelian). It is self-similar and strongly fractal. Considering
	\[
		[[\beta_1(\sigma), \beta_0(\sigma)], [\beta_0(\sigma), \beta_0(\sigma b)]] = ([\beta_0(\sigma), \beta_1(\sigma b)], [\beta_0(\sigma), \beta_1(b^{-1} \sigma)]) \neq \id,
	\]
	we see that the second Basilica $\bp_2(D)$ is not metabelian, and thus it is weakly branch over the second derived subgroup of $\bp_2(D)$.
\end{eg}

\section{Split groups, Layer Stabilisers and Hausdorff dimension}\label{sec:split}

The subgroup $\beta_i(G) \leq \bp_s(G)$, for $i \in \intseg{0}{s-1}$, has the property that its elements have non-trivial portrait only at vertices at levels $n \equiv_s i$ for $n \in \N$.

We consider an algebraic analogue of this property that will be used to determine the structure of the stabilisers of $\bp_s(G)$.

\begin{defn}\label{def:s-split}
	Let $G \leq \Aut(T)$ and $B \defeq \bp_s(G)$. Define:
	\begin{align*}
		S_{i} \defeq \LD \beta_j(G) \mid j \neq i \RD \leq B \text{ and } N_{i} \defeq (S_{i})^B \trianglelefteq B.
	\end{align*}
	We write $\phi_i: B \to B/N_i$ for the canonical epimorphism with kernel $N_{i}$. The quotient $B/N_i$ is isomorphic to the quotient of $G$ by the normal subgroup $K_i \defeq \beta_{i}^{-1}(\beta_i(G) \cap N_i)$. We call $K_i$ the \textit{$i$-th splitting kernel} of $G$.
	The group $G$ is called \emph{$s$-split} if its $s$-th Basilica group $B$ is a split extension of $N_{i}$ by $\beta_i(G)$ for all $i \in \intseg{0}{s-1}$, or equivalently if all splitting kernels of $G$ are trivial.
\end{defn}

\begin{prop}\label{prop:non-abelian not split}
	Let $G \leq \Aut(T)$ be a group that does not stabilise the vertex $0$. Then $\beta_i([G,G]) \leq N_i$ for $i \in \intseg{1}{s-1}$. In particular, an $s$-split group (for $s > 1$) is abelian.
\end{prop}

\begin{proof}
	Let $g, h \in G$, $k \in G\setminus\st(0)$ and let $i \in \intseg{1}{s-1}$. Write $\gamma = \beta_{i-1}(g), \eta = \beta_{i-1}(h), \overline{\gamma} = \beta_i(g), \overline{\eta} = \beta_i(h)$ and $\kappa = \beta_0(k)$. Then
	\begin{align*}
		\kappa^{-1} (\kappa)^{\overline{\gamma}^{-1}} (\kappa^{-1})^{\overline{\gamma}^{-1}\overline{\eta}}(\kappa)^{\overline{\eta}}|_x &=
		\kappa^{-1}|_{\kappa(x)} \overline{\gamma}|_{\kappa(x)} \kappa|_x (\overline{\gamma}^{-1} \overline{\eta}^{-1} \overline{\gamma})|_x \kappa^{-1}|_{\kappa(x)} \overline{\gamma}^{-1}|_{\kappa(x)} \kappa|_x \overline{\eta}|_x\\&=
		\kappa|_x^{-1} \overline{\gamma}|_{\kappa(x)} \kappa|_x (\overline{\gamma}^{-1} \overline{\eta}^{-1} \overline{\gamma})|_x \kappa|_x^{-1} \overline{\gamma}|_{\kappa(x)}^{-1} \kappa|_x \overline{\eta}|_x\\
		&=
		\begin{cases}
			[\gamma, \eta] &\text{ if }x = 0,\\
			\id &\text{ otherwise}.
		\end{cases}
	\end{align*}
	Thus $\kappa^{-1} (\kappa)^{\overline{\gamma}^{-1}} (\kappa^{-1})^{\overline{\gamma}^{-1}\overline{\eta}}(\kappa)^{\overline{\eta}} = ([\gamma, \eta], \id, \dots, \id) = [\overline{\gamma}, \overline{\eta}]$ is an element of $N_i \cap \beta_i(G)$.
\end{proof}

We remark that $[G, G] \leq K_0$ does not necessarily hold. For example, consider a group $G$ such that $[G, G] \not\leq \St_G(1)$. Since $N_0 \leq \St_{\bp_s(G)}(1)$, the zero-th splitting kernel can not contain $[G,G]$.

\begin{defn}\label{def:non-ab}
	We call a subgroup $H$ of a group $G$ \emph{non-absorbing in $G$} if for all $h_0, \dots, h_{m-1} \in H$ such that $\psi_1^{-1}(h_0, \dots, h_{m-1}) \in G$, implies $\psi_1^{-1}(h_0, \dots, h_{m-1}) \in H$. If $G$ is weakly branch over $H$, then $H$ is non-absorbing in $G$.
\end{defn}

\begin{prop}\label{prop:split-intersection_in_commutator}
	Let $G \leq \Aut(T)$ be self-similar and such that $G|^\epsilon$ acts regularly on $X$. Assume that $[G, G]$ is non-absorbing in $G$. Then for $i \in \intseg{1}{s-1}$ we have $K_i = [G, G]$, and $K_0 \leq [G, G]$. In particular, if $G$ is abelian, it is $s$-split for all $s \in \Nplus$.
\end{prop}

\begin{proof}
	The inclusion $[G, G] \leq K_i$ for $i \in \intseg{1}{s-1}$ is proven in \cref{Proposition}{prop:non-abelian not split}. Thus we prove $K_i \leq [G, G]$ for $i \in \intseg{0}{s-1}$.
	
	Set $B \defeq \bp_s(G)$ and define $\mathcal N \defeq \bigcup_{i = 0}^{s-1} (\beta_i(G) \cap N_i)$. We employ the decomposition in syllables, cf.\ \cref{Subsection}{sub:contraction}. For every $b \in \mathcal N$ there is an index $i \in \intseg{0}{s-1}$ such that $b$ can be written both as an element of the image of some $\beta_i$ and a word in $N_i$, i.e.
	\begin{equation}\label{eq:dec}
		b = \beta_i(g_0) = \prod_{j=1}^{\ell(b)} (h_j)^{\beta_i(g_j)}\tag{$\ast$}
	\end{equation}
	for suitable $\ell(b) \in \N, \, g_j \in G$ and $h_j \in S_i$. The minimal possible value of $\ell(b)$ is called the \textit{restricted syllable length}, and from here onwards we use the symbol $\ell$ for this invariant. Write $\mathcal C = \bigcup_{i = 0}^{s-1} \beta_{i}([G, G])$ (notice that this a union of subsets with pairwise trivial intersection), and define 
	$$
	\mathcal M \defeq \{ b \in \mathcal N\setminus\mathcal C \mid \ell(b) \leq \ell(c) \text{ for all } c \in \mathcal N\setminus\mathcal C\},
	$$
	the set of all non-commutator elements with minimal restricted syllable length.
	
	We shall prove that for every $b \in \mathcal M$ there exists a first level vertex $x_i \in X$ such that:
	\begin{enumerate}
		\item $b|_{x_i} \in \mathcal M$ and
		\item $b|_x = \id$ for all $x \in X\setminus\{x_i\}$.
	\end{enumerate}
	Furthermore we prove that
	\begin{enumerate}
		\setcounter{enumi}{2}
		\item $b \in \St_{B}(1)$, i.e. $\mathcal M \leq \St_B(1)$.
	\end{enumerate}
	Every subset $\mathcal M \subseteq \Aut(T)$ with these properties is empty. Indeed, if $b \in \mathcal M$, there is some vertex $u \in T$ such that $b|^u \neq \id$, since $b$ is not trivial. But by properties $(1)$ and $(2)$ $b|_u$ is either trivial or a member of $\mathcal M$, hence by property $(3)$ stabilises the first layer, a contradiction.
	
	But if $\mathcal M$ is empty, $\mathcal N$ is contained in $\mathcal C$, hence all splitting kernels are subgroups of $[G, G]$, finishing the proof.
	
	Assume that there is some $b \in \mathcal M$. We fix the decomposition and the type given by (\ref{eq:dec}), but write $\ell$ for $\ell(b)$ to shorten the notation.
	
	We first observe that $\ell \neq 1$. If $\ell = 1$, we have \(\beta_i(g_0) = h_1^{\beta_i(g_1)}\), consequently $h_1 \in \beta_i(G) \cap S_i$. But $h_1|^u = \id$ for all $u$ with $|u| \equiv_s i$, while $\beta_i(G)|^u = \{ \id \}$ for $u \notin \omega_i(T)$ by \cref{Lemma}{lem:portrait}. Thus $h_1 = \id = b \notin \mathcal M$, which is a contradiction.
	
	We split the proof of statements $(1)$ to $(3)$ into two cases: $i = 0$ and $i \neq 0$.\\
	
	\noindent\emph{Case $i = 0$}: Since $N_0 \leq \St_{B}(1)$, statement $(3)$ is fulfilled. We have $S_0|_0 = S_{s-1}$ and $S_0|_x = \{\id\}$ for $x \in X\setminus\{0\}$. Also $\beta_0(G)|_x \leq \beta_{s-1}(G|_x)$ for $x \in X$, hence $N_0|_x \leq N_{s-1}$. Thus all sections $b|_x$ are members of $\beta_{s-1}(G) \cap N_{s-1} \subseteq \mathcal N$.
	
	The first layer sections of $b$ are given by
	\[
		b|_x = \beta_{s-1}(g_0|_x)= \prod_{j \in L_x} (h_j|_0)^{\beta_{s-1}(g_j|_x)}, \quad \text{ for } x \in X,
	\]
	where $L_x = \{j \mid 1 \leq j \leq \ell \text{ and } g_j(x) = 0\}$. The sum $\sum_{x \in X} |L_x|$ equals $\ell$. By the minimality of $\ell$, either all sections of $b$ are contained in $\beta_{s-1}([G,G])$, or there is some $x_i \in X$ such that $\ell(b|_{x_i}) = |L_{x_i}| = \ell$. In the first case, since $[G,G]$ is non-absorbing in $G$, this implies $b \in \beta_0([G,G])$, a contradiction. In the second case, $L_x = \emptyset$ for $x \neq x_i$, i.e.\ $b|_x = \id$ for $x \neq x_i$. This proves statement (2). Furthermore, if $b|_{x_i} \notin \mathcal M$, it is contained in $\beta_{s-1}([G,G])$. Since $[G,G]$ is non-absorbing over $G$, this implies $b \in \beta_0([G,G])$. Thus $b|_{x_i} \in \mathcal M$, and statement (1) is true.\\
	
	\noindent\emph{Case $i \neq 0$}: Recall that $b|_x = \beta_i(g_0)|_x = \id$ for $x \neq 0$. This is statement $(2)$ with $x_i = 0$.
	
	We consider the first layer sections of $b$. For $x \in X$ and $1 \leq j \leq \ell$,
	\begin{equation}\label{eq:fls}
		\tag{$\dagger$}h_j^{\beta_i(g_j)}|_x = \begin{cases}
			(h_j|_x)^{\beta_{i-1}(g_j)} &\text{ if }x = 0 \text{ and } h_j \in \st_{B}(0),\\
			h_j|_x \beta_{i-1}(g_j) &\text{ if }x = 0 \text{ and } h_j \notin \st_{B}(0),\\
			\beta_{i-1}(g_j^{-1}) h_j|_x &\text{ if } h_j \not\in \st_B(0) \text{ and } x = h_j^{-1}(0),\\
			h_j|_x &\text{ otherwise}.
		\end{cases}
	\end{equation}
	Since $G|^\epsilon$ acts regularly, $\st_{B}(0) = \St_{B}(1)$. We divide the long product in (\ref{eq:dec}) into segments that stabilise the first layer: Let $x \in X$, and consider the subsequence $(j^{(k)}_x)_{k \in \intseg{1}{t_x}}$ of $\intseg{1}{\ell}$ consisting of all indices $j^{(k)}_x$ such that \((\prod_{j = j^{(k)}_x}^{\ell} h_j)(x) = 0\). Clearly $\sum_{x \in X} t_x = \ell$.
	
	Set $j^{(0)}_x = 1$ and $j^{(t_x+1)}_x = \ell + 1$. Then $\prod_{j = j^{(k)}_x}^{j^{(k+1)}_x-1} h_j \in \St_{B}(1)$ for all $k \in \intseg{1}{t_x}$, and one may write
	\begin{equation}\label{eq:segments}
		\tag{$\ddagger$} b = \prod_{k = 0}^{t_x} \prod_{j = j^{(k)}_x}^{j^{(k+1)}_x-1} (h_j)^{\beta_i(g_j)}.
	\end{equation}
	
	We now make another case distinction.\\
	
	\noindent\emph{Subcase $t_x = \ell$ for some $x \in X\setminus\{0\}$:} We will prove that this case can not occur. The equation $t_x = \ell$ implies $h_\ell(x) = 0$ and $h_j \in \St_{B}(1)$ for all $j \in \intseg{1}{\ell-1}$. We may assume $g_\ell = \id$, by passing to a conjugate if necessary. Looking at the second and fourth case of (\ref{eq:fls}), we obtain
	\[
		\beta_{i-1}(g_0) = b|_0 = \prod_{j = 1}^{\ell-1} (h_j|_{h_{\ell}(0)}) \cdot h_{\ell}|_0 \in N_{i-1}.
	\]
	Thus $\beta_{i-1}(g_0)$ is an element of $\mathcal N$ of restricted syllable length at most $1$, hence trivial. Consequently $g_0$ and $b$ are trivial, a contradiction.\\

	\noindent\emph{Subcase $t_0 = \ell$:} This implies $h_j \in \St_{B}(1)$ for all $j \in \intseg{1}{\ell}$, and statement (3) holds. By the first case of (\ref{eq:fls})
	\[
		b|_0 = \prod_{j = 1}^{\ell(b)} {h_j|_0}^{\beta_{i-1}(g_j)} \in N_{i-1} \cap \beta_{i-1}(G),
	\]
	which is of restricted syllable length at most $\ell$. As we previously argued in the case $i = 0$, we have $b|_0 \notin \beta_{i-1}([G,G])$ and consequently statement (1) holds, since otherwise $b \in \beta_i([G,G])$ because $[G,G]$ is non-absorbing over $G$.\\
	
	\noindent\emph{Subcase $t_x < \ell$ for all $x \in X$:} We shall prove that this case can not occur. Combining (\ref{eq:segments}) with (\ref{eq:fls}) for $x \in X$ we calculate
	\[
		b|_x = \prod_{k = 0}^{t_x-1} \left(\left(\Pi_{j = j^{(k)}_x}^{j^{(k+1)}_x-1} (h_j)^{\beta_i(g_j)}\right)|_0\right) \left(\Pi_{j = j^{({t_x})}_x}^{\ell(b)} (h_j)^{\beta_i(g_j)}\right)|_x
	\]
	and for $k \in \intseg{1}{t_x-1}$
	\begin{align*}
		\prod_{j = j^{(k)}_x}^{j^{(k+1)}_x-1} (h_j)^{\beta_i(g_j)}|_0 &= \beta_{i-1}(g_{j^{(k)}_x}^{-1})(\prod_{j = j^{(k)}_x}^{j^{(k+1)}_x-1} h_j|_{\prod_{i = j + 1}^{j^{(k+1)}_x-1} h_i(0)}) \beta_{i-1}(g_{j^{(k+1)}_x-1})\\
		&= \beta_{i-1}(g_{j^{(k)}_x}^{-1}g_{j^{(k+1)}_x-1}) (\prod_{j = j^{(k)}_x}^{j^{(k+1)}_x-1} h_j|_{\prod_{i = j + 1}^{j^{(k+1)}_x-1} h_i(0)})^{\beta_{i-1}(g_{j^{(k+1)}_x-1})}.
	\end{align*}
	Consequently, every segment $\prod_{j = j^{(k)}_x}^{j^{(k+1)}_x-1} (h_j)^{\beta_i(g_j)}$ of $b$ contributes at most one syllable of $N_{i-1}$ and a member of $\beta_{i-1}(G)$ to $b|_x$. We obtain
	\begin{align*}	   
	   b|_x \equiv_{N_{i-1}} \begin{cases}
	   \beta_{i-1} \left(g_1^{-1} \prod_{k = 1}^{t_x} \left(g_{j^{(k)}_x-1}g_{j^{(k)}_x}^{-1}\right) g_\ell \right) &\text{ if } x = 0,\\
	   \beta_{i-1} \left(\phantom{g_1^{-1}} \prod_{k = 1}^{t_x} \left(g_{j^{(k)}_x-1}g_{j^{(k)}_x}^{-1}\right)\phantom{g_\ell}\right) &\text{ otherwise}.
	   \end{cases}
	\end{align*}
	Write $b|_x = \beta_{i-1}(f_x) n_x$ with $n_x \in N_i$ and $f_x$ equal to the corresponding product in $G$ in the last equation. Since the subsequences form a partition, every $\beta_{i-1}(g_{j^{(k)}_x})$ and its inverse appear in precisely one section of $b$, and we have
	\[
		\prod_{x \in X} b|_x \equiv_{N_{i-1}} \prod_{x \in X} \beta_{i-1}(f_x) \equiv_{\beta_{i-1}([G,G])} \prod_{j = 1}^{\ell}{\beta_{i-1}(g_jg_j^{-1})} = 1.
	\]
	Now we look at $n_x$. Since every segment $\prod_{j = j^{(k)}_x}^{j^{(k+1)}_x-1} (h_j)^{\beta_i(g_j)}$ contributes at most one syllable, and $h_j \notin \St_{B}(1)$ for some $j \in \intseg{1}{\ell}$, we have $\ell(n_x) \leq t_x < \ell$. Also \(\beta_{i-1}(f_x)n_x = b|_x = \id\) for $x \neq 0$, hence $n_x = \beta_{i-1}(f_x^{-1}) \in \mathcal N$. By minimality, $f_x \in [G,G]$. Then also \(f_0 \equiv_{[G,G]} \prod_{x \in X} f_x \equiv_{[G,G]} \id\), and $\beta_{i-1}(f_0^{-1}g_0) = \beta_{i-1}(f_0^{-1})b|_0 = n_0 \in \mathcal N$. Again, by minimality, $f_0^{-1}g_0 \in [G,G]$, thus $g_0 \in [G,G]$, a contradiction.
	
	This completes the proof.
\end{proof}

\begin{eg}
	Let $\bp_s(\mathcal O_m^d)$ be a generalised Basilica group (cf.\ \cref{Definition}{def:gen_basilica}). Since $\mathcal O_m^d$ is free abelian and self-similar, and $\mathcal O_m^d|^\epsilon$ is cyclic of order $m$, by \cref{Proposition}{prop:split-intersection_in_commutator}, the group $\mathcal O_m^d$ is $s$-split.
\end{eg}

\begin{qst}
	Motivated by the small gap between \cref{Proposition}{prop:split-intersection_in_commutator} and \cref{Proposition}{prop:non-abelian not split} we ask:
	
	\textit{Is every abelian group $G \leq \Aut(T)$ acting spherically transitive $s$-split for all $s > 1$?}
\end{qst}

\begin{coro}\label{cor:abelisation}
	Let $G \leq \Aut(T)$ be a self-similar $s$-split group. Then the abelianisation $\bp_s(G)$ is
	 \[
	 	\bp_s(G)^{\ab} \cong G^s.
	 \]
\end{coro}

\begin{proof}
	Consider the normal subgroup $H \defeq \LD [\beta_i (G), \beta_j (G)] \mid i, j \in \intseg{0}{s-1}, i \neq j \RD^{\bp_s(G)}$ and observe that $H \leq N_i$ for all $i \in \intseg{0}{s-1}$. We obtain an epimorphism $G^s \to \bp_s(G)/H$, mapping the $i$-th component of $G^s$ to $\beta_i(G) (H)$, for $i \in \intseg{0}{s-1}$. This map is also injective. Let $\prod_{i \in \intseg{0}{s-1}}\beta_i(g_i) \equiv_H \prod_{i \in \intseg{0}{s-1}}\beta_i(h_i)$ for some $g_i, h_i \in G$. Then for all $x \in X$
	\[
		\beta_x(g_xh_x^{-1}) \equiv_H \prod_{i \in \intseg{0}{s-1}\setminus\{x\}} \beta_i(g_i^{-1}h_i) \in N_x
	\]
	and $\beta_x(g_xh_x^{-1}) \in N_x$. Since $G$ is $s$-split, this implies $g_x = h_x$. Thus $\bp_s(G)/H \cong G^s$. But from \cref{Proposition}{prop:non-abelian not split} $G$ is abelian and consequently $H = [\bp_s(G), \bp_s(G)]$.
\end{proof}

\begin{prop}\label{prop:torfree}
	Let $G \leq \Aut(T)$ be a torsion-free self-similar group such that the quotient $G/K$ with $K = \beta_{0}^{-1}(\beta_0(G) \cap N_0)$ is again torsion-free. Then $\bp_s (G)$ is torsion-free.
\end{prop}

\begin{proof}
	Let $b \in \bp_s(G)$ be a torsion element. Since $G/K$ is torsion-free, we obtain $b \in \ker\phi_0 = N_0 \leq \St_{\bp_s(G)}(1)$. Thus the first layer sections of $b$ are again torsion elements of $\bp_s(G)$, because $\bp_s(G)$ is self-similar by \cref{Lemma}{lem:self-sim}. Hence an iteration of the argument yields $b = \id$.
\end{proof}

\begin{qst}
	On the other end of the spectrum, the group ${\bp_2(\mathcal G)}$ (cf.\ \cref{Example}{eg:grigorchuks}) is periodic as is $\mathcal G$, which can be proven analogous to \cite{BGS03}*{Theorem 6.1}, and the second Basilica groups of the periodic Gupta-Sidki-$p$-groups (cf.\ \cref{Example}{eg:gupta_sidki}) are periodic by \cite{Pet21}. Motivated by this observation we ask:
	
	\emph{Is there a periodic group $G \leq \Aut(T)$ acting spherically transitive such that $\bp_s(G)$ is not periodic for some $s \in \Nplus$?}
\end{qst}

\cref{Proposition}{prop:torfree} and \cref{Corollary}{cor:abelisation} prove \cref{Theorem}{thm:split_transfer}.

\subsection{Layer Stabilisers}\label{sub:stabilisers_and_dimension}

For an $s$-split group $G \leq \Aut(T)$ the $s$-th Basilica decomposes as $\bp_s(G) = N_{i} \rtimes \beta_i(G)$. Recall from \cref{Definition}{def:s-split} that $\phi_i$ denotes the map to $\bp_s(G)/N_i$, identified with the quotient $G/K_i$, such that $\phi_i(n \beta_i(g)) = gK_i$ for all $g \in G, n \in N_i$.

\begin{lemma}\label{lem:stab}
	Let $G \leq \Aut(T)$ be a strongly fractal group and let $B = \bp_s(G)$. Let $b_0, \dots, b_{m-1} \in B$. Then $\psi_1^{-1}(b_0, \dots, b_{m-1})$ is an element of $\St_{B}(1)$ if and only if there is an element $g \in \St_G(1)$ such that for all $x \in X$
	$$
		\phi_{s-1} (b_x) = g|_x K_{s-1}.
	$$
\end{lemma}

\begin{proof}
	If there is some element $g \in \St_G(1)$ of the required form, clearly
	\[
		\beta_0(g) \equiv_{\psi_1^{-1}(N_{s-1}^m)} (b_0, \dots, b_{m-1}).
	\]
	Now we claim that $\psi_1(N_0) \geq N_{s-1}^m$. Let
	\[
		b = \prod_{j = 0}^{\ell-1} h_j^{\beta_{s-1}(g_j)} \in N_{s-1},
	\]
	with $h_j \in S_{s-1}$. Then there are elements $\hat h_j = (h_j, \id, \dots, \id) \in S_0$ by the definition of $S_{s-1}$. Furthermore, since $G$ is strongly fractal, there are elements $\hat g_j \in \St_G(1)$ such that $\beta_0(\hat g_j)|_{0} = \beta_{s-1}(g_j)$, yielding
	\[
		\prod_{j = 0}^{\ell-1} \hat h_j^{\beta_0(\hat g_j)} = (b, \id, \dots, \id).
	\]
	Since $G$ acts spherically transitively, the claim follows by \cref{Lemma}{lem:fractal_branching}. Thus there is an element in $N_0\beta_0(g) \leq \St_B(1)$ with sections $(b_0, \dots, b_{m-1})$.
	
	Let now $b = \psi_1^{-1}(b_0, \dots, b_{m-1}) \in \St_{B}(1)$. Then $b$ decomposes as a product $n\beta_0(g)$ with $n \in N_0$ and $g \in \St_G(1)$. This implies, for any $x \in X$,
	\[
		\phi_{s-1}(b_x) = \phi_{s-1}((n\beta_0(g))|_x) = \phi_{s-1}(\beta_{s-1}(g|_x)) = g|_x K_{s-1}.\qedhere
	\]
\end{proof}

\begin{lemma}\label{lem:geomcont}
	Let $G$ be fractal and self-similar and let $B = \bp_s(G)$. Let $n \in \N$.
	\begin{enumerate}
		\item[(i)] $\psi_1(\beta_{i}(\St_G(n))^B) = (\beta_{i-1}(\St_G(n))^B)^m$ for all $i \neq 0$.
	\end{enumerate}
	Assuming further that $G$ is very strongly fractal,
	\begin{enumerate}
		\item[(ii)] $\psi_1([\beta_0(\St_G(n+1)), N_0]^B) = ([\beta_{s-1} (\St_G(n)), N_{s-1}]^B)^m$.
	\end{enumerate}
\end{lemma}

\begin{proof}
	(i) The inclusion $\psi_1(\beta_{i}(\St_G(n))^B) \leq (\beta_{i-1}(\St_G(n))^B)^m$ is obvious. We prove the other direction. Let $g \in \St_G(n)$ and $b \in B$. Since $B$ is fractal by \cref{Lemma}{lem:fractality}, there is an element $c \in \st_B(0)$ such that $c|_0 = b$. Now
	\[
		(\beta_{i}(g))^{c} = (\beta_{i-1}(g), \id, \dots, \id)^c = (\beta_{i-1}(g)^b, \id, \dots, \id),
	\]
	yielding statement (i), by \cref{Lemma}{lem:fractal_branching}.

	(ii) The inclusion $\psi_1([\beta_0(\St_G(n+1)), N_0]^B) \leq ([\beta_{s-1} (\St_G(n)), N_{s-1}]^B)^m$ follows directly from $N_0|_x \leq N_{s-1}$ and $\beta_{0}(\St_G(n + 1))|_x \leq \beta_{s-1}(\St_G(n))$, where $x \in X$. Thanks to \cref{Lemma}{lem:fractal_branching}, for the other inclusion it is enough to prove that
	$
		([\beta_{s-1}(g), k], \id, \dots, \id)$ is contained in $\psi_1([\beta_0(\St_G(n+1)), N_0]^B)
	$
	for all $g \in \St_G(n)$ and $k \in N_{s-1}$.
	Let 
	\[
		k = \prod_{j = 0}^\ell (\beta_{i_j} (k_{j}))^{\beta_{s-1}(k_j')} \in N_{s-1}.
	\]
	Since $G$ is strong fractal there are elements $t_j \in \St_G(1)$ such that $t_j|_0 = k_j'$. Furthermore, since $G$ is very strongly fractal there is an element $h \in \St_G(n+1)$ such that $h|_0 = g$. Then
	\begin{align*}
		[\beta_0 (h), \prod_{j = 0}^\ell (\beta_{i_j+1}(k_j))^{\beta_{0}(t_j)}] \in [\beta_0(\St_G(n+1)), N_0]^B
	\end{align*}
	and
	\begin{align*}
		[\beta_0 (h), \prod_{j = 0}^\ell (\beta_{i_j+1}(k_j))^{\beta_{0}(t_j)}]|_x &= [(\beta_0 (h))|_x, \prod_{j = 0}^\ell ((\beta_{i_j+1}(k_j))|_x)^{(\beta_{0}(t_j))|_x}]\\
		&= \begin{cases}
			[\beta_{s-1} (g), k] &\text{ if }x = 0,\\
			[\beta_{s-1} (h|_x), \prod_{j = 0}^\ell \id^{\beta_{s-1}(t_j|_x)}] = \id &\text{ otherwise.}
		\end{cases}\qedhere
	\end{align*}
\end{proof}

\begin{proof}[Proof of \cref{Theorem}{thm:stab}]
	Let $B = \bp_s(G)$. For any $n \in \N$, write $n = s q + r$ with $q \geq 0$ and $r \in \intseg{0}{s-1}$. We have to prove
	\[
		\St_B(n) = \LD \beta_i(\St_G(q+1)), \beta_{j}(\St_G(q)) \mid 0 \leq i < r \leq j < s\RD^B.
	\] 
	For convenience, we will denote the right-hand side of this equation by $H_n$. It is clear that $H_n \leq \St_B(n)$ for all $n \in \N$. It remains to establish the other inclusion. For $n = 0$ the statement is clearly true, so we proceed by induction and assume that the statement is true for some fixed $n = s q + r$ with $q \geq 0$ and $r \in \intseg{0}{s-1}$.
	Define
	\begin{align*}
		J \defeq \LD \beta_i(\St_G(q + 1)), \beta_j(\St_G(q)), [\beta_{s-1}(\St_G(q)), N_{s-1}]^B \mid 0 \leq i \leq r - 1 < j < s-1 \RD^B,
	\end{align*}
	and observe that by \cref{Lemma}{lem:geomcont} we find $J^m \leq \psi_1(H_{n+1})$, which yields
	\[
		(\St_B(n))^m/\psi_1(H_{n+1}) = (\beta_{s-1}(\St_G(q)))^m\psi_1(H_{n+1})/\psi_1(H_{n+1}).
	\]
	Hence for every $g \in \St_B(n+1)$, there are elements $t_0,\dots,t_{m-1} \in \St_G(q)$ such that
	\[
		\psi_1(g) \equiv_{\psi_1(H_{n+1})} (\beta_{s-1}(t_0), \dots, \beta_{s-1}(t_{m-1})).
	\]
	Since $\phi_{s-1}\beta_{s-1}(t_x)=t_x K_{S-1}$ for all $x \in X$, $g \in \St_B(1)$ and $H_{n+1} \leq \St_B(1)$, by \cref{Lemma}{lem:stab} there are elements $k_0, \dots, k_{m-1} \in K_{s-1}$ and $h \in \St_G(1)$ such that
	\[
		\psi_1^{-1}(h|_0k_0, \dots, h|_{m-1}k_{m-1}) =
		\psi_1^{-1}(t_0, \dots, t_{m-1}).
	\]
	Define $\widetilde h = h\psi_1^{-1}(k_0, \dots, k_{m-1})$. Now $G$ is weakly regular branch over $K_{s-1}$, hence $\psi_1^{-1}(K_{s-1}^m) \leq \St_{K_{s-1}}(1)$, and consequently $\widetilde h \in \St_G(1)$.
	But $\widetilde h|_x = t_x \in \St_G(q)$ for $x \in X$, whence $\widetilde h \in \St_{G}(q+1)$ and
	\[
		(\beta_{s-1}(t_0), \dots, \beta_{s-1}(t_{m-1})) = \psi_1(\beta_0(\widetilde h)) \in \psi_1(\beta_0(\St_G(q+1))) \leq \psi_1(H_{n+1}),
	\]
	implying $g \in H_{n+1}$. This completes the proof.
\end{proof}

\subsection{Hausdorff Dimension}\label{sub:hausdorff_dimension}

We remind the reader that $\Gamma$ is the subgroup of $\Aut(T)$ consisting of all automorphisms whose labels are elements of $\langle \sigma \rangle$, with $\sigma$ being a fixed $m$-cycle in $\Sym(X)$.

\begin{defn}
	Let $G \leq \Gamma$. The \emph{Hausdorff dimension of} $G$ \emph{relative to} $\Gamma$ is defined by
	\begin{align*}
		\dimH G \defeq \liminf_{n\to\infty}\frac{\log_{m}|G/\St_G(n)|}{\log_{m}|\Gamma/\St_{\Gamma}(n)|}
		= (m-1)\liminf_{n\to\infty}\frac{\log_{m}|G/\St_G(n)|}{m^n}.
	\end{align*}
	This relates to the usual definition of Hausdorff dimension over arbitrary spaces by taking the closure, i.e.\ using this definition, the group $G$ has the same Hausdorff dimension as its closure $\overline G$ in $\Gamma$, cf.\ \cite{BS97}.
	We drop the base $m$ in $\log_m$ from now on. Denote the quotient \(\St_G(n)/\St_G(n+1)\) by $L_{G}(n)$.
	The integer series (for $n > 0$) obtained by
	\[
		o_{G}(n) \defeq \log(|L_G(n-1)|^m)-\log|L_G(n)|
	\]
	is called the \emph{series of obstructions of} $G$. We set $o_G(0) = -1$ for convenience.
\end{defn}

The series of obstructions of a group $G$ determines its Hausdorff dimension, precisely how we will see in \cref{Lemma}{lem:obstruction_to_hdim}. Nevertheless, one might wonder why it is necessary to define this seemingly impractial invariant. We will demonstrate in \cref{Proposition}{prop:obstruction} that it is (to some degree) preserved under $G \mapsto \bp_s(G)$. Furthermore, many well-studied subgroups of $\Gamma$ have a well-behaved series of obstructions. For example, it is easy to see that $\Gamma$ itself has 
\begin{align*}
	o_\Gamma(n) &= \log |\wr^{n} \mathrm C_m/\wr^{n-1} \mathrm C_m|^m - \log |\wr^{n+1} \mathrm C_m/\wr^{n} \mathrm C_m| \\
	&= m \log m^{m^{n}} - \log m^{m^{n+1}} = 0,
\end{align*}
for $n \in \Nplus$, where $\wr^n A$ is the $n$-times iterated wreath product of $A$, with the convention that $\wr^0 A$ is the trivial group. On the other hand, since the layer stabiliser of $\mathcal O_m^d$ are the subgroups generated by $\LD \pi_0(a)^{m^{k+1}},\dots,\pi_{l-1}(a)^{m^{k+1}},\pi_l(a)^{m^k},\dots,\pi_{d-1}(a)^{m^k}\RD$, the quotients $L_{\mathcal O_m^d}(n)$ are all cyclic of order $m$, and
\begin{align*}
	o_{\mathcal O_m^d}(n)
	& = m - 1.
\end{align*}
A Gupta--Sidki $p$-group $G$ has precisely two terms unequal to $0$, a consequence of $\St_G(n) = \St_G(n-1)^p$ for $n \geq 3$, cf.\ \cite{FAZR14}. Similarly, the series of obstructions of the Grigorchuk group has only one non-zero term.

\begin{lemma}\label{lem:obstruction_to_hdim}
	Let $G \leq \Gamma$ act spherically transitive. Then
	\[
		\dimH G = 1 - \limsup_{n \to \infty} \sum_{i = 1}^{n} (m^{-i} - m^{-(n+1)}) o_G(i).
	\]
\end{lemma}

\begin{proof}
	By definition $\log |L_G(0)| = 1$ and
	\(
		\log |L_G(n)| = m \log|L_{G}(n-1)| - o_G(n)
	\)
	for $n \geq 1$. An inductive argument yields
	\begin{align*}
		\log|G/\St_G(n+1)|
		&= \log|G/\St_G(n)| - \sum_{k = 0}^{n} m^{n-k} o_G(k)
		&= - \sum_{k = 0}^n \frac{m^{k+1}-1}{m-1} o_G(n - k).
	\end{align*}
	This gives
	\begin{align*}
		\liminf_{n\to\infty}\frac{(m-1)}{m^{n + 1}}\log\left|\frac{G}{\St_G(n + 1)}\right| &= - \limsup_{n \to \infty} \sum_{i = 0}^{n} (m^{i - n} - m^{-(n+1)}) o_G(n - i)\\
		&= 1 - \limsup_{n \to \infty} \sum_{i = 1}^{n} (m^{-i} - m^{-(n+1)}) o_G(i).\qedhere
	\end{align*}
\end{proof}

\begin{lemma}\label{lem:self-sim_obstructions}
	Let $G \leq \Gamma$ be self-similar. Then for all $n > 0$
	\[
		o_G(n) = \log [\St_G(n-1)^m : \psi_1(\St_G(n))] - \log [\St_G(n)^m : \psi_1(\St_G(n+1))].
	\]
\end{lemma}

\begin{proof}
	We have, for $n > 0$,
	\begin{align*}
		\left|\frac{\St_G(n-1)^m}{\psi_1(\St_G(n))}\right|
		&= \frac{|\St_G(n-1)^m/\psi_1(\St_G(n+1))|}{|L_{G}(n)|}\\
		&= \frac{|L_{G}(n-1)|^m}{|L_{G}(n)|}\left|\frac{\St_G(n)^m}{\psi_1(\St_G(n+1))}\right|,
	\end{align*} hence
	\[
		o_{G}(n) = \log[\St_G(n-1)^m:\psi_1(\St_G(n))] - \log[\St_G(n)^m:\psi_1(\St_G(n+1))].\qedhere
	\]
\end{proof}

\begin{lemma}\label{lem:betterGeomcont}
	Let $G$ be very strongly fractal, self-similar and weakly regular branch over the splitting kernel $K_{s-1}$. Then for all $\ell, n \in \Nplus$
	\begin{align*}
		&\psi_1(\beta_0(\St_G(\ell + 1)) \cap [\beta_0(\St_G(n+1)), N_0]^B)\\
		= \ &(\beta_{s-1}(\St_G(\ell)) \cap [\beta_{s-1}(\St_G(n)), N_{s-1}]^B)^m.
	\end{align*}
\end{lemma}

\begin{proof}
	The left-hand set is clearly contained in the right-hand set. We prove the other inclusion. Let $(b_0, \dots, b_{m-1}) \in (\beta_{s-1}(\St_G(\ell)) \cap [\beta_{s-1}(\St_G(n)), N_{s-1}]^B)^m$. By \cref[(ii)]{Lemma}{lem:geomcont} there exists $b \in [\beta_0(\St_G(n+1)), N_0]^B \leq \St_B(1)$ such that $\psi_1(b) = (b_0, \dots, b_{m-1})$. It remains to prove that $b \in \beta_0(\St_G(\ell + 1))$.
	
	Since the set $\beta_{s-1}(\St_G(\ell)) \cap [\beta_{s-1}(\St_G(n)), N_{s-1}]^B$ is contained in $\beta_{s-1}(\St_{K_{s-1}}(1))$ and since $G$ weakly regular branch over $K_{s-1}$, there is an element $g \in K_{s-1}$ such that
	\[
		\psi_1(g) = (\beta_{s-1}^{-1}(b_0), \dots, \beta_{s-1}^{-1}(b_{m-1})) \in \St_G(\ell)^m.
	\]
	Consequently, $\psi_1(\beta_0(g)) = (b_0, \dots, b_{m-1}) = \psi_1(b)$, and $b = \beta_0(g)$ is a member of $\psi_1(\beta_0(\St_G(\ell + 1)) \cap [\beta_0(\St_G(n+1)), N_0]^B)$.
\end{proof}

\begin{prop}\label{prop:obstruction}
	Let $G \leq \Gamma$ be very strongly fractal, self-similar and weakly regular branch over the splitting kernel $K_{s-1}$. Then the series of obstructions for $B = \bp_s(G)$ fulfills
	\begin{align*}
		o_{B}(n) = \begin{cases}
			0 &\text{ if }n \not\equiv_s 0,\\
			o_{G}(\frac{n}{s}) &\text{ otherwise.}
		\end{cases}
	\end{align*}
\end{prop}

\begin{proof}
	Consider first the case $n \equiv_s k \neq 0$. By \cref{Theorem}{thm:stab} the quotient
	\(
		L_{B}(n)
	\)
	is normally generated in $B$ by images of elements of $\beta_{k}(\St_G(\lfloor n/s \rfloor))$. Similarly the images of $\beta_{k-1}(\St_G(\lfloor n/s \rfloor))$ are the normal generators of $L_{B}(n-1)$. Thus \cref[(i)]{Lemma}{lem:geomcont} shows that $o_B(n) = 0$.	
	
	Now consider the case $n = qs$. To shorten the notation, we abbreviate
	\begin{align*}
		R_q &\defeq \beta_0(\St_G(q)) \text{ for } q \in \N \text{ and}\\
		T_q &\defeq \beta_{s-1}(\St_G(q)) \text{ for } q \in \N.
	\end{align*}
	Define the normal subgroups
	\begin{align*}
		U &= \LD \St_B(n+1) \cup [R_q, N_{0}]^B\RD \trianglelefteq B \quad\text{and}\\
		V &= \LD \St_B(n) \cup [T_{q-1}, N_{s-1}]^B\RD \trianglelefteq B.
	\end{align*}
	Using \cref{Theorem}{thm:stab}, we see that $U$ and $V$, respectively, are normally generated by the sets
	\begin{align*}
		R_{q+1} \cup \bigcup_{i = 1}^{s-1} \left( \beta_{i}(\St_G(q)) \right) \cup [R_q, N_{0}] \quad \text{and} \quad
		T_q \cup \bigcup_{i = 0}^{s-2} \left( \beta_{i}(\St_G(q)) \right) \cup [T_{q-1}, N_{s-1}].
	\end{align*}
	Let $g \in \St_G(q+1)$ and $b \in B$. We write $b = \beta_0(g_b)n_b$ for $g_b \in G$ and $n_b \in N_0$. Then
	\[
		\beta_0(g)^b
		= \beta_0(g^{g_b})^{n_b} = \beta_0(g^{g_b})[\beta_0(g^{g_b}), n_b] \in R_{q+1} [R_{q+1}, N_0].
	\]
	Consequently, we drop the conjugates of $R_{q+1}$ in our generating set for $U$, and write
	\[
		U = \LD R_{q+1} \cup \bigcup_{i = 1}^{s-1} \left( \beta_{i}(\St_G(q))^B \right) \cup [R_q, N_{0}]^B \RD.
	\]
	Similarly, the subgroup $V$ is generated by
	\[
		T_q \cup \bigcup_{i = 0}^{s-2} \left( \beta_{i}(\St_G(q))^B \right) \cup [T_{q-1}, N_{s-1}]^B.
	\]
	Using \cref{Theorem}{thm:stab}, it is now easy to see that
	\begin{align*}
		\St_B(n)/U
		&\cong R_q/(R_q \cap U).
	\end{align*}
	Since $\beta_{i}(\St_G(q)) \leq \St_B(n+1)$ for $i \neq 0$, we see that the intersection
	\[
		\LD \beta_1(\St_G(q)) \cup \dots \cup \beta_{s-2}(\St_G(q))\cup T_q \RD^B \cap R_q \leq R_{q+1}
	\]
	is contained in $R_{q+1}$. We conclude
	\[
		R_q \cap U = R_q \cap R_{q+1} [R_q, N_{0}]^B.
	\]
	Now
	\begin{align*}
		R_q \cap R_{q+1} [R_q, N_{0}]^B = R_{q+1}(R_q \cap [R_q, N_{0}]^B)
	\end{align*}
	and
	\begin{align*}
		[R_q \cap R_{q+1} [R_q, N_{0}]^B : R_{q+1}]
		 = [R_q \cap [R_q, N_{0}]^B: R_{q+1} \cap [R_q, N_{0}]^B].
	\end{align*}
	Consequently, the order of $\St_B(n)/U$ equals
	\begin{align*}
		|L_G(q)| \cdot [R_q \cap [R_q, N_{0}]^B: R_{q+1} \cap [R_q, N_{0}]^B]^{-1}.
	\end{align*}
	A similar computation shows that the order of $\St_B(n-1)/V$ is
	\[
		|L_G(q-1)| \cdot [T_{q-1} \cap [T_{q-1}, N_{s-1}]^B: T_q \cap [T_{q-1}, N_{s-1}]^B]^{-1}.
	\]
	We now apply \cref{Lemma}{lem:betterGeomcont} in the cases $\ell = n = q - 1$ and $\ell = n + 1 = q$, i.e.\ we have
	\begin{align*}
		\psi_1(R_q \cap [R_q, N_{0}]^B) &= (T_{q-1} \cap [T_{q-1}, N_{s-1}]^B)^m \text{ and}\\
		\psi_1(R_{q+1} \cap [R_q, N_{0}]^B) &= (T_q \cap [T_{q-1}, N_{s-1}]^B)^m.
	\end{align*}
	We see that the second factor in the formula for the order of $\St_B(n)/U$ is the $m$-th power of the corresponding factor for $\St_B(n-1)/V$, and obtain
	\[
		\frac{|\St_B(n-1)/V|^m}{|\St_B(n)/U|} = \frac{|L_G(q-1)|^m}{|L_G(q)|} = m^{o_G(q)}.
	\]
	Now we compare $V^m$ and $\psi_1(U)$. By \cref[(i)]{Lemma}{lem:geomcont} and \hyperref[lem:geomcont]{(ii)}, $\psi_1(U)$ is generated by
	\begin{align*}
		\psi_1(R_{q+1}) \cup \bigcup_{i = 0}^{s-2} \left( ((\beta_i(\St_G(q)))^B)^m  \right) \cup ([T_{q-1}, N_{s-1}]^B)^m.
	\end{align*}
	We define yet another subgroup
	\[
		W = \LD \ \bigcup_{i = 0}^{s-2} \left( ((\beta_i(\St_G(q)))^B)^m  \right) \cup ([T_{q-1}, N_{s-1}]^B)^m \rangle \leq \psi_1(U) \leq B^m.
	\]
	Evidently $W \trianglelefteq B^m$, $W \leq N_{s-1}^m$, and $W \trianglelefteq \psi_1(U) \leq V^m$. We have
	\begin{align*}
		\psi_1(U)/W &\cong \psi_1(R_{q+1})/(\psi_1(R_{q+1}) \cap W) \text{ and}\\
		V^m/W &\cong  T_q^m/(T_q^m \cap W).
	\end{align*}
	The two divisors are equal: Clearly $\psi_1(R_{q+1}) \cap W$ is contained in $T_q^m \cap W$. Let $(\beta_{s-1}(g_0), \dots, \beta_{s-1}(g_{m-1})) \in T_q^m \cap W \leq (T_q \cap N_{s-1})^m$. Since $T_q \cap N_{s-1} \leq \beta_{s-1}(K_{s-1})$, the elements $g_0, \dots, g_{m-1}$ are members of $K_{s-1}\cap\St_G(q)$. Now since $G$ is weakly regular branch over $K_{s-1}$, there is an element $k \in K_{s-1} \cap \St_G(q+1)$ such that
	\(
		\psi_1(k) = (g_0, \dots, g_{m-1}),
	\)
	and consequently $\beta_0(k) \in R_{q+1}$ fulfills
	\[
		\psi_1(\beta_0(k)) = (\beta_{s-1}(g_0), \dots, \beta_{s-1}(g_{m-1})) \in \psi_1(R_{q+1}) \cap W.
	\]
	We compute
	\begin{align*}
		[V^m:\psi_1(U)] &= [V^m/W : \psi_1(U)/W]\\
		&= [T_q^m: \psi_1(R_{q+1})]\\
		&= [(\beta_{s-1} \times \dots \times \beta_{s-1})(\St_G(q)^m) :
		(\beta_{s-1} \times \dots \times \beta_{s-1})(\psi_1(\St_G(q+1)))]\\
		&= [\St_G(q)^m : \psi_1(\St_G(q+1))].
	\end{align*}
	This implies
	\begin{align*}
		[\St_B(n-1)^m : \psi_1(\St_B(n))] &= [{\St_B(n-1)^m/\psi_1(U)}:{\psi_1(\St_B(n))/\psi_1(U)}]\\
		&= \frac{[\St_B(n-1)^m : V^m][V^m : \psi_1(U)]}{[\psi_1(\St_B(n)) : \psi_1(U)]}\\
		&= \frac{|L_G(q-1)|^m}{|L_G(q)|} \cdot [\St_G(q)^m : \psi_1(\St_G(q+1))].
	\end{align*}
	Since $o_{B}(k) = 0$ for $k \not\equiv_s 0$, by \cref{Lemma}{lem:self-sim_obstructions},
	$$
		\log[\St_B(n)^m : \psi_1(\St_B(n+1))] = \log[\St_B(n + s - 1)^m : \psi_1(\St_B(n + s)],
	$$
	hence
	\begin{align*}
		o_{B}(n) 
		=& \log[\St_B(n-1)^m : \psi_1(\St_B(n))] - \log[\St_B(n+s-1)^m : \psi_1(\St_B(n+s))]\\
		=& o_{G}(q) + \log\left|\frac{\St_G(q)^m}{\psi_1(\St_G(q+1))}\right|
		- o_{G}(q+1) - \log\left|\frac{\St_G(q+1)^m}{\psi_1(\St_G(q+2))}\right|\\
		=& o_{G}(q) - o_{G}(q+1) + o_{G}(q+1)\\
		=& o_{G}(q).\qedhere
	\end{align*}
\end{proof}

\begin{proof}[Proof of \cref{Corollary}{cor:hausdorff_bounded}]
	By \cref{Lemma}{lem:obstruction_to_hdim} and \cref{Proposition}{prop:obstruction}
	\begin{align*}
		\dimH G &= 1 - \limsup_{n \to \infty} \sum_{i = 1}^n (m^{-i} - m^{-(n + 1)})o_G(i)\quad\text{and}\\
		\dimH \bp_s(G) &= 1 - \limsup_{n \to \infty} \sum_{i = 1}^n (m^{-i} - m^{-(n + 1)})o_{\bp_s(G)}(i)\\
		&= 1 - \limsup_{n \to \infty} \sum_{i = 1}^n (m^{-si} - m^{-(sn + 1)})o_{G}(i).
	\end{align*}
	We prove $m^{-i} - m^{-(n+1)} > m^{-si} - m^{-(sn+1)}$, equivalently $m^{sn+1 - i}+1 > m^{s(n-i)+1} + m^{(s-1)n}$. This is a consequence of $sn + 1 - i - (s(n-i) + 1) = (s-1)i \geq 1$ and $sn+1-i - (s-1)n = n-i + 1\geq 1$, with equality precisely when $i=1, s=2$, resp.\ $n=i$. Therefore at least one of the differences is greater than $1$, and the limit of $\sum_{i = 1}^n (m^{-si} - m^{-(sn + 1)})o_{G}(i)$ is strictly greater than the limit of $\sum_{i = 1}^n (m^{-i} - m^{-(n + 1)})o_{G}(i)$. The statement follows.
\end{proof}

\begin{eg}\label{eg:ggs-hausdorff}
	Let $G \leq \Aut(T_p)$, $p$ a prime, be a $\mathsf{GGS}$-group defined by the triple $(\mathsf C_p, \mathsf C_p, \omega)$, cf. \cref{Definition}{def:spinal}, where $\mathsf C_p$ denotes the cyclic group of order $p$ acting regularly on $X$. To be a $\mathsf{GGS}$-group means $\omega_i = \omega_j$ for $i, j \in \N$, thus we write $\omega$ for $\omega_1$. This is a $(p-1)$-tuple of endomorphisms of $\mathsf C_p$. Every such endomorphism is a power map, hence we may identify $\omega$ with an element $(e_1, \dots, e_{p-1})$ of $\mathbb F_p^{p-1}$. Assume that
	\begin{equation}\label{eq:torsion}
		\tag{$\star$}e_1 + \dots + e_{p-1} \equiv_p 0
	\end{equation}
	and that 
	there is some $i \in \intseg{1}{p-1}$
	\begin{equation}\label{eq:regular}
		\tag{$\diamond$}e_i \neq e_{p-i}.
	\end{equation}
	In \cite{FAZR14} the order of the congruence quotients $G/\St_G(n)$ is explicitly calculated in terms of the rank $t$ of the circulant matrix associated to the vector $(0, e_1, \dots, e_{p-1})$, i.e.\ the matrix with rows being all cyclic permutations of the given vector. Under our assumptions (\ref{eq:torsion}) and (\ref{eq:regular}), for all $n \in \Nplus$
	\[
		\log_p(G/\St_G(n+1)) = t p^{n-1} + 1,
	\]
	and $\log_p(G/\St_G(1)) = 1$. Additionaly, (\ref{eq:torsion}) is equivalent to $t < p$. By \cref{Lemma}{lem:self-sim_obstructions}, for $n > 2$,
	\begin{align*}
		o_{G}(n) &= p \cdot \log_p(|L_G(n-1)|) - \log_p(|L_G(n)|)\\
		&= p \cdot \log_p \frac{p^{t \cdot p^{n-2} + 1}}{p^{t\cdot p^{n-3} + 1}} - \log\frac{p^{t \cdot p^{n-1} + 1}}{p^{t\cdot p^{n-2} + 1}} = 0
	\end{align*}
	and
	\begin{align*}
		o_{G}(2) &= p\cdot\log\frac{p^{t + 1}}{p} - \log\frac{p^{t\cdot p + 1}}{p^{t + 1}} = tp - t(p-1) = t \quad\text{and}\\
		o_{G}(1) &= p \cdot\log p - \log\frac{p^{t + 1}}{p} = p - t.
	\end{align*}
	Consequently, $\dimH G = t(p-1)/p^2$ (cf.\ \cite{FAZR14} for a more general formula).
	
	We aim to apply \cref{Proposition}{prop:obstruction}.
	Condition (\ref{eq:regular}) is equivalent to $G$ being weakly regular branch (in fact, regular branch) over $[G,G]$, by \cite{FAZR14}*{Lemma 3.4}. More precisely, we have
	\[
		\psi_1([\St_G(1),\St_G(1)]) = [G,G]^p.
	\]
	By \cref{Proposition}{prop:split-intersection_in_commutator} this implies that $K_{s-1} = [G, G]$.
	We now prove that $G$ is very strongly fractal. It is easy to see that $\St_G(1)|_x = G$ for all $x \in X$, and by \cite{FAZR14}*{Lemma 3.3} $\psi_1(\St_G(n)) = \St_G(n-1)^p$ for all $n \geq 3$. Thus it remains to check if $\St_G(2)|_x = \St_G(1)$ for all $x \in X$. By the fact that $[\St_G(1),\St_G(1)]|_x = [G,G]$ for all $x \in X$ and $[\St_G(2) : [\St_G(1), \St_G(1)]] = p^{p - t} \geq p$ (see again \cite{FAZR14}), we see that $\St_G(2)$ contains an element $g$ such that $\psi_1(g) \in \St_G(1)^p \setminus [G,G]^p$. Hence at least for one $x \in X$
	\[
		\St_G(1) \geq \St_G(2)|_x > [G,G].
	\]
	But since $[\St_G(1): [G,G]] = p$ by \cite{FAZR14}*{Theorem 2.1}, this implies $\St_G(2)|_x = \St_G(1)$, and since $G$ is spherically transitive, this holds for all $x \in X$, and $G$ is very strongly fractal. We remark that by \cite{UA16}*{Proposition 5.1} the condition (\ref{eq:torsion}) alone implies that $G$ is super strongly fractal, but our argument additionally needs (\ref{eq:regular}), since otherwise $[[G,G]^p : \psi_1([\St_G(1),\St_G(1)])] = p$ (cf.\ \cite{FAZR14}*{Lemma 3.5}).
	
	Now we may apply \cref{Proposition}{prop:obstruction} to calculate the Hausdorff dimension of $\bp_s(G)$:
	\begin{align*}
		o_{\bp_s(G)}(s) = p-t\quad\text{and}\quad o_{\bp_s(G)}(2s) = t
	\end{align*}
	and $o_{\bp_s(G)}(n) = 0$ for all other values $n \in \Nplus$, hence
	\begin{align*}
		\dimH \bp_s(G) &= 1 - \limsup_{n \to \infty} \sum_{i=1}^n \left(\frac 1 {p^i} - \frac 1 {p^{n+1}}\right)o_{\bp_s(G)}(i)\\
		&= 1 - \limsup_{n \to \infty} \left(\frac {p-t} {p^s} + \frac {t} {p^{2s}} - \frac {p - t + t} {p^{n+1}}\right)\\
		&= 1 - \left(\frac {p - t} {p^s} + \frac t {p^{2s}}\right) = \frac{p^{s-1}-1}{p^{s-1}} + \frac{t(p^{s}-1)}{p^{2s}}.
	\end{align*}
\end{eg}

\section{The generalised Basilica groups}\label{sec:gen_bas}

\noindent Let $d,\,m,\,s \in \Nplus$ with $m,\,s \geq 2$. In the subsequent \cref{sections}{sec:gen_bas}, \cref{}{sec:gen_bas_Lpres}, \cref{}{sec:gen_bas_prep} and \cref{}{sec:gen_bas_mCSP} we study the generalised Basilica groups, $\bp_s(\mathcal O_m^d)$, where $\mathcal O_{m}^d = D_d(\mathcal O_m) = \LD \pi_i(a) \mid i \in \intseg{0}{d-1} \RD$ (cf.\ \cref{Proposition}{prop:direct_prod_of_self_sim} and \cref{Definition}{def:gen_basilica}). For convenience, we use the following notation for the generators of $\bp_s(\mathcal O_m^d)$: let $i \in \intseg{0}{d-1}$ and $j \in \intseg{0}{s-1}$, and
\[\begin{array}{llll}
	a_{i,j} &\defeq \beta_j(\pi_i(a)) &= (a_{i, j-1}, \id, \dots, \id), &\text{ for } j\neq 0\\
	a_{i,0} &\defeq \beta_0(\pi_i(a)) &= (a_{i-1, s-1}, \dots, a_{i-1, s-1}), &\text{ for } i \neq 0\\
	a_{0,0} &\defeq \beta_0(\pi_0(a)) &= \sigma(a_{d-1, s-1}, \id, \dots, \id), &
\end{array}\]
where $\sigma$ is the $m$-cycle $(0\;1\;\dots\;m-1)$. For any fixed $j$, the elements $a_{i,j}$ commute and are of infinite order.

Now we prove \cref{Theorem}{thm:gen_bas_properties}, which is obtained as corollaries of results from \cref{Section}{sec:basics} and \cref{Section}{sec:split}. 

\begin{proof}[Proof of \cref{Theorem}{thm:gen_bas_properties}]
   The statements (i) and (ii) follow directly from \cref{Lemma}{lem:levtra}, \cref{Lemma}{lem:self-sim} and  \cref{Lemma}{lem:fractality}. \cref{Proposition}{prop:finite_state_bounded_invariant} together with \cref{Corollary}{coro:amenability} imply the statement (ii). The statement (iii) is a consequence of \cref{Proposition}{prop:contrating} and \cref{Corollary}{cor:word_problem}.
   Thanks to \cref{Proposition}{prop:split-intersection_in_commutator}, the group $\mathcal{O}_m^{d}$ is $s$-split. Therefore the statements (iv), (v) and (vi) follow from \cref{Corollary}{cor:abelisation}, \cref{Proposition}{prop:torfree} and \cref{Proposition}{prop:weakly_branch_by_verbal}. The proof of (vii) can easily be generalised from \cite[Proposition 4]{GZ02}. For the special case $\bp_p(\mathcal O_p)$, where $p$ is a prime, see \cite{Sas18}.
\end{proof}

We use \cref{Theorem}{thm:stab} to provide a normal generating set for the layer stabilisers of the generalised Basilica groups. This description of layer stabilisers is crucial in proving the $p$-congruence subgroup property of the generalised Basilica groups (see \cref{Section}{sec:gen_bas_mCSP}).

\begin{theorem}\label{thm:gen_bas_stab}                                                                    
	Let $n \in \N$. Write $n = s q + r$ with $r \in \intseg{0}{s-1}$ and $q = d k + l \geq 0$ with $l \in \intseg{0}{d-1}$. Then the $n$-th layer stabiliser of $B = \bp_s(\mathcal{O}_m^{d})$ is given by 
		\[
			\St_B(n) = \LD a_{i,j}^{m^{k+1}},a_{i',j'}^{m^k} \mid 0 \leq is+j \leq ls+r-1 < i' s + j' \leq ds-1\RD^B.
		\] 
	\end{theorem}
	\begin{proof}
    Let $a$ be the generator of the $m$-adic odometer $\mathcal{O}_m$. Set $G = D_d(\mathcal O_m) \cong \Z^d$. For every $i \in \intseg{0}{d-1}$, denote by $a_i = \pi_i(a)$ the generators of $G$. Since powers of the elements $a_0,\dots,a_{d-1}$ act on vertices of disjoint levels of the $m$-regular rooted tree $T$ and they commute with each other, we have
    \begin{equation*}\label{eq:stab_m-adic odometer}
       \St_G(q) = \LD a_0^{m^{k+1}},\dots,a_{l-1}^{m^{k+1}},a_l^{m^k},\dots,a_{d-1}^{m^k}\RD.
    \end{equation*}
    Now observe that for every vertex $x \in X$, $i \in \intseg{0}{d}$ and $k \in \N$,
	\begin{align*}
		a_i^{m^k}|_x &= a_{i-1}^{m^k}\\
		a_0^{m^k}|_x &= a_{d-1}^{m^{k-1}}.
	\end{align*}
	Therefore $\St_G(q)|_x = \St_G(q-1)$ and hence $G$ is very strongly fractal.
	A straightforward calculation using \cref{Theorem}{thm:stab} yields the result.
\end{proof}

Using the description of the layer stabilisers of $G$, we obtain \cref{Theorem}{thm:gen_bas_haus} as a direct application of \cref{Lemma}{lem:obstruction_to_hdim} and \cref{Proposition}{prop:obstruction}.

\begin{proof}[Proof of \cref{Theorem}{thm:gen_bas_haus}]
	The series of obstructions of $G = \mathcal O_m^d$ is constant $m-1$ for all $n \in \Nplus$, signifying Hausdorff-dimension $0$ (cf.\ \cref{Lemma}{lem:obstruction_to_hdim}). We have seen in the proof of \cref{Theorem}{thm:gen_bas_stab} that $\bp_s(G)$ is very strongly fractal. Therefore, by \cref{Proposition}{prop:obstruction} one has ${o_{\bp_s(G)}(qs)} = m-1$ for all $q \in \Nplus$ and $o_{\bp_s(G)}(n) = 0$ for all other levels.
	
	According to \cref{Lemma}{lem:obstruction_to_hdim} it is
	\begin{align*}
		\dimH \bp_s(G)
		&= 1 - \limsup_{n \to \infty} \sum_{i = 1}^n (m^{-i} - m^{-(n+1)}) o_{\bp_s(G)}(i)\\
		&= 1 - (m-1)\limsup_{n \to \infty} \left(m^{-s}\frac{1-m^{-s\lfloor n/s \rfloor}}{1-m^{-s}} - \lfloor n/s \rfloor m^{-(n+1)}\right)\\
		&= 1 - (m-1)\frac{m^{-s}}{1-m^{-s}}\\
		&= \frac{m^s - m}{m^s - 1}.
	\end{align*}
	In particular, the Hausdorff dimension is independent of $d$.
\end{proof}

\section{An $L$-presentation for the generalised Basilica group}\label{sec:gen_bas_Lpres}

\noindent Let $d,\,m,\,s \in \Nplus$ with $m,\,s \geq 2$.
In this section we will provide a concrete $L$-presentation for the generalised Basilica group $\bp_s(\mathcal{O}_m^{d})$, hence proving \cref{Theorem}{thm:gen_bas_Lpres}. We will later use this presentation to prove that all generalised Basilica groups $\bp_s(\mathcal O_p^d)$ with $p$ a prime have the $p$-congruence subgroup property.

\begin{defn}\cite{Bar03}*{Definition 1.2} \label{def:L-ptn} 
	An $L$-presentation (or an endomorphic presentation) is an expression of the form
	\begin{equation*} \label{eq:end prestn}
	   L = \LD Y \mid Q \mid \Phi \mid R \RD,
	\end{equation*}
	where $Y$ is an alphabet, $Q,R \subset F_Y$ are sets of reduced words in the free group $F_Y$ on $Y$ and $\Phi$ is a set of endomorphisms of $F_Y$. The expression $L$ gives rise to a group $G_L$ defined as 
	\begin{align*}
		G_L = F_Y / \LD Q \cup \LD\Phi\RD(R) \RD^{F_Y},
	\end{align*}
	where $\LD\Phi\RD(R)$ denotes the union of the images of $R$ under every endomorphism in the monoid $\LD\Phi\RD$ generated from $\Phi$. An $L$-presentation is finite if $Y,Q,\Phi, R$ are finite. 
\end{defn}

We now set out to prove \cref{Theorem}{thm:gen_bas_Lpres}. To do this, we follow the strategy from \cite{GZ02} which is motivated from \cite{Gri98}: let 
 \begin{equation} \label{eq:gen_bas_Y}
   Y = \{a_{i,j} \mid i\in \intseg{0}{d-1},j \in \intseg{0}{s-1}\}.
 \end{equation}
 For convenience, we do not distinguish notationally between the generators of $\bp_s(\mathcal{O}_m^{d})$ and the free generators for the presentation. 
 Observe that for a fixed $j$ the generators $a_{i,j}$ and $a_{i',j}$ of $\bp_s(\mathcal{O}_m^{d})$ commute for all $i,i' \in \intseg{0}{d-1}$. Write
  \begin{equation} \label{eq:gen_bas_Q}
     Q = \{[a_{i,j},a_{i',j}] \mid i,i' \in \intseg{0}{d-1}, \, j \in \intseg{0}{s-1}\} \subseteq F_Y
  \end{equation}
 and denote by $F$ the quotient of $F_Y$ by the normal closure of $Q$ in $F_Y$. We identify $F$ with a free product of free abelian groups
  \begin{equation*}\label{eq:free product}
	F  = \bigast_{j \in \intseg{0}{s-1}} \langle a_{i,j}\mid i \in \intseg{0}{d-1} \rangle \cong \Z^d *\cdots * \Z^d. 
  \end{equation*}
The group $\bp_s(\mathcal{O}_m^{d})$ is a quotient of $F$. Let $\mathrm{proj}: F \rightarrow \bp_s(\mathcal{O}_m^{d})$ be the canonical epimorphism. Now observe that the subgroup 
  \begin{equation}\label{eq:Delta} 
     \Delta = \langle a_{i,j}^{a_{0,0}^k}, \, a_{0,0}^m \, \mid (i,j) \in \intseg{0}{d-1} \times \intseg{0}{s-1} \backslash \{(0,0)\}, \, k \in \intseg{0}{m-1} \rangle,
  \end{equation}
is normal of index $m$ in $F$ and it is the full preimage of $\St_{\bp_s(\mathcal{O}_m^{d})}(1)$ under the epimorphism $\mathrm{proj}$ (cf.\  \cref{Theorem}{thm:gen_bas_stab}). We define a homomorphism $\Psi: \Delta \to F^m$ modelling the process of taking sections as follows:
\[
	\begin{array}{rllll}
		\Psi(a_{0,0}^m) &= (a_{d-1,s-1},\dots,a_{d-1,s-1})	 &\eqdef z_0,	&	\\
		\Psi(a_{i,0}^{a_{0,0}^k}) = \Psi(a_{i,0}) &= (a_{i-1,s-1},\dots,a_{i-1,s-1})	 &\eqdef z_{i}	     &\text{for } i \neq 0,\\
		\Psi(a_{i,j}^{a_{0,0}^k}) &= (\id^{*k}, a_{i,j-1}, \id^{*(m-k-1)})
			&\eqdef x_{i,j,k} &\text{for } j \neq 0,\\
		\Psi(a_{i,j}^{a_{0,0}^{-k}}) &= (\id^{*(m-k)}, a_{i,j-1}^{a_{d-1,s-1}^{-1}}, \id^{*(k-1)}),	
	\end{array}
\]
 where the ranges of $i,\, j$ and $k$ are as in \eqref{eq:Delta}.
 Clearly, $\ker(\Psi) \leq \ker(\mathrm{proj})$. Define
   \begin{align} 
   \alpha(v,k)  = a_{0,0}^{^{mv_0+k}}a_{1,0}^{v_1} \cdots a_{d-1,0}^{v_{d-1}} \text{ for } v = (v_0,\dots,v_{d-1}) \in \Z^d \text{ and } k \in \intseg{0}{m-1}, \label{eq:gen_bas_alpha(v,k)}\\
   R  = \{
			[a_{i,j},a_{i',j'}^{\alpha(v,k)}]
			\mid i,i' \in \intseg{0}{d-1},
			\, 
			j,j' \in \intseg{1}{s-1} , 
			\,
			k \in \intseg{1}{m-1}, 
			\,
			v \in \Z^d
		\} \label{eq:gen_bas_R},
  \end{align}
where by abuse of notation we interpret $\alpha(v,k)$ and $r \in R$ both as elements of $F_Y$ and their images in $F$. We will prove in \cref{Proposition}{prop:ptn_delta} that the kernel of $\Psi$ is normally generated from the image of $R$ in $F$, implying that the set $R$ belongs to the set of defining relators of $\bp_s(\mathcal{O}_m^{d})$. By definition of the elements $a_{i,j}$, we may obtain the elements of the set $R$ as vertex sections. To incorporate these elements to the set of defining relators we introduce the following endomorphism of $F_Y$ defined as 
 \begin{align} \label{eq:gen_bas_phi}
		\Phi : \begin{cases}
			a_{i,j} & \mapsto a_{i,j+1} \text{ for } j \neq s-1,\\
			a_{i,s-1} & \mapsto a_{i+1,0} \text{ for } i \neq d-1,\\
			a_{d-1,s-1} & \mapsto a_{0,0}^m,
		\end{cases}
 \end{align}
 where $i \in \intseg{0}{d-1}$ and $j \in \intseg{0}{s-1}$.
 
  \begin{theorem} \label{thm:bas_l_pres}
    The generalised Basilica group admits the $L$-presentation 
    \[
     L = \LD Y \mid Q \mid \Phi \mid R \RD
    \]
    where $Y,\, Q, \, R$ and $\Phi$ are given by \eqref{eq:gen_bas_Y}, \eqref{eq:gen_bas_Q}, \eqref{eq:gen_bas_R} and \eqref{eq:gen_bas_phi}.
  \end{theorem}

 Observe that for any $g \in Q$ and $r \in \N$, it holds that $\Phi^r(g) \in \langle Q^{F_Y} \rangle$.
	Considering the presentation defining $F$ we may assume that $\Phi$ is an endomorphism of $F$ and that $R$ is a subset of $F$. To prove \cref{Theorem}{thm:bas_l_pres}, it is enough to show that $\ker(\Psi) = \langle R^F \rangle$ and $\ker(\mathrm{proj}) = \bigcup_{r \in \N} \Phi^r(R)$. We will obtain the first part from \cref{Proposition}{prop:ptn_delta} and the latter from \cref{Lemma}{lem:bas_lpres_lemma_2} to \cref{Lemma}{lem:bas_lpres_lemma_4}.

\begin{prop}\label{prop:ptn_delta}
	Let $\Tilde{\Delta}$ be the image of $\Delta$ under $\Psi$. Let $z^v$ be the product $z_0^{v_0}\cdots z_{d-1}^{v_{d-1}}$ for every $v = (v_0,\dots,v_{d-1})  \in \Z^d$. Then $\Tilde{\Delta}$ admits the presentation
	\[\LD\ \mathcal{S} \mid \mathcal{R}\ \RD\] 
	where $\mathcal{S} = \{x_{i,j,k},\, z_i \mid i \in \intseg{0}{d-1},
			\, 
			j \in \intseg{1}{s-1}, 
			\, 
			k \in \intseg{0}{m-1}\}$
		and 
		\begin{align*}
		\mathcal{R} = 
		\left\langle 
		\begin{array}{c|l}
			[x_{i,j,k},x_{i',j,k}], \, [x_{i,j,k},x_{i',j',k'}^{z^{v}}],
		    & i,i' \in \intseg{0}{d-1},
			\, 
			j,j' \in \intseg{1}{s-1},\\ 
			{[}z_i,z_{i'}{]}
			&
			k, k'\in \intseg{0}{m-1} \text{ with } k \neq k',
			\,
			v \in \Z^d
		\end{array}\right\rangle.
		\end{align*}
	As a consequence, we obtain that 
	\[
		\ker(\Psi) = \LD \{
		    [a_{i,j},a_{i',j'}^{\alpha(v,k)}]
			\mid i,i' \in \intseg{0}{d-1},
			\, 
			j,j' \in \intseg{1}{s-1} , 
			\,
			k \in \intseg{1}{m-1}, 
			\,
			v \in \Z^d \} \RD^F ,
	\] 
	where $\alpha(v,k)$ is given by \eqref{eq:gen_bas_alpha(v,k)}.
\end{prop}

\begin{proof}
	Let $A = \langle a_{i,j} \mid i \in \intseg{0}{d-1}, \, j \in \intseg{0}{s-2} \rangle^F$ and $Z  = \langle z_0,\dots,z_{d-1} \rangle \cong \Z^d$ be subgroups of $F$ and $\Tilde{\Delta}$ respectively. Notice that $\Tilde{\Delta}$ is a sub-direct product of $m$ copies of $F$ and the elements $x_{i,j,k}$ and $ x_{i',j',k'}$ commute if $k \neq k'$ or if $k =k'$ and $j = j'$. It follows from the definition of $\Psi$ that 
	\begin{align*}
		A^m = \left\langle 
		\begin{array}{l|l}
			x_{i,j,k}  
			& i \in \intseg{0}{d-1},
			\, 
			j \in \intseg{1}{s-1} , 
			\,
			k \in \intseg{0}{m-1}
		\end{array}
		\right\rangle^{\Tilde{\Delta}} \leq \Tilde{\Delta}.
	\end{align*}
	Hence 
 	$\tilde{\Delta} = A^mZ$, yielding $\tilde{\Delta} = A^m \rtimes Z$. Now, since $F$ is a free product of free abelian groups, the group $A$ is freely generated from the elements of the form 
 	  \[
 	   a_{i,j}^{a_{d-1,s-1}^{v_0}a_{0,s-1}^{v_1}\cdots a_{d-2,s-1}^{v_{d-1}}},
 	  \]
 	where $v_i \in \Z, \, i \in \intseg{0}{d-1}$ and $j \in \intseg{0}{s-2}$. Therefore,
	the group $A^m$ is generated from the elements   \[
	   x_{i,j,k}^{z^v} = (\id^{*k}, a_{i,j-1}^{{a_{d-1,s-1}^{v_0}a_{0,s-1}^{v_1}\cdots a_{d-2,s-1}^{v_{d-1}}}}, \id^{*(m-k-1)}),
	 \]
	where $i \in \intseg{0}{s-1}, \, j \in \intseg{1}{s-1}, \, k \in \intseg{0}{m-1}$ and
	\[
	  z^v = z_0^{v_0}\cdots z_{d-1}^{v_{d-1}} = ({a_{d-1,s-1}^{v_0}a_{0,s-1}^{v_1}\cdots a_{d-2,s-1}^{v_{d-1}}}, \dots, {a_{d-1,s-1}^{v_0}a_{0,s-1}^{v_1}\cdots a_{d-2,s-1}^{v_{d-1}}}),
	\]
	with $v_i \in \Z$.
	We obtain a presentation of $A^m$ as 
	\begin{align*}
	    \left\langle 
		\begin{array}{l|l}
			x_{i,j,k} ^{z^v}  
			& [x_{i,j,k},x_{i',j,k}] = [x_{i,j,k}^{z^v},x_{i',j',k'}^{z^{v'}}]=\id, 
			\, i,i' \in \intseg{0}{d-1},
			\, j,j' \in \intseg{1}{s-1},\\
			& 
			k,k' \in \intseg{0}{m-1} \text{ with } k \neq k',
			\, v,v' \in \Z^d 
		\end{array}
		\right\rangle.
	\end{align*}
	Hence $\Tilde{\Delta}$, being a semi-direct product, admits the presentation $\LD \mathcal{S}\mid \mathcal{R} \RD$, since conjugating an element $x_{i,j,k}$ by $z_i$ does not yield a new relation.
	Therefore, the kernel of $\Psi$ is normally generated from the preimage of the set of defining relators for $\tilde{\Delta}$. Notice that the preimages of the elements $[z_i,z_{i'}]$ and $[x_{i,j,k},x_{i',j,k}]$ are trivial in $\Delta$. Hence,
	\begin{align*}
		\ker(\Psi) & = \left\langle 
		\begin{array}{l|l}
			[a_{i,j}^{\alpha(v,k)},a_{i',j'}^{\alpha(v',k')}]  
			&  
			\, i,i' \in \intseg{0}{d-1},
			\, j,j' \in \intseg{1}{s-1}, \\
			&
			k,k' \in \intseg{0}{m-1} \text{ with } k \neq k',
			\, v,v' \in \Z^d 
		\end{array}
		\right\rangle^{\Delta}.
	\end{align*}
	Indeed, $\ker(\Psi)$ is normal in $F$. Given $v \in \Z^d$ and $k \in \intseg{0}{m-1}$, define
	\begin{align*}
		\underline{v} &= (\lfloor (mv_0 + k + 1)/m \rfloor, v_1, \dots, v_{d-1}) \in \Z^d \quad\text{and}\\
		\underline{k} &= k + 1 \pmod{m} \in \intseg{0}{m-1}.
	\end{align*}
	Then
	\begin{align*}
	   \alpha(v,k)a_{0,0} 
	   &
	   = a_{0,0}^{^{mv_0+k + 1}}a_{1,0}^{v_1} \cdots a_{d-1,0}^{v_{d-1}} = \alpha(\underline{v},\underline{k})\\
	   \alpha(v',k')a_{0,0}
	   &
	   = a_{0,0}^{^{mv'_0+k' + 1}}a_{1,0}^{v'_1} \cdots a_{d-1,0}^{v'_{d-1}} = \alpha(\underline{v}',\underline{k}')
	\end{align*}
	implies
	\[
	   [a_{i,j}^{\alpha(v,k)},a_{i',j'}^{\alpha(v',k')}]^{a_{0,0}} = [a_{i,j}^{\alpha(v,k)a_{0,0}},a_{i',j'}^{\alpha(v',k')a_{0,0}}] = [a_{i,j}^{\alpha(\underline{v},\underline{k})},a_{i',j'}^{\alpha(\underline{v'},\underline{k'})}] \in \ker(\Psi).
	\]
	A similar calculation shows $[a_{i,j}^{\alpha(v,k)},a_{i',j'}^{\alpha(v',k')}]^{a_{0,0}^{-1}} \in \ker(\Psi)$.
	We get
	\[
		\ker(\Psi) = \left\langle 
		\begin{array}{l|l}
			[a_{i,j},a_{i',j'}^{\alpha(v,k)}]  
			&  
			\, i,i' \in \intseg{0}{d-1},
			\, j,j' \in \intseg{1}{s-1}, \, 
			k \in \intseg{1}{m-1},
			\, v \in \Z^d 
		\end{array}
		\right\rangle^{F}.\qedhere
	\]
\end{proof}

\notation\label{notation tau(v,k)} 
	Let $i,i'\in \intseg{0}{d-1},\, j,j' \in \intseg{1}{s-1},\, k \in \intseg{1}{m-1}, \,  v \in \Z^d$ and $n \in \N$. Define
	\begin{align*}
		\Omega_0 &\defeq \ker (\Psi), &&
		\Omega_n \defeq \Psi^{-1}(\Omega_{n-1}^m) \text{ for } n \geq 1,\\
		\tau_{v,k}(i,j,i',j') &\defeq [a_{i,j},a_{i',j'}^{\alpha(v,k)}], && X_n \defeq \langle \Phi^r(\tau_{v,k}(i,j,i',j')) \mid r \in \intseg{0}{n} \rangle^{F},
	\end{align*}
	where $\alpha(v,k)$ is given by \eqref{eq:gen_bas_alpha(v,k)}. Denote further by $\Omega$ the kernel of the epimorphism $\mathrm{proj}: F \rightarrow \bp_s(\mathcal{O}_m^{d})$. We will prove $\Omega_n = X_n$ and $\Omega = \bigcup \limits ^{\infty}_{n=0}\Omega_n$, proving \cref{Theorem}{thm:bas_l_pres}.

\begin{lemma} \label{lem:bas_lpres_lemma_2}
	For $w \in F'$ the identity  $\Psi(\Phi(w)^{a_{0,0}^k})= (\id^{*k},w,\id^{*(m-k-1)})$ holds for every $k \in \intseg{0}{m-1}$.
\end{lemma}

 \begin{proof}
	Observe from the definition of $\Phi$ that
	 \[
	   \Phi(F) = \LD a_{i,j},a_{0,0}^m \mid (i,j) \in \intseg{0}{d-1} \times \intseg{0}{s-1} \backslash\{(0,0)\} \RD \leq \Delta.
	 \]
	 Then by direct calculation using the definition of the homomorphism $\Psi$ and $\Phi$ we get the desired identity.
 \end{proof}
 
 \begin{lemma}\label{lem:bas_lpres_lemma_3}
	The equality $\Omega_n = X_n$ holds for all $n \in \N$.
 \end{lemma}

 \begin{proof}
	It follows from \cref{Proposition}{prop:ptn_delta} that $\Omega_0 = \ker(\Psi)=X_0$.
	The proof proceeds by induction on $n$.
	Since $\Phi(F) \leq \Delta$, for every $r \in \mathbb{N}_0$, we have $\Phi^r(\Delta)\leq \Delta$. Hence $X_n \leq \Delta$ for all $n \in \N$. Assume for some $n \geq 1$ that $\Omega_{n-1}=X_{n-1}.$  We will prove that 
 \[
   \Psi(X_n)= \Omega_{n-1}^m=\Psi(\Omega_n).
 \]
	Let $i,i' \in \intseg{0}{d-1}, \, j,j' \in \intseg{1}{s-1} , \, k \in \intseg{1}{m-1}, \, r \in \intseg{1}{n}$ and $v \in \Z^d$. For every $\Phi^r(\tau_{v,k}(i,j,i',j')) \in X_n$ and for every $\ell \in \intseg{0}{m-1}$, since $\Phi^{r-1}(\tau_{v,k}(i,j,i',j')) \in F'$, we obtain from \cref{Lemma}{lem:bas_lpres_lemma_2} that
	\begin{align*}
		\Psi((\Phi^r(\tau_{v,k}(i,j,i',j')))^{a_{0,0}^{\ell}}) 
		&= 		(\id^{*\ell}, \Phi^{r-1}(\tau_{v,k}(i,j,i',j')), \id^{*(m-\ell-1)}).
	\end{align*}
	Since $\Tilde{\Delta}$ is a sub-direct product of $m$ copies of $F$ and $X_{n-1}$ is normally generated from the elements of the form $\Phi^{r-1}(\tau_{v,k}(i,j,i',j'))$, we obtain that $\Psi(X_n) = \Omega_{n-1}^m = \Psi(\Omega_n)$. 
	
	But since $\ker(\Psi|_{\Omega_n})=\ker(\Psi)\cap \Omega_n=\Omega_0= X_0 = \ker(\Psi)\cap X_n=\ker(\Psi|_{X_n})$, we get $\Omega_n = X_n$, and the result follows by induction.
\end{proof}

\begin{lemma}\label{lem:bas_lpres_lemma_4}
	We have $\Omega = \bigcup \limits_{n=0}^{\infty}\Omega_n$.
\end{lemma}

\begin{proof}
	Write $B$ for $\bp_s(\mathcal{O}_m^d)$ and recall that $\mathrm{proj} : F \rightarrow B$ is the canonical epimorphism. Notice that $\St_{B}(1)$ is a quotient of $\Delta$ and further $\Omega_0 = \ker(\Psi)\leq \ker(\mathrm{proj}) = \Omega$. Proceeding by induction on $n$, we will prove that $\bigcup \limits_{n=0}^{\infty}\Omega_n \leq \Omega$. Assume that $\Omega_{n-1} \leq \Omega$ for some $n \geq 1$. Let $w \in \Omega_n$ and let $w_k$ be the $k$-th component of $\Psi(w)$. Then $w_k \in \Omega_{n-1}$ for all $k \in \intseg{0}{m-1}$. Then the first layer sections of $\mathrm{proj}(w)\in \St_{B}(1)$ act trivially on the subtrees hanging from the vertices of level one of the $m$-regular rooted tree. Hence $\mathrm{proj}(w)$ acts trivially and $\mathrm{proj}(w) = \id$ in $B$. It follows by induction that $\Omega_n \leq \Omega$ for all $n \in \N$. Since $\Omega_{n-1} \leq \Omega_n$ for all $n \in \Nplus$, we obtain $\bigcup \limits_{n=0}^{\infty}\Omega_n \leq \Omega$.
	
	Now, to see the converse choose an arbitrary element $w \in F$ such that $\mathrm{proj}(w) = \id$ in $B$. Then by \cref{Theorem}{thm:gen_bas_stab} $\mathrm{proj}(w) \in \St_B(1)$ and hence $w \in \Delta$. Denote by $w_k$ the $k$-th component of $\Psi(w)$. Then $\mathrm{proj}(w) = \id$ if and only if $\mathrm{proj}(w_k) = \id$ for all $k \in \intseg{0}{m-1}$, implying that $w_k \in \Delta$ for all $k \in \intseg{0}{m-1}$. Now repeat this process of taking sections by replacing $w$ with $w_k$. This process is equivalent to the algorithm solving the word problem for $B$, cf. \cite[Proposition 5]{GZ02}. Thanks to \cref{Corollary}{cor:word_problem}, the word problem for $B$ is solvable and hence this process terminates in a finite number of steps. This implies the existence of an element $n \in \mathbb{N}_0$ such that $w \in \Omega_n$, completing the proof.
\end{proof}

To conclude this section, we want to point out that akin to \cite[Proposition 11]{GZ02}, one can introduce a set of $d$ endomorphisms, each corresponding to a generator $a_{i,0}$, and obtain a finite $L$-presentation for $\bp_s(\mathcal{O}_m^{d})$. We omit the proof of \cref{Theorem}{thm:gen_bas_endopres} below due to its technicality, but a rigorous proof can be found in the PhD Dissertation of the second author.

 \begin{theorem}\label{thm:gen_bas_endopres}
 The group $\bp_s(\mathcal{O}_m^{d})$ admits the following $L$-presentation:
 \[
	\left\langle
	\begin{array}{c|c|c|c}
		a_{i,j} & [a_{i,j},a_{i',j}] && [a_{i,j},a_{i',j'}^{\alpha(v,k)}],i,i'\in \intseg{0}{d-1},\\
		i \in \intseg{0}{d-1} & i,i'\in \intseg{0}{d-1} 
		& \Phi,\Theta_0,\dots,\Theta_{d-1} 
		&  j,j' \in \intseg{1}{s-1}, \, k \in \intseg{1}{m-1}\\
		j \in \intseg{0}{s-1} & j \in \intseg{0}{s-1} 
		&& v \in \{0\} \times  \{0,1\}^{d-1}
	\end{array}
	\right\rangle
\]
 where $\alpha(v,k)$ and $\Phi$ are given by \eqref{eq:gen_bas_alpha(v,k)} and \eqref{eq:gen_bas_phi}, respectively, and $\Theta_{i'}$ are endomorphisms of the free group on the set of generators defined as
  \begin{align*}
    \Theta_{i'} : \begin{cases}                    
              a_{i,j} \mapsto a_{i,j}a_{i,j}^{a_{i',0}} \text{ for } j \neq 0, i' \neq 0,\\
              a_{i,j} \mapsto a_{i,j}a_{i,j}^{a_{0,0}^m} \text{ for } j \neq 0, i' = 0,\\
              a_{i,0} \mapsto a_{i,0}.
             \end{cases}
    \end{align*}
\end{theorem}

\section{Structural properties of the generalised Basilica groups}\label{sec:gen_bas_prep}

\noindent Let $d,\,m,\,s \in \Nplus$ with $m,\,s \geq 2$. Here we prove some structural properties of the generalised Basilica groups $\bp_2(\mathcal{O}_m^{d})$. These result reflect a significant structural dissimilarity between $\bp_2(\mathcal{O}_m^{d})$ and $\bp_s(\mathcal{O}_m^{d})$ for $s>2$.
This structural dissimilarity plays a vital role when we consider the $p$-congruence subgroup property of the generalised Basilica groups, see \cref{Figure}{fig:proofMCSP}, which is treated in \cref{Section}{sec:gen_bas_mCSP}.

For convenience, we omit the subscript from $\psi_1$ and identify an element $g \in \St_{B}(1)$ with its image under the map $\psi_1$.

\begin{prop}\label{prop:gen_bas_commsbgp}
  Let $B$ be the generalised Basilica group $\bp_s(\mathcal{O}_m^{d})$. Then $\psi^{-1}((B')^m)$ is a subgroup of $B'$ and
  \begin{align*}
    B'/\psi^{-1}((B')^m) 
        & =\left\langle 
        \begin{array}{l|l}
         c_{i,j,k} \, \psi^{-1}((B')^m)  
         & i \in \intseg{0}{d-1},  j \in \intseg{1}{s-1}, 
         \, k \in \intseg{1}{m-1}
         \end{array}
         \right \rangle  \\
         & \cong \Z^{d(m-1)(s-1)},
         \end{align*}
         where $c_{i,j,k} = [a_{i,j},a_{0,0}^{k}]$. In particular, it holds that $\psi ^{-1}((B')^m) \geq B''.$
\end{prop}

\begin{proof} 
 Notice that $B' = \langle [a_{i,j},a_{i',j'}] \mid i,i' \in \intseg{0}{d-1}, \, j,j' \in \intseg{0}{s-1} \rangle ^B$. For $i,i' \in \intseg{0}{d-1}$ and $j,j' \in \intseg{1}{s-1}$, we have $[a_{i,j},a_{i',j}] = \id$ and for $j \neq j'$
  \begin{align*}
    [a_{i,j},a_{i',j'}] &= ([a_{i,j-1},a_{i',j'-1}],\id^{*(m-1)}) \\
    [a_{i,j},a_{i',0}] &= ([a_{i,j-1},a_{i'-1,s-1}],\id^{*(m-1)}) \text{ for } i' \neq 0,\\
    [a_{i,j},a_{0,0}^m] &=
    ([a_{i,j-1},a_{d-1,s-1}],\id^{*(m-1)}).
  \end{align*}
 Therefore, we obtain
  \[
    \langle [a_{i,j},a_{i',j'}] \mid i,i' \in \intseg{0}{d-1}, \, j,j' \in \intseg{0}{s-1} \rangle \times \{\id \} \times \cdots \times \{\id \} \leq \psi(B'),
  \]
  yielding that $(B')^m \leq \psi(B')$ by \cref{Lemma}{lem:fractal_branching}.
  
 Now, recall our definition $c_{i,j,k} = [a_{i,j},a_{0,0}^{k}]$ and
	\[
		C = \langle c_{i,j,k} \mid i \in \intseg{0}{d-1}, \, j \in \intseg{1}{s-1}, \, k \in \intseg{1}{m-1} \rangle.
	\]
	We claim that $B'/\psi^{-1}((B')^m) = \overline{C}$, where $\overline{C}$ denotes the image of $C$ in the quotient group. For convenience, we will write the equivalence $\equiv_{\psi^{-1}((B')^m)}$ without the subscript. Observe that, for $i,i' \in \intseg{0}{d-1}, \, j,j' \in \intseg{1}{s-1}$ and $k \in \intseg{1}{m-1}$,
	\begin{align*}
	   [a_{i,j},a_{i',j'}] \equiv \id, &&  [a_{i,j},a_{i',0}] \equiv \id \text{ for } i' \neq 0, && [a_{i,j},a_{0,0}] = c_{i,j,1},
	\end{align*}
	and 
	\begin{align*}
		c_{i,j,k}=[a_{i,j},a_{0,0}^{k}]\equiv (a_{i,j-1}^{-1},\id^{*(k-1)},a_{i,j-1},\id^{*(m-k-1)}).
	\end{align*}
 Therefore, to prove the claim, it suffices to show that $\overline{C}$ is normal in $B/\psi^{-1}((B')^m)$. Let $i,i' \in \intseg{0}{d-1}, \, j,j' \in \intseg{1}{s-1}$ and $k \in \intseg{1}{m-1}$. An easy calculation yields
	\begin{align*}
		c_{i,j,k}^{a_{i',j'}^{\pm 1}} 
		\equiv c_{i,j,k} && \text{and} &&
		c_{i,j,k}^{a_{i',0}^{\pm 1}}  
		\equiv c_{i,j,k} \text{ for } i' \neq 0.
	\end{align*}
	Furthermore, 
	\[\begin{array}{lll}
		c_{i,j,k}^{a_{0,0}} \equiv&   (\id, a_{i,j-1}^{-1},\id^{*(k-1)},a_{i,j-1},\id^{*(m-k-2)})	& \equiv c_{i,j,1}^{-1}c_{i,j,k+1}	
		\hspace{1.2cm} \text{ if } k \neq m-1,\\
		c_{i,j,k}^{a_{0,0}} \equiv&   (a_{i,j-1}, a_{i,j-1}^{-1},\id^{*(m-2)})	& \equiv
		c_{i,j,1}^{-1}	           
		\hspace{2.3cm} \text{ if } k = m-1,\\
		c_{i,j,k}^{a_{0,0}^{-1}} \equiv& (\id^{*(k-1)},a_{i,j-1},\id^{*(m-k-1)},a_{i,j-1}^{-1}) &\equiv
		\begin{cases}
		c_{i,j,m-1}^{-1}c_{i,j,k-1}   &\text{ if } k \neq 1, \\
		c_{i,j,m-1}^{-1}     
        & \text{ if } k = 1,
		\end{cases}
	\end{array}\]
	implying that $B'/\psi^{-1}((B')^m) = \overline{C}$. Observe that, for a fixed $i \in \intseg{0}{d-1}$ and $j \in \intseg{1}{s-1}$,
	\[
		\Z^{m-1} \cong \{(a_{i,j-1}^{x_1},\dots,a_{i,j-1}^{x_m}) \mid  x_r \in \Z, \sum \limits_{r=1}^{m} x_r=0\} = \langle \overline{c}_{i,j,k} \mid k \in \intseg{1}{m-1} \rangle \leq \overline{C}.
	\]
	Since $B/B' \cong \Z^{ds}$ (\cref[(iv)]{Theorem}{thm:gen_bas_properties}), this yields 
	\[
		B'/\psi^{-1}((B')^m) = \overline{C} = \prod \limits_{(i,j) \in \intseg{0}{d-1}\times\intseg{1}{s-1}}\langle \overline{c}_{i,j,k} \mid k \in \intseg{1}{m-1} \rangle \cong \Z^{d(m-1)(s-1)}.
	\]
\end{proof}

Now we prove \cref{Theorem}{thm:gamma}. In addition, we provide a generating set for the quotient group $\gamma_2(\bp_s(\mathcal{O}_m^{d}))/\gamma_3(\bp_s(\mathcal{O}_m^{d}))$.
\begin{theorem}\label{thm:gen_bas_gammasbgp}
 Let $B$ be the generalised Basilica group $\bp_s(\mathcal{O}_m^{d})$. We have:
 \begin{enumerate}
    \item[(i)] For $s = 2$, $B'/\gamma_3(B) = \langle [a_{i,0},a_{i',1}] \, \gamma_3(B) \mid i,i' \in \intseg{0}{d-1} \rangle \cong \Z^{d^2}.$
    \item[(ii)] For $s > 2$, the quotient group $B'/\gamma_3(B) \cong \mathrm C_m^{ds-2} \times \mathrm C_{m^2}$. Moreover, it is generated from the set 
     \begin{align*}
      \{[a_{i,j},a_{0,0}] \, \gamma_3(B), [a_{0,1},a_{i',0}] \, \gamma_3(B)
	    \mid 
	  i \in \intseg{0}{d-1}, \, i' \in \intseg{1}{d-1}, \, j \in \intseg{1}{s-1}\}.
	 \end{align*}
 \end{enumerate}
\end{theorem} 

\begin{proof}
(i) We use \cref{Theorem}{thm:bas_l_pres} to obtain a presentation for $B/\gamma_3(B)$. Take $Y, \, Q, \, \Phi$ and $R$ as given in \cref{Theorem}{thm:bas_l_pres} and set $Q' = Q \cup \gamma_3(F_Y)$, where $F_Y$ is the free group on $Y$. If $s=2$, the set $R$ becomes 
    \[
    R  = \{
			[a_{i,1},a_{i',1}^{\alpha(v,k)}]
			\mid i,i' \in \intseg{0}{d-1},
			\,
			k \in \intseg{1}{m-1}, 
			\,
			v \in \Z^d
		\} 
    \]
 and for every $[a_{i,1},a_{i',1}^{\alpha(v,k)}] \in R$,    
	\[
		[a_{i,1},a_{i',1}^{\alpha(v,k)}]\equiv_{\gamma_3(F_Y)}[a_{i,1},a_{i',1}]\in \LD Q' \RD^{F_Y},
	\]
	where $\alpha(v,k)$ is given by \eqref{eq:gen_bas_alpha(v,k)}. Since $\LD Q'\RD$ is invariant under $\Phi$, the presentation $\langle\ Y \mid Q'\ \rangle$ defines the group $B/\gamma_3(B)$, yielding that
	\[
		B'/\gamma_3(B) = \langle [a_{i,0},a_{i',1}] \mid i,i' \in \intseg{0}{d-1} \rangle  \cong \Z^{d^2}.
	\]
 (ii) Consider again $Y, \, Q, \, \Phi$ and $R$ as given in \cref{Theorem}{thm:bas_l_pres} and $Q' = Q \cup \gamma_3(F_Y)$. First observe that the element
	\[
		[a_{i,j},a_{i',j'}^{\alpha(v,k)}]\equiv_{\gamma_3(F_Y)}[a_{i,j},a_{i',j'}] 
	\]
	 belongs to $\LD Q' \RD^{F_Y}$ if and only if $j = j'$. Setting
	 \[S = \{ [a_{i,j},a_{i',j'}] \mid i,i'\in \intseg{0}{d-1}, \, j,j' \in \intseg{1}{s-1} \text{ with } j \neq j'\} \subseteq F_Y,\]
	 we notice that the group $B/\gamma_3(B)$ admits the $L$-presentation $\langle\ Y \mid Q' \mid \Phi \mid S\ \rangle$. 
	 Now, define
	 \begin{align*}
	 	T = \left\{
		\begin{array}{c|c}
			[a_{i,j},a_{i',0}], [a_{i'',1},a_{i',0}], & i \in \intseg{0}{d-1}, \\
			\left[a_{i,j}, a_{0,0}\right]^m, [a_{i',1},a_{0,0}]^m,[a_{0,1},a_{i',0}]^m, & i', i'' \in \intseg{1}{d-1}, \\
			\left[a_{0,1},a_{0,0}\right]^{m^2} &  j \in \intseg{2}{s-1}
	 	\end{array}\right\}
	 \end{align*}
	and $N =  Q' \cup S \cup T$ as subsets of $F_Y$. We claim that $\Phi^r(S)\subseteq  N^{F_Y}$ for all $r \in \N$, and hence the presentation $\langle\ Y \mid N\ \rangle$ defines the group $B/\gamma_3(B)$.
	Therefore, the commutator subgroup of $B/\gamma_3(B)$ is generated from the set 
	\[
		\left\{
		\begin{array}{c|c}
			[a_{i,j},a_{0,0}], [a_{i',1},a_{0,0}], & i \in  \intseg{0}{d-1},\, i' \in  \intseg{1}{d-1},\\
			\left[a_{0,1},a_{i',0}\right],[a_{0,1},a_{0,0}] & j \in \intseg{2}{s-1}
		\end{array}
		\right\},
	\]
	yielding that:
	\[
		B'/\gamma_3(B)  \cong \mathrm C_m^{d(s-2)} \times \mathrm C_m^{d-1}\times \mathrm C_m^{d-1}\times \mathrm C_{m^2}=\mathrm C_m^{ds-2} \times \mathrm C_{m^2}. 
	\]
	Now, let $i,i' \in \intseg{0}{d-1}$. Observe first that, for $j ,j' \in \intseg{1}{s-2},$ 
	\[
		\Phi([a_{i,j},a_{i',j'}])=[a_{i,j+1},a_{i',j'+1}]\in S.
	\]
	
	To prove the claim, it is enough to consider the elements of the form $\Phi^r([a_{i,j},a_{i',j'}])$ with either $j$ or $j'$, but not both, equal to $s-1$. Without loss of generality suppose that $1 \leq j \leq s-2$ and $j' =s-1$. Since $\gamma_3(F_Y) \leq N^{F_Y}$, we work modulo $\gamma_3(F_Y)$. We have 
	\begin{align*}
		\Phi([a_{i,j},a_{i',s-1}])   
		   & \equiv \begin{cases}
		   	[a_{i,j+1},a_{i'+1,0}]^{m} &\text{ if }i' = d-1\\
		   	[a_{i,j+1},a_{i'+1,0}] &\text{otherwise.}
		   \end{cases}
	\end{align*}	  
	For convenience, the images of $\Phi^2([a_{i,j},a_{i',s-1}])$ and $\Phi^3([a_{i,j},a_{i',s-1}])$ are given in the tabular form, see \cref{Table}{tab:images of Phi^2} and \cref{Table}{tab:images of Phi^3}.
	 
	 \begin{table}[H]
	 \scalebox{1.13}{
	    \centering
	    \renewcommand{\arraystretch}{1.35}
	 \begin{tabular}{ |c|c|c|c| }
    \hline
     &  & $j \neq s-2$ & $j = s-2$ \\
    \hline 
    
    \multirow{2}{*}{$i' \neq d-1$}&  $i \neq d-1$ & \multirow{2}{*}{$[a_{i,j+2},a_{i'+1,1}]$} & $[a_{i+1,0},a_{i'+1,1}]$ \\ \cline{2-2}\cline{4-4}
    
    & $i = d-1$ &  &  $[a_{0,0},a_{i'+1,1}]^m$ \\ 
    
    \hline
    
    \multirow{2}{*}{$i' = d-1$}&  $i \neq d-1$ & \multirow{2}{*}{$[a_{i,j+2},a_{0,1}]^m$} & $[a_{i+1,0},a_{0,1}]^m$ \\ \cline{2-2}\cline{4-4}
    
    & $i = d-1$ &  &  $[a_{0,0},a_{0,1}]^{m^2}$ \\ 
    
    \hline
    
    \end{tabular}}
    \caption{Images of $\Phi^2([a_{i,j},a_{i',s-1}])$.}
	    \label{tab:images of Phi^2}
	\end{table}
	
	 \begin{table}[H]
	    \centering
	    \renewcommand{\arraystretch}{1.35}
	    \scalebox{1.13}{
	 \begin{tabular}{ |c|c|c|c|c| }
    \hline
    
     &  & $j \notin \{s-3,s-2\}$ & $j = s-2$ & $j = s-3$ \\
     
    \hline
    
    \multirow{2}{*}{$i' \neq d-1$}&  $i \neq d-1$ & \multirow{2}{*}{$[a_{i,j+3},a_{i'+1,2}]$} & $[a_{i+1,1},a_{i'+1,2}]$ & $[a_{i+1,0},a_{i'+1,2}]$ \\ \cline{2-2}\cline{4-5}
    
    & $i = d-1$ &  &  $[a_{0,1},a_{i'+1,2}]^m$ &  $[a_{0,0},a_{i'+1,2}]^m$ \\ 
    
    \hline\multirow{2}{*}{$i' = d-1$}&  $i \neq d-1$ & \multirow{2}{*}{$[a_{i,j+3},a_{0,2}]^m$} & $[a_{i+1,1},a_{0,2}]^m$ & $[a_{i+1,0},a_{0,2}]^m$ \\ \cline{2-2}\cline{4-5}
    
    & $i = d-1$ &  &  $[a_{0,1},a_{0,2}]^{m^2}$ &  $[a_{0,0},a_{0,2}]^{m^2}$ \\ 
    \hline
    \end{tabular}}
    \caption{Images of $\Phi^3([a_{i,j},a_{i',s-1}])$.}
	    \label{tab:images of Phi^3}
	\end{table}
	
	Observe that the element $\Phi^r([a_{i,j},a_{i',s-1}])\in N^{F_Y}$ for $r\in \intseg{1}{3}$. By iterating the process we see that $\Phi^r([a_{i,j},a_{i',j'}]) \in N^{F_Y}$, for all $r \in \N$ and $[a_{i,j},a_{i',j'}] \in S$.
\end{proof}

\begin{lemma}\label{lem:gen_bas_secondcommsbgp}
	Let $B$ be the generalised Basilica group $\bp_s(\mathcal{O}_m^{d})$. The following assertions hold:
	\begin{enumerate}
		\item[(i)] For $s = 2$, $B'' = \psi ^{-1}(\gamma_3(B)^m)$.
		\item[(ii)] For $s > 2$, $B'' \gneq \psi ^{-1}(\gamma_3(B)^m)$.
	\end{enumerate}
\end{lemma}

\begin{proof}
   We first prove that $\gamma_3(B)^m \leq \psi(B'')$ for all $s \geq 2$. From \cref{Lemma}{lem:fractal_branching}, since
	$$
		\gamma_3(B) = \langle [[a_{i_1,j_1},a_{i_2,j_2}],a_{i_3,j_3}] \mid i_1, i_2, i_3 \in \intseg{0}{d-1}, \, j_1, j_2, j_3 \in \intseg{0}{s-1} \rangle ^{B},
	$$
	and $B$ is self-similar and fractal (\cref[(ii)]{Theorem}{thm:gen_bas_properties}), it is enough to prove that the set
	\begin{equation}
		\{([[a_{i_1,j_1},a_{i_2,j_2}],a_{i_3,j_3}],\id^{*(m-1)}) \mid i_1,i_2,i_3 \in \intseg{0}{d-1}, \, j_1,j_2,j_3 \in \intseg{0}{s-1} \} \tag{$\ast $}
	\end{equation}
	is contained in $\psi(B'').$
	Let $i_1, i_2, i_3 \in \intseg{0}{d-1}$ and $j_1, j_2, j_3 \in \intseg{0}{s-1}$. We split the proof into four cases.
    
    Case 1: $j_1 = j_2 = j_3 = s-1$. Clearly, $[[a_{i_1,s-1}, a_{i_2,s-1}], a_{i_3,s-1}] = \id.$
    
	Case 2: $j_3 \neq s-1$. In light of  \cref{Proposition}{prop:gen_bas_commsbgp}, the elements $([a_{i_1,j_1},a_{i_2,j_2}],\id^{*(m-1)})$ and $(a_{i_3,j_3},a_{i_3,j_3}^{-1},\id^{*(m-2)}) = [a_{i_3,j_3+1}, a_{0,0}]^{-1}$ belong to $\psi(B')$, implying that 
	\[
		([[a_{i_1,j_1},a_{i_2,j_2}],a_{i_3,j_{3}}],\id^{*(m-1)}) \in \psi(B'').
	\]
   
    Now, observe from \cref{Proposition}{prop:gen_bas_commsbgp} that 
    $\psi(B'')\geq (B'')^m$. Therefore, if there exist $g =(g_0,\dots,g_{m-1}), h = (h_0,\dots,h_{m-1}) \in B$ such that $g_i \equiv_{B''} h_i$ for all $i \in \intseg{0}{m-1}$ then $g \equiv_{\psi(B'')} h$.
    
	Case 3: $j_3 = s-1,\, j_1 \neq s-1 $ and $j_2 \neq s-1$.
	Now, from the Hall--Witt identity (see \cite{Rob96}*{p.\ 123}), we can easily derive that
	\[
		[[y,x],z][[z,y],x][[x,z],y] \equiv_{B''} [[y,x],z^y][[z,y],x^z][[x,z],y^x] = \id,
	\]
	for all $x, y, z \in B$.  Setting $x=a_{i_1,j_1}, \, y =  a_{i_2,j_2}$ and $z = a_{i_3,j_3}$, we get that the element
	\begin{align*}
		([[y,x],z],\id^{*(m-1)})^{-1} \equiv_{\psi(B'')} ([[z,y],x][[x,z],y],\id^{*(m-1)})
	\end{align*}
	belongs to $\psi(B''),$ as the right-hand side product belongs to $\psi(B'')$ by Case 2.

	Case 4: $j_3 = s-1= j_1,\, j_2 \neq s-1$ or $j_3 = s-1= j_2,\, j_1 \neq s-1$. Notice that
	\[
		[[a_{i_1,j_1},a_{i_2,s-1}],a_{i_3,s-1}] \equiv_{B''} [[a_{i_2,s-1},a_{i_1,j_1}],a_{i_3,s-1}]^{-1},
	\]
	thus, it is enough to consider the first case. We claim that, for every $j \in \intseg{0}{s-1}$, it holds $[[a_{i_1,j},a_{i_2,0}],a_{i_3,j}]\equiv_{B''} \id$. Then by taking the $j_2$-th projection of the element $[[a_{i_1,s-1},a_{i_2,j_2}],a_{i_3,s-1}]$ we obtain,
	\begin{align*}
		\psi_{j_2}([[a_{i_1,s-1},a_{i_2,j_2}],a_{i_3,s-1}]) & = ([[a_{i_1,(s-1-j_2)},a_{i_2,0}],a_{i_3,(s-1-j_2)}],\id^{*(m^{j_2}-1)}) \\
		& \equiv_{\psi_{j_2}(B'')} \id, 
	\end{align*}
	implying $[[a_{i_1,s-1},a_{i_2,j_2}],a_{i_3,s-1}] \equiv_{B''} \id,$ and hence $(\ast)$ follows.
	
	If $i_2 = 0$ or $j = 0$, it is then immediate that $[[a_{i_1,j}, a_{i_2,0}],a_{i_3,j}] = \id$. Assume that $i_2 \neq 0$ and $j \neq 0$. From the presentation of $B$ given in \cref{Theorem}{thm:bas_l_pres}, we have 
	\[
		[[a_{i_1,j},\alpha(v,k)],a_{i_3,j}] = [a_{i_1,j}^{-1}a_{i_1,j}^{\alpha(v,k)},a_{i_3,j}] = [a_{i_1,j}^{-1},a_{i_3,j}]^{a_{i_1,j}^{\alpha(v,k)}} [a_{i_1,j}^{\alpha(v,k)},a_{i_3,j}] = \id,
	\]
 where $\alpha(v,k)$ is given by \eqref{eq:gen_bas_alpha(v,k)}. Now, by setting $v = (0^{\ast(i_2-1)}, 1, 0^{\ast(m-i_2-1)})$ and $k = 1$, we get $\alpha(v,k) = a_{0,0}a_{i_2,0}$ and consequently
	\begin{align*}
		\id &= [[a_{i_1,j},a_{0,0}a_{i_2,0}],a_{i_3,j}] = [[a_{i_1,j},a_{i_2,0}][a_{i_1,j}, a_{0,0}]^{a_{i_2,0}},a_{i_3,j}]\\
		& \equiv_{B''} [[a_{i_1,j},a_{i_2,0}],a_{i_3,j}][[a_{i_1,j}, a_{0,0}]^{a_{i_2,0}},a_{i_3,j}] \equiv_{B''} [[a_{i_1,j}, a_{i_2,0}],a_{i_3,j}].
	\end{align*}
	
	Next we prove (i). Assume that $s = 2$ and notice that it suffices to prove that $B'/ \psi^{-1}(\gamma_3(B)^m )$ is abelian. We use the fact that the commutator subgroup can be described by $B' = \langle [a_{i_1,1},a_{i_2,0}] \mid i_1,i_2 \in \intseg{0}{d-1} \rangle^B$ as $s=2$.
	
	Looking at the section decomposition of these generators,
	\begin{align*}
		[a_{i_1,1},a_{i_2,0}] &= (\ [a_{i_1,0},a_{i_2-1,1}], \id^{\ast(m-1)}) \quad \text{for } i_2 \neq 0, \text{ and}\\
		[a_{i_1,1},a_{0,0}] &= (a_{i_1,0}^{-1}, a_{i_1,0}, \id^{\ast(m-2)}),
	\end{align*}
	we immediately see that they commute modulo $\gamma_3(B)^m$. Thus, $B'/ \psi^{-1}(\gamma_3(B)^m )$ is abelian.
    
	(ii) The inclusion $\psi^{-1}(\gamma_3(B)^m) \leq B''$ has been already proven above. We prove that $\psi^{-1}(\gamma_3(B)^m)$ is a proper subgroup of $B''$, by showing that $B'/\psi^{-1}(\gamma_3(B)^m)$ is non-abelian. Suppose to the contrary $B'/\psi^{-1}(\gamma_3(B)^m)$ is abelian. Then, for every $i \in \intseg{0}{d-1}$ and $j \in [2,s-1]$
	\begin{align*}
	    \id \equiv_{\psi^{-1}(\gamma_3(B)^m)} [[a_{i,j},a_{0,0}], [a_{0,1},a_{0,0}]] = ([a_{i,j-1}^{-1},a_{0,0}^{-1}], [a_{i,j-1},a_{0,0}],\id^{*(m-2)}).
	\end{align*}
	This implies $[a_{i,j-1},a_{0,0}] \equiv _{\gamma_3(B)} \id$, which is a contradiction to \cref[(ii)]{Theorem}{thm:gen_bas_gammasbgp}.
\end{proof}

\section{Congruence properties of The Generalised Basilica groups}\label{sec:gen_bas_mCSP}

\noindent Here we prove that the generalised Basilica group $\bp_s(\mathcal{O}_p^{d})$ has the $p$-CSP for $d, \, s \in \Nplus$ with $s > 2$ and $p$ a prime. We follow the strategy from \cite{GUA19}, where it is proved that the original Basilica group $\mathcal{B} = \bp_2(\mathcal{O}_2)$ has the $2$-congruence subgroup property. However, on account of \cref{Theorem}{thm:gen_bas_gammasbgp} and \cref{Lemma}{lem:gen_bas_secondcommsbgp}, our reasoning must be different, and we will use \cref{Theorem}{thm:gen_bas_stab}. 

Let $G$ be a subgroup of the automorphism group of the $p$-regular rooted tree $T$ and let $\mathcal{C}$ be the class of all finite $p$-groups.

\begin{defn} [\cite{GUA19}*{Definition 5}] \label{defn GaUA} 
	A subgroup $G$ of $\Aut(T)$ has the $p$-congruence subgroup property ($p$-CSP) if every normal subgroup $N \trianglelefteq G$ satisfying $G/N \in \mathcal{C}$ contains some layer stabiliser in $G$. The group $G$ has the $p$-CSP modulo a normal subgroup $M \trianglelefteq G$ if every normal subgroup $N\trianglelefteq G$ satisfying $G/N \in \mathcal{C}$ and $M \leq N$ contains some layer stabiliser in $G$.
\end{defn}

By setting $\mathcal{C}$ as the class of all finite $p$-groups in \cite[Lemma 6]{GUA19}, we obtain the following result:
 
\begin{lemma} \label{lemma GaUA} 
	Let $G$ be a subgroup of $\Aut(T)$ and $N \trianglelefteq M \trianglelefteq G$. If $G$ has the $p$-CSP modulo $M$ and $M$ has the $p$-CSP modulo $N$ then $G$ has the $p$-CSP modulo $N$.
\end{lemma}

Let $d, \, s  \in \Nplus$ with $s>2$ and let $p$ be a prime. Set $B = \bp_s(\mathcal{O}_p^{d})$. From \cref[(vi)]{Theorem}{thm:gen_bas_properties} $B$ is weakly regular branch over its commutator subgroup $B'$ and from \cref{Lemma}{lem:gen_bas_secondcommsbgp} 
\[
	B' \geq \gamma_3(B) \geq B'' > \psi^{-1}(\gamma_3(B)^p).
\]  
We will prove that
 \begin{enumerate}
     \item $B$ has the $p$-CSP modulo $\gamma_3(B)$, and,
     \item $\gamma_3(G)$ has the $p$-CSP modulo $\psi^{-1}(\gamma_3(B)^p)$.
 \end{enumerate}
Then \cref{Theorem}{thm:gen_bas_mCSP} follows by a direct application of \cite[Theorem 1]{GUA19}. Applying \cref{Lemma}{lemma GaUA} to \cref{Proposition}{prop CSP 1} and \cref{Proposition}{prop CSP 2} we will obtain step (1). Similarly, by applying \cref{Lemma}{lemma GaUA} to \cref{Proposition}{prop CSP 6} and \cref{Proposition}{prop CSP 7} yields step (2). Now, set $M \defeq \psi^{-1}((B')^p)$ and $N \defeq \psi^{-1}(\gamma_3(B)^p)$. Considering \cref{Proposition}{prop:gen_bas_commsbgp}, \cref{Theorem}{thm:gen_bas_gammasbgp} and \cref{Lemma}{lem:gen_bas_secondcommsbgp}, we summarise the proof of \cref{Theorem}{thm:gen_bas_mCSP} in \cref{Figure}{fig:proofMCSP}.

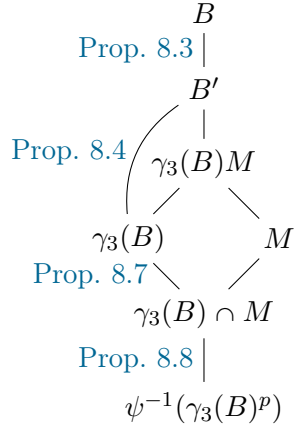
\begin{figure}
	\centering
	\begin{tikzpicture}[scale=1]
		\node (B)   at (0,3) {$B$};
		\node (der)   at (0,2) {$B'$};
		\node (gamM)   at (0,1) {$\gamma_3(B)M$};
		\node (gam) at (-1,0) {$\gamma_3(B)$};
		\node (M) at (1,0) {$M$};
		\node (gamIntM)   at (0,-1) {$\gamma_3(B)\cap M$};
		\node (N) at (0,-2.25) {$\psi^{-1}(\gamma_3(B)^p)$};
		
		\path
			(B)			edge	node	[left]	{\cref{Prop.}{prop CSP 1}}	(der)
			(der)		edge	[bend right]	node	[left]	{\cref{Prop.}{prop CSP 2}}	(gam)
			(der)		edge	(gamM)
			(gamM)		edge	(M)
			(M)			edge	(gamIntM)
			(gamIntM)	edge	node	[left]	{\cref{Prop.}{prop CSP 7}}	(N)
			(gamM)		edge	(gam)
			(gam)		edge	node	[left]	{\cref{Prop.}{prop CSP 6}\;}	(gamIntM);
	\end{tikzpicture}
	\caption{The steps of the proof of \cref{Theorem}{thm:gen_bas_mCSP}, where $M:= \psi^{-1}((B')^p)$}\label{fig:proofMCSP}
\end{figure}

\begin{prop} \label{prop CSP 1}
	The group $B$ has the $p$-CSP modulo $B'$.
\end{prop}

\begin{proof}
	Set $b_{is+j} = a_{i,j}$ for all $i \in \intseg{0}{d-1}$ and $j \in \intseg{0}{s-1}$. Define, for $r \in \intseg{0}{ds-1}$, $A_r = \langle b_r,\dots,b_{ds-1}\rangle B'$ and set $A_{ds} = B'$. We will prove that $A_r$ has the $p$-CSP modulo $A_{r+1}$ for all $r \in  \intseg{0}{ds-1}$. Then the result follows from the \cref{Lemma}{lemma GaUA}.
	
	Clearly, $A_r / A_{r+1}\St_{A_r}(n) \in \mathcal{C}$ and by \cref[(iv)]{Theorem}{thm:gen_bas_properties} we have $A_r/A_{r+1} = \langle b_r \rangle \cong \Z$. In $\Z$, the subgroups of index a power of $p$ are totally ordered, whence it suffices to prove that $\vert A_r: A_{r+1}\St_{A_{r}}(n) \vert$ tends to infinity when $n$ tends to infinity. In fact, we prove that $b_r^{p^n} \notin A_{r+1}\St_{A_r}(n d s+ r+1)$ for $n \in \mathbb{N}_0$. Assume to the contrary that $b_r^{p^n} \in A_{r+1}\St_{A_r}(n d s +r+1)$. In particular, $b_r^{p^n} \in A_{r+1}\St_B(n d s + r + 1)$. Thanks to  \cref{Theorem}{thm:gen_bas_stab}, we have
	$
 	   \St_{B}(n d s +r+1) = \langle b_0^{p^{n+1}}, \dots, b_{r}^{p^{n+1}}, b_{r+1}^{p^n}, \dots, b_{ds-1}^{p^{n}}\rangle^{B}.
	$
	Thus, there exists $x_0,\dots,x_{ds-1} \in \Z$ such that
	\begin{equation*}
		b_r^{p^n}  \equiv_{B'} b_0^{x_0p^{n+1}} \cdots b_r^{x_rp^{n+1}}b_{r+1}^{x_{r+1}} \cdots b_{ds-1}^{x_{ds-1}},
	\end{equation*}
	contradicting \cref[(iv)]{Theorem}{thm:gen_bas_properties}. 
\end{proof}

\begin{prop} \label{prop CSP 2}
	 The group $B'$ has the $p$-CSP modulo $\gamma_{3}(B)$.
\end{prop}

\begin{proof}
	Notice from \cref[(ii)]{Theorem}{thm:gen_bas_gammasbgp} that $\gamma_3(B)$ is a subgroup of index a power of $p$ in $B'$ and hence it suffices to prove that $\St_{B'}(n)$ is contained in $\gamma_3(B)$ for some $n$, equivalently 
	$\vert B'/ \gamma_3(G)\St_{B'}(n)\vert = \vert B'/ \gamma_3(B)\vert$.
	Observe that,
	\begin{align*}
		B'/ \gamma_3(B)\St_{B'}(n)
		\cong 		
		B'\St_{B}(n)/\gamma_3(B)\St_{B}(n).
	\end{align*}
	Now, in light of \cref[(ii)]{Theorem}{thm:gen_bas_gammasbgp}, we choose $n \in \Nplus$ such that the set
	\begin{equation*}
	  \{[a_{i,j},a_{0,0}] \mid i \in  \intseg{0}{d-1}, \, j \in \intseg{1}{s-1}\} \cup 
	  \{ [a_{0,1},a_{i',0}] \mid i' \in \intseg{1}{d-1}\} \cup \{[a_{0,1},a_{0,0}]^p\},
	\end{equation*}
	has trivial intersection with $\St_B(n)$. 
	One can easily compute from the description of the stabilisers in \cref{Theorem}{thm:gen_bas_stab} that $n = ds +2$ is the smallest number with this property.
	We construct a group $H$ which admits an epimorphism from the group $B/\gamma_3(B)\St_{B}(ds +2)$ and see that commutator subgroup $H'$ has the desired size.
		
	Now fix $n = ds + 2$ and set $\Gamma = B/\gamma_3(B)\St_B(n)$. Again from \cref{Theorem}{thm:gen_bas_stab} we have $\St_B(n) = \langle b_0^{p^2}, b_1^{p^2}, b_2^{p},  \dots, b_{ds-1}^{p} \rangle ^B$, where $b_{is+j} = a_{i,j}$ as in the proof of \cref{Proposition}{prop CSP 1}. By a straightforward calculation using the presentation of $B / \gamma_3(B)$, given in the proof of \cref[(ii)]{Theorem}{thm:gen_bas_gammasbgp}, we obtain the following presentation for $\Gamma$:
	\begin{equation}\label{pres of gamma}
	    \LD\ \mathcal{S} \mid \mathcal{R}\ \RD,
	\end{equation}
	where $\mathcal{S} = \{b_r \mid r \in \intseg{0}{ds-1}\}$ and
	\begin{align*}
		\mathcal{R} = \left\langle 
		\begin{array}{l|l}
		    b_{0}^{p^2},b_{1}^{p^2},b_{t}^p, {[}b_{t},b_{t'}{]},
		    & t,t' \in \intseg{2}{ds-1} \\
			{[}b_{1},b_{t''}{]}, 
			&t'' \in \intseg{2}{ds-1}, \text{ not a multiple of } s\\
			{[}b_{0},b_{is}{]},  \gamma_3(F)
			& i \in \intseg{1}{d-1}
		\end{array}\right\rangle,
	\end{align*}
	where $F$ is the free group on the set of generators of $\Gamma$.
 
	Let $R$ be the ring $\Z/ p^2 \Z$. Let $\mathrm{UT}_{ds+1}(R) \leq \mathrm{GL}_{ds+1}(R)$ be the group of all upper triangular matrices over $R$ with entries $1$ along the diagonal. Denote by $E_{i,j}(\ell)$ the element of $\mathrm{UT}_{ds+1}(R)$ with the entry $\ell \in R$ at the position $(i,j)$. For $i \in \intseg{1}{d(s-1)-1}$ and $j \in \intseg{1}{d-1}$, define 
	\begin{align*}
		x_i &= E_{i,ds-1}(p),
		& y_j & = E_{d(s-1)+j,ds}(p),\\
		  y & = E_{ds-1,ds}(1),
		& z & = E_{ds, ds+1}(1),
	\end{align*} 
	and define $\mathcal{H}$ to be the subgroup of $\mathrm{UT}_{ds+1}(R)$ generated by the set $\{x_i,y_j,y,z\}$. By abuse of notation denote the image of the set of generators of $\mathcal{H}$ in the quotient group $\mathcal{H}/\gamma_3(\mathcal{H})$ by the same symbols and set $H= \mathcal{H}/\gamma_3(\mathcal{H})$. By an easy computation, we obtain
	\begin{equation*}
		x_i^p  = y_j^p = y^{p^2} = z^{p^2} = [x_i,x_{i'}] = [y_j,y_{j'}] = [y,y_j] = [x_i,y_j] = [x_i,z]  =  \id,
	\end{equation*}
	for all $i,i' \in \intseg{1}{d(s-1)-1}$ and $j,j' \in \intseg{1}{d-1}$. 
	Now, fix a bijection $\alpha$ from the set $\{b_{r} \mid r \in \intseg{2}{ds-1} \setminus \{s,2s,\dots,(d-1)s\}$ to the set $\{x_i \mid i \in \intseg{1}{d(s-1)-1}\}$. Define a map $\varphi$ from the set of generators of $\Gamma$ to the set of generators of $H$ by
	\[
		\begin{array}{lllll}
		\varphi(b_{0}) &= y && \varphi(b_{1}) &= z\\
		\varphi(b_{js}) &= y_j \text{ for } j \in \intseg{1}{d-1} && \varphi(b_{r}) &= \alpha(b_{r}), \text{ otherwise}.
		\end{array}
	\]
	Then $\varphi$ extends to an epimorphism $\Gamma \to H$, since as seen above, $\varphi(b_{r})$ satisfies all the relations of the given presentation \eqref{pres of gamma} of the group $\Gamma$. Furthermore, observe that the commutator subgroup of $H$ is generated by the union of the sets 
	 \[
	\{\ [x_i,y] \mid i \in \intseg{1}{d(s-1)-1}\ \} \cup  \{\ [y_j,z]  \mid j \in \intseg{1}{d-1}\ \} \cup \{\ [y,z]\ \}.
	 \]
	Hence,
	\begin{align*}
		\vert \Gamma' \vert \geq \vert \varphi(\Gamma') \vert = \vert H' \vert = p^{d(s-1)-1}p^{d-1}p^2=p^{ds}.
	\end{align*}
	Indeed $\vert \Gamma' \vert \leq \vert B'/\gamma_3(B) \vert = p^{ds}$, and thus $\vert \Gamma' \vert = p^{ds}$, completing the proof.    
\end{proof}

We now need two general lemmata.

\begin{lemma}\label{lemma CSP 4}
	Let $H \leq \Aut(T)$ and $L, K \trianglelefteq H$ with $L \leq K$ and let $\mathcal{C}$ be the class of all finite $p$-groups. Assume further that $H/K \in \mathcal{C}$ and $H/L$ is abelian. If $H$ has the $p$-CSP modulo $L$, then $K$ has the $p$-CSP property modulo $L$.
\end{lemma}
 
\begin{proof}
	Let $\Tilde{K}$ be a normal subgroup of $K$ satisfying $L \leq \Tilde{K}$ and $K/\Tilde{K} \in \mathcal{C}.$ Since $H/L$ is abelian, $\Tilde{K}/L$ is normal in $H/L$ and hence $\tilde{K}$ is normal in $H$. Also notice that $H/\Tilde{K} \in \mathcal{C}$. As $H$ has the $p$-CSP there exists $n \in \N$ such that $\St_H(n) \leq \Tilde{K}$. In particular $\St_K(n) = \St_H(n) \cap K \leq \St_H(n) \leq \Tilde{K}$, completing the proof.
\end{proof}

\begin{lemma}\label{lemma CSP 5}
	Let $H \leq \Aut(T)$ and $L, K \trianglelefteq H$. If $KL$ has the $p$-CSP modulo $L$, then $K$ has the $p$-CSP property modulo $K \cap L$.
\end{lemma}

\begin{proof}
	Choose $\Tilde{K} \trianglelefteq K$ with $K \cap L \leq \Tilde{K}$ and $K/\Tilde{K} \in \mathcal{C}$. Then, $\Tilde{K}L \trianglelefteq KL$ and 
	$
		KL/\Tilde{K}L \cong K/\Tilde{K} \in \mathcal{C}.
	$
	As $KL$ has the $p$-CSP property modulo $L$, it holds that $\St_{KL}(n)\leq \Tilde{K}L$ for some $n$. Thus, $\St_K(n)=\St_{KL}(n)\cap K \leq \Tilde{K}L \cap K = \Tilde{K}$.
\end{proof}

\begin{prop} \label{prop CSP 6}
	The group $\gamma_3(B)$ has the $p$-CSP modulo $\gamma_3(B) \cap M$.
\end{prop}

\begin{proof}
	We prove that $\gamma_3(B)M$ has the $p$-CSP modulo $M$. Then by \cref{Lemma}{lemma CSP 5} we obtain the result. It follows from \cref{Proposition}{prop:gen_bas_commsbgp} and \cref[(ii)]{Theorem}{thm:gen_bas_gammasbgp} that $B'/M$ is abelian and that $B'/\gamma_3(B)M \in \mathcal{C}$, respectively. Thanks to \cref{Lemma}{lemma CSP 4}, it is enough to prove that $B'$ has the $p$-CSP modulo $M$. 
 
	Let $ i \in \intseg{0}{d-1}, j \in \intseg{1}{s-1}$ and $k \in \intseg{1}{p-1}$. Define $c_{i(s-1)+j} \defeq b_{is+j} \defeq a_{i,j}$. Set $t = i(s-1)+j$ and $r = is+j$ and note that $c_t$ is a relabeling of the elements $b_r$ (defined in the proof of \cref{Proposition}{prop CSP 1}) by excluding the elements of the form $b_{is}$ for $i \in \intseg{0}{d-1}$.
	From \cref{Proposition}{prop:gen_bas_commsbgp}, we have
	\[
		B'/M = \langle\ [a_{i,j}, a_{0,0}^k] \mid  i \in \intseg{0}{d-1}, \, j \in \intseg{1}{s-1}, \, k \in \intseg{1}{p-1}\ \rangle.
	\]
	Set $\ell = (k-1)(ds-d)+t $ and $ e_{\ell} = [c_{t}, a_{0,0}^{k}]$. Then,
	\begin{align*}
		\psi(e_{\ell})= \psi([c_{t}, a_{0,0}^{k}]) = \psi([b_{r}, a_{0,0}^{k}]) =
             (b_{r-1}^{-1},\id^{*(k-1)},b_{r-1},\id^{*(p-k-1)}).
	\end{align*}  
	For $\ell \in \intseg{1}{(p-1)(ds-d)}$, set $M_{\ell}=\langle e_{\ell} ,\dots,e_{(p-1)(ds-d)}\rangle M$ and $M_{(p-1)(ds-d)+1} = M$. It follows from \cref[(iv)]{Theorem}{thm:gen_bas_properties} that $M_{\ell}/M_{\ell+1} = \langle e_{\ell}\rangle \cong \Z$. We will prove that $\vert M_{\ell} : M_{\ell+1} \St_{M_\ell}(n) \vert $ tends to infinity as $n$ tends to infinity. 
	Assume to the contrary that there are $n, n' \in \Nplus$ such that for all $\tilde{n} \geq n'$, $e_{\ell}^{p^n} \in M_{\ell+1}\St_{M_{\ell}}(\tilde{n})$. There exist $x_{\ell+1},\dots,x_{(p-1)(ds-d)}\in \Z$ such that
	\[
	e_{\ell}^{p^n}e_{\ell+1}^{x_1} \cdots e_{(p-1)(ds-d)}^{x_{(p-1)(ds-d)}} \in M \St_{{M_{\ell}}}(\tilde{n}) \leq M \St_B(\tilde{n}),
	\]
	hence
	\begin{align*}
		\psi(e_{\ell}^{p^n}e_{\ell+1}^{x_1} \cdots e_{(p-1)(ds-d)}^{x_{(p-1)(ds-d)}}) \in (B')^p\cdot (\St_{B}(\tilde{n}-1))^p.
	\end{align*}
	Consider the $k$-th coordinate, $xb_{r-1}^{p^n} \in B' \St_{B}(\tilde{n}-1)$, where $x$ is a product of elements of the form $b_{r'}$ such that $r' > r-1$. Then $x \in A_r$, where $A_r$ is defined as in the proof of \cref{Proposition}{prop CSP 1}. This implies 
	\(
		b_{r-1}^{p^n} \in A_{r} \St_{B}(\tilde{n}) \text{ for all } \tilde{n} \geq n'-1,
	\)
	which contradicts \cref{Proposition}{prop CSP 1}.
\end{proof}

\begin{prop} \label{prop CSP 7}
	The group $\gamma_3(B)\cap M$ has the $p$-CSP modulo $N$.
\end{prop}

\begin{proof}
	It is straightforward from \cref[(ii)]{Theorem}{thm:gen_bas_gammasbgp} that the group $M/N$ is a finite abelian and $M/N \in \mathcal{C}$. By \cref{Lemma}{lemma CSP 4}, it suffices to prove that $M$ has the $p$-CSP modulo $N$. From \cref{Proposition}{prop CSP 2}, it follows that  $\text{St}_{B'}(n) \leq \gamma_3(B)$ for some $n$. Therefore, 
	\begin{align*}
		\psi(\St_{M}(n+1)) \leq  (\St_{B'}(n))^p \leq \gamma_3(B)^p,
 \end{align*}
  and hence 
  $
    \St_{M}(n+1) \leq  \psi^{-1}((\St_{B'}(n))^p)\leq N.
  $
\end{proof}

\begin{proof}[Proof of \cref{Theorem}{thm:gen_bas_mCSP}]
	By applying \cref{Lemma}{lemma GaUA} to \cref{Proposition}{prop CSP 1} and \cref{Proposition}{prop CSP 2} we obtain that
the group $B$ has the $p$-CSP modulo $\gamma_{3}(B)$.
 Further application of \cref{Lemma}{lemma GaUA} to \cref{Proposition}{prop CSP 6} and \cref{Proposition}{prop CSP 7} yields that $\gamma_{3}(G)$ has the $p$-CSP modulo $N$.
Now, the result follows by \cite[Theorem 1]{GUA19}.
\end{proof}

\end{document}